\documentclass[11pt]{article}
\usepackage{geometry, booktabs, subcaption}
\usepackage{titlesec}
\usepackage[T1]{fontenc}
\usepackage{authblk}
\usepackage{amssymb,amsmath}
\usepackage{graphics}
\usepackage[square]{natbib}
\RequirePackage[colorlinks,citecolor=blue,urlcolor=blue,linkcolor=black]{hyperref}
\usepackage{mathrsfs,amsthm,amsmath}
\usepackage{amssymb,multirow}
\usepackage{amsfonts}
\usepackage{stmaryrd}
\usepackage{dsfont}

\usepackage{ucs}
\usepackage{amsmath}
\usepackage{amsfonts}
\usepackage{amssymb}
\usepackage{amsthm}
\usepackage{fontenc}
\usepackage{graphicx}
\usepackage{cleveref}
\usepackage{bm}

\usepackage{bbm}

\usepackage{mathrsfs}
\usepackage{amssymb,multirow}
\usepackage{amsfonts}
\usepackage{stmaryrd}
\usepackage{dsfont}
\usepackage{subcaption}
\usepackage{placeins}

\usepackage{ucs}
\usepackage{amsmath}
\usepackage{amsfonts}
\usepackage{amssymb}
\usepackage{fontenc}
\usepackage{graphicx}
\usepackage[export]{adjustbox}
\usepackage{cleveref}

\RequirePackage[OT1]{fontenc}

\newcommand{\mutilde}{\tilde \mu}

\newcommand{\prob}{\mathrm{Pr}}

\newcommand{\dd}{\mathrm{d}}
\newcommand{\indicator}{\mathrm{I}}

\renewcommand{\mid}{\,|\,}

\newcommand{\X}{\mathbb{X}}

\newcommand{\Sd}{\mathbb{S}}

\newcommand{\R}{\mathbb{R}}
\newcommand{\E}{\mathds{E}}

\newcommand{\N}{\mathbb{N}}

\newcommand{\calC}{\mathscr{C}}

\def\simind{\overset{\mathrm{ind}}{\sim}}
\def\simiid{\overset{\mathrm{iid}}{\sim}}

\newcommand{\ddr}{\mathrm{d}}

\newtheorem{thm}{Theorem}

\newtheorem{defi}[thm]{Definition}

\newtheorem{prp}[thm]{Proposition}

\newtheorem{corollary}[thm]{Corollary}

\usepackage{fancyhdr}

\title{\Huge Random measure priors in Bayesian recovery from sketches}
\author[1]{Mario Beraha}
\author[2]{Stefano Favaro}
\author[3]{Matteo Sesia}
\affil[1]{\normalsize{Department of Mathematics, Polytechnic University of Milan, Milan, Italy}} 
\affil[2]{\normalsize{Department of Economics and Statistics, University of Torino and \newline Collegio Carlo Alberto, Torino, Italy}}
\affil[3]{\normalsize{University of Southern California, Los Angeles, California, United States}}

\begin{document}

\maketitle

\begin{abstract}
This paper introduces a Bayesian nonparametric approach to frequency recovery from lossy-compressed discrete data, leveraging all information contained in a sketch obtained through random hashing. By modeling the data points as random samples from an unknown discrete distribution endowed with a Poisson-Kingman prior, we derive the posterior distribution of a symbol's empirical frequency given the sketch. This leads to principled frequency estimates through mean functionals, e.g., the posterior mean, median and mode. We highlight applications of this general result to Dirichlet process and Pitman-Yor process priors. Notably, we prove that the former prior uniquely satisfies a sufficiency property that simplifies the posterior distribution, while the latter enables a convenient large-sample asymptotic approximation. Additionally, we extend our approach to the problem of cardinality recovery, estimating the number of distinct symbols in the sketched dataset. Our approach to frequency recovery also adapts to a more general ``traits'' setting, where each data point has integer levels of association with multiple symbols, typically referred to as ``traits''.  By employing a generalized Indian buffet process, we compute the posterior distribution of a trait's frequency using both the Poisson and Bernoulli distributions for the trait association levels, respectively yielding exact and approximate posterior frequency distributions.
\end{abstract}



\section{Introduction}

\subsection{Background and motivation}

Information recovery about a large discrete dataset given a (lossy) compressed representation, or {\em sketch}, of that data is a classical problem at the crossroad of computer science and information theory \citep{Mis(82),Alo(99),Man(02),Kar(03),Cha(04),Cor(05),Ind(06)}. 
 Sketching is often driven by the need to manage memory constraints, as handling vast numbers of symbols can be computationally intensive, and by privacy concerns, especially when dealing with sensitive data \citep{Blu(20), Cor(20),Med(22)}.
Two well-known challenges based on sketched data are {\em frequency recovery} and {\em cardinality recovery}. For concreteness, we begin by discussing frequency recovery.

In the typical ``species'' setting, where the original dataset comprises $n \geq 1$ points $(x_1, \ldots, x_n)$ with each $x_i$ corresponding to a symbol or ``species'' label taking values in a set $\mathbb{S}$, the goal of {\em frequency recovery} is to estimate the number of occurrences of a new object $x_{n+1}$ in $(x_{1},\ldots,x_{n})$, denoted as $f_{x_{n+1}}$. Formally, this can be written as:
\begin{displaymath}
f_{x_{n+1}}=\sum_{i=1}^{n}I(x_{i}=x_{n+1}),
\end{displaymath}
with $I(\cdot)$ denoting the indicator function.
This problem  is relevant for many applications, including machine learning in high-dimensional feature spaces \citep{Shi(09),Agg(10)}, cybersecurity in tracking password popularity \citep{Sch(10)}, web and social network data analysis \citep{Son(09),Cor(17)}, natural language processing \citep{Goy(12)}, sequencing analysis in biological sciences \citep{Zha(14),Ber(16),Sol(16),Mar(19),Leo(20)}, and privacy-protecting data analysis \citep{Dwo(10),Mel(16),Cor(17),Cor(18),Koc(20)}. 

The count-min sketch (CMS) is a popular algorithm for frequency recovery \citep{Cor(05)}. It relies on a sketch of $(x_{1},\ldots,x_{n})$ obtained by $D\geq1$ independent $J$-wide random hash functions $h_{k}:\mathbb{S}\rightarrow [J] := \{1,\ldots,J\}$, for $k \in [D] := \{1,\ldots, D\}$ and $J\geq1$. Each hash function maps the $x_{i}$'s into $J$ buckets, defining the sketch  $\mathbf{C}_{D,J}\in\mathbb{N}_{0}^{D\times J}$, with $\mathbb{N}_{0}=\mathbb{N}\cup\{0\}$, whose $(k,j)$-th element $C_{k,j}$ counts the data points mapped by the $k$-th hash function into the $j$-th bucket. Based on $\mathbf{C}_{D,J}$, the CMS bounds $f_{x_{n+1}}$ from above by taking the smallest count among the $D$ buckets into which $x_{n+1}$ is mapped, i.e.,
\begin{equation} \label{cms}
\hat{f}_{x_{n+1}}=\min\{C_{1,h_{1}(x_{n+1})},\ldots,C_{D,h_{D}(x_{n+1})}\}.
\end{equation}
We refer to Figure~\ref{fig:diagram-cms} in Section~\ref{sec2} for a schematic visualization of this procedure, focusing on the special case where the data are sketched using a single hash function. 
In general, while the CMS upper bound is remarkable for its simplicity and robustness, it becomes loose if hash collisions (different objects mapped into the same bucket) are frequent, due to the pessimistic assumption that data are fixed \citep[Chapter 3]{Cor(20)}.

This challenge has recently motivated statistical approaches that treat the sketched data as random and rely on modeling assumptions to obtain more informative estimates. 
The first Bayesian nonparametric (BNP) approach to frequency recovery was introduced by \citet{Cai(18)}. They modeled the $x_{i}$'s as a random sample $(X_{1},\ldots,X_{n})$ from an unknown discrete distribution endowed with a Dirichlet process (DP) prior \citep{Fer(73)} and obtained estimates of $f_{X_{n+1}}$ as mean functionals of the posterior distribution of $f_{X_{n+1}}$ given $(C_{1,h_{1}(X_{n+1})},\ldots,C_{D,h_{D}(X_{n+1})})$, e.g., the mean, median and mode. In addition to enabling the inclusion of prior knowledge on the data distribution, the BNP approach allows assessing uncertainty using posterior distributions. See  \citet{Dol(21),Dol(23)} for an extension of \citet{Cai(18)} to the Pitman-Yor process (PYP) prior \citep{Pit(97)}.

As outlined in Section~\ref{sec:intro-preview}, this paper extends the foundations laid by previous research to establish a more versatile and comprehensive BNP framework. 
In particular, our novel framework allows: conditioning on all information contained in the sketched data, employing a wider array of prior distributions, and tackling other information retrieval challenges beyond frequency recovery in the ``species'' setting.

Specifically, we focus on analyzing a sketch $\mathbf{C}_{J}=(C_{1},\ldots,C_{J})\in\mathbb{N}_{0}^{J}$ obtained by a {\em single hash function}. 
The work of \citet{Cai(18)} was motivated by the goal developing a learning-augmented version of the CMS, hence their interest in a posterior distribution with respect to the same information from a sketch $\mathbf{C}_{D,J}$ obtained from $D$ distinct hash functions as in \eqref{cms}. 
By contrast, we develop a purely Bayesian approach, separate from the CMS, including in the posterior distribution all information from the sketch $\mathbf{C}_{J}$. 
From a Bayesian statistical perspective, our approach is arguably more natural than that of \citet{Cai(18)} because: i) the practical usefulness of combining a statistical model for the $x_{i}$'s with the sketch $\mathbf{C}_{D,J}$ from multiple hash functions would be unclear, being such a sketch designed for a (model-free) recovery algorithm; ii) the use of a posterior distribution with respect to the sole $C_{k,h_{k}(X_{n+1})}$'s may determine a loss of information, unless we can verify that the $C_{k,h_{k}(X_{n+1})}$'s are sufficient to estimate $f_{X_{n+1}}$.

\subsection{Preview of our contributions} \label{sec:intro-preview}

\subsubsection{BNP frequency recovery}

Our first contribution is a novel BNP approach to frequency recovery that utilize the {\em full} posterior distribution of $f_{X_{n+1}}$ given a sketch obtained with a single random hash function. This approach departs from that of \citet{Cai(18)}, which did not condition on the information contained in the sketch outside the bucket into which $X_{n+1}$ is hashed.
We compute the full posterior distribution of $f_{X_{n+1}}$ given $\mathbf{C}_{J}$ and the bucket in which $X_{n+1}$ is hashed, denoted as $h(X_{n+1})$. We derive this posterior under a broad class of Poisson-Kingman (PK) priors \citep{Pit(03)}, encompassing both the DP and PYP priors. For the DP prior, we find that the posterior depends on $\mathbf{C}_{J}$ only through $C_{h(X_{n+1})}$; this is consistent with the approach of \citet{Cai(18)}, proving that it relies on a sufficient statistic. However, our findings also reveal that the DP prior is the only PK prior satisfying this sufficiency property. Additionally, we show that the posterior distribution is computationally intractable under the PYP prior if $n$ is large. Thus, we develop a large-sample asymptotic approximation. 

\subsubsection{BNP cardinality recovery}

Our second contribution extends the BNP approach to frequency recovery to address the \textit{cardinality recovery} problem \citep[Chapter 2]{Cor(20)}---another classical problem in computer science. Here, the objective is to utilize the same sketch obtained with a single hash function to estimate the number of distinct symbols in $(x_{1},\ldots,x_{n})$, i.e., 
\begin{displaymath}
k_{n}=|{x_{1},\ldots,x_{n}}|,
\end{displaymath}
where $|\cdot|$ denotes the cardinality set-function. Prior works treated frequency and cardinality recovery as separate problems, each requiring a different sketching algorithm. For instance, the hyperloglog algorithm is widely used for cardinality recovery \citep{Fla(07),Fla(83),Fla(85)}. See also \citet{Cha(06)}, \citet{Che(11)}, \citet{Tin(14),Tin(16)}, and \citet{Pet(21)} for recent contributions, some of which rely on modeling assumptions for the data. We show that the BNP approach enables a direct connection between the frequency and cardinality recovery problems, yielding a posterior mean estimate of $k_{n}$ from the sketch $\mathbf{C}_{J}$. To the best of our knowledge, this is the first approach enabling both frequency and cardinality recovery from the same sketch $\mathbf{C}_{J}$. More generally, for any $l \in [n]$, we obtain a posterior mean estimate of the number $m_{l,n}$ of distinct symbols with frequency $l$ in $(x_{1},\ldots,x_{n})$, commonly referred to as the {\em $l$-cardinality}.  This enables the recovery of the partition structure of $(x_{1},\ldots,x_{n})$ from $\mathbf{C}_{J}$.  

\subsubsection{BNP frequency recovery in the ``traits" setting}

Our third contribution is a BNP approach to a new ``traits'' setting of frequency recovery, in which each data point may be associated with more than one symbol, called ``traits'', and exhibits levels of association with each of those traits \citep{Cam(18)}. Multi-trait data arise in many domains: single-cell data encompass multiple genes with their expression levels, members of social networks connect with multiple friends to whom they send messages, and documents encompass different topics with their words. We consider $n\geq1$ data points $(x_{1},\ldots,x_{n})$, where each $x_{i}$ takes on a value in a set $\mathbb{S}^{\infty}\times\mathbb{N}_{0}^{\infty}$ representing traits and their levels. We assume the $x_{i}$'s to be modeled as a random sample $(X_{1},\ldots,X_{n})$ from the generalized Indian buffet process \citep{Jam(17)}. In this setting, the $X_{i}$'s are random counting measures, i.e. $X_{i}=\sum_{k\geq1}A_{i,k}\delta_{w_{k}}$, where $A_{i,k}\in\mathbb{N}_{0}$ is the level of association of the trait $w_{k}\in\mathbb{S}$ for the $i$-th data point. The distribution of $X_{i}$ is determined by: i) a distribution $G_{A}$ for the level of association $A_{i,k}$, which depends on a parameter $J_{k}>0$, for $k\geq1$; ii) the law of a completely random measure (CRM) serving as a prior distribution for the $J_{k}$'s. Based on a sketch $\mathbf{C}_{J}\in\mathbb{N}_{0}^{J}$ of $(X_{1},\ldots,X_{n})$ generated by a random hash function, we compute the posterior distribution of the empirical frequency level of a trait for a new $X_{n+1}$, given $\mathbf{C}_{J}$ and the bucket in which such a trait is hashed. It emerges that any CRM prior leads to a posterior distribution that depends on $\mathbf{C}_{J}$ solely through $C_{h(X_{n+1})}$, exhibiting a sufficiency property akin to that found in the ``species'' setting under the DP prior. Applications to Poisson and a Bernoulli level $G_{A}$ are presented in details. Our findings illustrate how the BNP approach facilitates frequency recovery in the broader ``traits'' setting, for which we are unaware of any existing algorithms. 

\subsection{Related work}

Several recent studies have studied from a statistical perspective the frequency recovery problem under the ``species'' setting. On the frequentist front, \cite{Tin(18)} introduced a bootstrap method tailored to the CMS algorithm, focusing on asymptotic guarantees. \cite{Ses(22)} and \cite{Ses(23)} proposed an alternative frequentist approach based on conformal inference ideas \citep{Vov(05)}, obtaining uncertainty estimates with finite-sample guarantees for any (possibly non-linear) sketch. It is worth noting that these different approaches can be viewed as complementary to each other. Specifically, conformal inference can be combined with both the bootstrap \citep{Ses(22)} and our Bayesian approach. This combination can, for example, yield ``calibrated'' Bayesian credible intervals with finite-sample frequentist properties \citep{Ses(22)}. Recently, \citet{Ber(23)} proposed a ``smoothed'' BNP approach aimed at mitigating some of the computational limitations inherent in Bayesian methods, which we highlight in this paper. While the smoothed approach yields computationally simpler estimators compared to our full Bayesian approach, it is limited to the case of normalized random measures, thus excluding the Pitman-Yor process, and only addresses the ``species'' setting. In the case of the DP prior, it turns our that the ``smoothed'' and BNP estimators coincide, yet the DP stands out as the only random measure for which this alignment occurs.

\subsection{Organization of the paper}

The paper is structured as follows. Section~\ref{sec2} considers the ``species" setting for frequency recovery, developing our BNP approach under a general PK prior, and in the special cases of the DP and PYP priors.  Also within the ``species" setting, Section~\ref{sec21} presents a BNP approach to cardinality recovery under the DP and PYP priors, obtaining estimators for the $l$-cardinality and the cardinality of the dataset.  In Section~\ref{sec3} we consider the ``traits" setting for frequency recovery, developing our BNP approach under a general CRM prior and a general distribution for the level of association, as well as for the special cases of Poisson and Bernoulli level distribution. Section~\ref{sec:numerics} contains an empirical validation of our methods on synthetic and real data, whereas Section~\ref{sec4} discusses some directions for future work. Proofs and other technical derivations are deferred to the Appendices. A software implementation of our methods is available at \url{https://github.com/mberaha/BNPSketching}.


\section{BNP frequency recovery in the ``species" setting}\label{sec2}

\subsection{Problem statement and setup}

\begin{figure}[!htb]
\centering
\includegraphics[width=12cm,height=6cm,keepaspectratio]{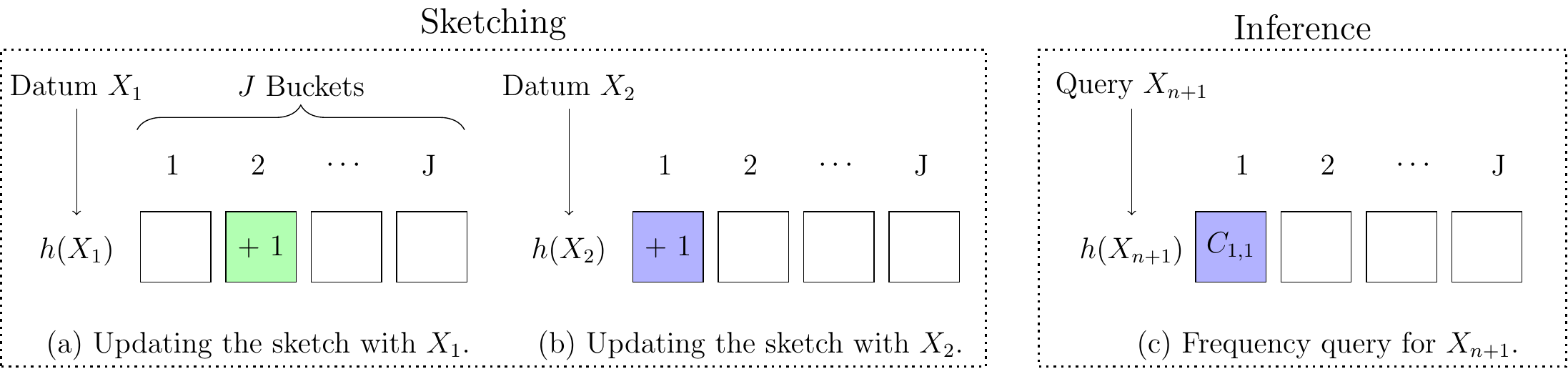}
\caption{\small{Data sketching in the ``species'' setting. Each object $X_i$, for $i \in \{1,\ldots,n\}$, is assigned to one of $J$ possible ``buckets'' (shown in different colors) by a random hash functions $h$, and the corresponding counters are incremented by one. One interesting problem is to estimate the empirical frequency of a ``query'' $X_{n+1}$ based on the sketch and the bucket assignments for $X_{n+1}$.}} \label{fig:diagram-cms}
\end{figure}

For $n\geq1$, consider data points  $(x_{1},\ldots,x_{n})$, with each $x_{i}$ representing a symbol (or ``species") label from a dictionary $\mathbb{S}$. 
Consider having access only to a sketch of $(x_{1},\ldots,x_{n})$, obtained through random hashing \citep[Chapter 5 and Chapter 15]{Mit(17)}. For an integer $J \geq 1$, let $h$ be a random hash function of width $J$, defined as a random mapping from $\mathbb{S}$ to $[J]$, chosen from a pairwise independent hash family $\mathcal{H}_{J}$. That is, $h: \mathbb{S} \to [J]$, and, for any $j_1,j_2 \in [J]$ and fixed $x_1,x_2 \in \mathbb{S}$ such that $x_1 \neq x_2$,
\begin{displaymath}
\text{Pr}[h(x_{1})=j_{1},\,h(x_{2})=j_{2}]= \frac{1}{J^{2}}.
\end{displaymath}
The pairwise independence of $\mathcal{H}_{J}$, also known as strong universality, implies uniformity, meaning that $\Pr[h(x)=j]=J^{-1}$ for $j\in[J]$. 
Hashing $(x_{1},\ldots,x_{n})$  through $h$ produces a random vector $\mathbf{C}_{J}=(C_{1},\ldots,C_{J})\in\mathbb{N}_{0}^{J}$, termed ``sketch'', whose $j$-th element (bucket) is
\begin{displaymath}
C_{j} = \sum_{i=1}^{n} I( h(x_i) = j),
\end{displaymath}
so that $\sum_{1\leq j\leq J} C_{j} = n$. This sketch generally has a smaller (physical) size than $(x_{1},\ldots,x_{n})$ due to the collisions of the $x_{i}$'s induced by random hashing \citep[Chapter 3]{Cor(20)}. The sketch $\mathbf{C}_{J}$ is a special version of the sketch $\mathbf{C}_{D,J}\in\mathbb{N}_{0}^{D\times J}$ at the basis of the CMS \citep{Cor(05)}, which simultaneously applies a collection of $D\geq1$ independent random hash functions from $\mathcal{H}_{J}$ \citep[Chapter 3]{Cor(20)}.
Based on a sketch $\mathbf{C}_{J}$ of $(x_{1},\ldots,x_{n})$, we study the BNP estimation of the empirical frequency $f_{x_{n+1}}$ of $x_{n+1}$ in $(x_{1},\ldots,x_{n})$; see Figure~\ref{fig:diagram-cms} for a schematic visualization.

\subsection{The BNP model}

Inspired by \citet{Cai(18)} and \citet{Dol(23)}, we rely on two modeling assumptions: i) the data $(x_{1},\ldots,x_{n})$ are modeled as a random sample $\mathbf{X}_{n}=(X_{1},\ldots,X_{n})$ from an unknown discrete distribution $P$ endowed by a prior $\mathscr{P}$; ii) the hash family $\mathcal{H}_{J}$ is independent of $P$. Then, we write the BNP model as:
\begin{align}\label{eq:exchangeable_model_hash}
\begin{split}
  C_{j} &\,=\, \sum_{i=1}^{n} I( h(X_i) = j), \qquad j\in[J], \\
  h&\,\sim\,\mathcal{H}_{J}, \\
  X_1,\ldots,X_{n}\,|\,P &\,\simiid\, P, \\
  P&\,\sim\,\mathscr{P}.
\end{split}
\end{align}
Under this model, the problem of recovering the empirical frequency $f_{X_{n+1}}$ from the sketch consists of computing the posterior distribution of $f_{X_{n+1}}$ given $\mathbf{C}_{J}$ and the bucket in which $X_{n+1}$ is hashed, i.e., $h(X_{n+1})$. Below, we compute this posterior assuming $\mathscr{P}$ belongs to the broad class of PK priors, and then we specialize our result to DP and PYP priors.

\subsection{Background on PK priors}

To define a PK prior \citep{Pit(03)}, we consider a completely random measure (CRM) $\mutilde$ on $\mathbb{S}$, which is a random element with values on the space of bounded measures on $(\mathbb{S},\mathcal{S})$, such that, for $k\geq1$ and for disjoint Borel sets $A_1, \ldots, A_k \in \mathcal{S}$, the random variables $\mutilde(A_{1}),\ldots,\mutilde(A_{k})$ are independent \citep{Kin(67)}. We focus on CRMs of the form
$\mutilde(\cdot) = \int_{\R_+} s \tilde N(\dd s,\, \cdot) = \sum_{k \geq 1} J_k \delta_{W_k}(\cdot)$,
where $\tilde N = \sum_{k \geq 1} \delta_{(J_k, W_k)}$ is a Poisson random measure on $\R_+ \times \mathbb{S}$ with L\'evy intensity $\nu(\dd s, \dd x)$, which characterizes the distribution of $\tilde{\mu}$ in terms its random jumps $J_{k}$'s and random locations $W_{k}$'s \citep{Kin(67),Kin(93)}. We focus on homogeneous L\'evy intensities, in the form $\nu(\dd s, \, \dd x) = \theta \rho(s) \dd s \, G_0(\dd x)$, where $\theta > 0$ is a parameter, $G_0$ is a nonatomic probability measure on $\mathbb{S}$, and $\rho(s) \dd s$ is a measure on $\R_+$ such that $\int_{\R_+} \rho(s) \dd s = +\infty$ and
\begin{equation}\label{eq_cond}
\psi(u) = \int_{\R_+}(1 - e^{-us}) \rho(s) \dd s < +\infty
\end{equation}
for all $u>0$, ensuring $0 < \mutilde (\mathbb{S}) < +\infty$ almost surely \citep{Pit(03),Reg(03)}. We write $\mutilde \sim \mbox{CRM}(\theta, \rho, G_0)$ to denote a homogeneous CRM on $\mathbb{S}$. A PK prior is the law of a ``normalization" of a CRM with respect to its total mass \citep[Chapter 4]{Pit(06)}.

\begin{defi}\label{pk}
Let $\mutilde \sim \mbox{CRM}(\theta, \rho, G_0)$ with total mass $T = \mutilde(\mathbb{S}) \sim f_T$. Let $P_k := J_k / T$ be the normalized random jumps of $\mutilde$, and denote by $\text{PK}(\theta, \rho \mid T=t)$ the conditional distribution of $(P_k)_{k\geq1}$ given $T=t$. If $g(\mutilde) \equiv g(T)$ such that $\E[g(T)] = 1$, then a PK prior with parameters $(\theta, \rho,  g f_T, G_0)$ is the law of the (discrete) random probability measure
$P(\cdot) = \sum_{k\geq1} \tilde P_k \delta_{W_{k}}(\cdot)$, on $\mathbb{S}$, where $(\tilde P_k)_{k\geq1}$ is distributed as the PK distribution $\text{PK}(\theta, \rho, g f_T) = \int_{\R_+} \text{PK}(\theta, \rho \mid T=t) g(t) f_T(t) \dd t$, and the $W_{k}$'s are independent and identically distributed as $G_{0}$.
\end{defi}

\subsection{General posterior distribution for PK priors} \label{sec:frequency-general}

Under the model in \eqref{eq:exchangeable_model_hash}, with $\mathscr{P}$ being a PK prior, i.e., $P\sim\text{PK}(\theta, \rho,  g f_T, G_0)$, the next theorem gives the posterior distribution of $f_{X_{n+1}}$ given $\mathbf{C}_{J}$ and $h(X_{n+1})$.  This is result can be applied upon suitable specifications of the measure $\rho$ and the function $g$. For instance, PK priors include the class of (homogeneous) normalized CRM priors, obtained by setting $g$ as the identity function \citep{Jam(02),Pru(02),Pit(03),Reg(03)}. It is sufficient to consider $g(t) \propto t^{-\gamma} e^{-\beta t}$ to recover from Definition~\ref{pk} the most popular priors in BNPs. Common choices of $\mutilde$ are the Gamma CRM and the $\alpha$-Stable CRM \citep{Kin(75)}, as they provide PK priors with a flexible tail behaviour, ranging from geometric tail to heavy power-law tails, respectively \citep[Chapter 3 and Chapter 4]{Pit(06)}. Under the Gamma CRM, Definition~\ref{pk} generalizes the DP prior, which is obtained by setting $g$ as the identity function. Under the $\alpha$-Stable CRM, Definition~\ref{pk} provides a generalization of the normalized $\alpha$-Stable prior \citep{Kin(75)}, which is obtained by setting $g$ as the identity function, and it also includes the PYP prior and the normalized Gamma process prior \citep{Jam(02),Lij(05),Lij(07)}. See \cite[Chapter 3]{Pit(06)} and \citet{Lij(10)} for other examples of PK priors. The application of the next theorem to the DP and the PYP priors will be considered below, showing their peculiar features.

\begin{thm}\label{main_pk}
Let $\mathbf{C}_{J}$ be a sketch of $\mathbf{X}_{n}$ under \eqref{eq:exchangeable_model_hash}, with $P \sim \text{PK}(\theta, \rho, g f_T,G_{0})$ and $g(t) \propto t^{-\gamma} e^{-\beta t}$, and consider an additional (unobservable) $X_{n+1}$. 
With $\psi$ defined as in \eqref{eq_cond},
\begin{multline}\label{post_pk}
\prob[f_{X_{n+1}}=l\,|\,\mathbf{C}_{J}=\mathbf{c},h(X_{n+1})=j]\\
\quad=\frac{\theta}{J} \binom{c_j}{l} \frac{\int_{\R_+} u^{n + \gamma}  \phi^{(c_j - l)}(u + \beta)
 \prod_{k\neq j} \phi^{(c_k)}(u + \beta)  \, \kappa(u + \beta, l+1) \dd u
 }{
  \int_{\R_+ } u^{n+\gamma}   \phi^{(c_j + 1)}(u + \beta)
  \prod_{k\neq j} \phi^{(c_k)}(u + \beta)  \dd u},
\end{multline}
for all $l \in \{0,1, \ldots, c_j\}$, where $\phi^{(n)}(u) = (-1)^{n} \frac{\dd^n}{\dd u^n} e^{-\theta/J \psi(u)}$ and \\ $\kappa(u, n) = \int_{\R_+} e^{-us} s^n \rho(s) \dd s$.
\end{thm}

See Appendix~\ref{app:proof_main} for the proof of Theorem~\ref{main_pk}. This result provides a BNP solution to the frequency recovery problem, as it leads to an estimator of $f_{X_{n+1}}$, with respect to a suitable loss function, as a mean functional of \eqref{main_pk}. Credible intervals may be also derived. The use of the squared loss leads to the posterior mean as an estimator of  $f_{X_{n+1}}$, i.e.,
\begin{equation}\label{eq:post-mean}
\hat{f}_{X_{n+1}}=\sum_{l=0}^{c_j}l \cdot \text{Pr}[f_{X_{n+1}}=l\,|\,\mathbf{C}_{J}=\mathbf{c},h(X_{n+1})=j].
\end{equation}
Regarding the use of the sketch, \citet[Proposition 1]{Cai(18)} and \citet[Theorems 1 and 2]{Dol(23)} provide posterior distributions of $f_{X_{n+1}}$ with respect to the sole bucket $C_{h(X_{n+1})}$ in which $X_{n+1}$ is hashed. By contrast, Theorem~\ref{main_pk}  considers the entire sketch $\mathbf{C}_{J}$ and the hash bucket of $X_{n+1}$, denoted as $h(X_{n+1})$. Therefore, our posterior distribution provides a principled BNP estimator also in those situations where $C_{h(X_{n+1})}$ may not be a sufficient statistic, as elaborated below. Regarding the specification of the prior distribution, \citet{Cai(18)} and \citet{Dol(23)} focused on the DP and PYP priors, obtaining posterior distributions through conjugacy or quasi-conjugacy properties. 
However, Theorem~\ref{main_pk} considers a general PK prior. As discussed in \citet{Dol(23)}, relying on conjugacy poses a limitation when aiming to consider more diverse prior distributions because the DP and PYP priors are the only quasi-conjugate PK priors.  Our proof of Theorem~\ref{main_pk} overcomes this limitation by avoiding the use of any form of conjugacy for the prior.

\subsection{Results under the DP prior}\label{sec:dp_freq}

We consider Definition~\ref{pk} with $\mutilde$ being a Gamma CRM, i.e., $\rho(s)=s^{-1}\exp\{-s\}$, and $g$ the identity function. 
Then, $P\sim\text{PK}(\theta, \rho,  g f_T, G_0)$ is a DP with mass $\theta>0$ and base measure $G_0$ \citep{Fer(73),Pit(03)}; in short, $P\sim\text{DP}(\theta,G_{0})$. In this case, a simplified expression for the posterior of $f_{X_{n+1}}$ given $\mathbf{C}_{J}$ and $h(X_{n+1})$ is obtained from Theorem~\ref{main_pk} as follows.  Denote by $(a)_{(n)}$ the rising factorial of $a$ of order $n$, i.e., $(a)_{(n)}=\prod_{0\leq i\leq n-1}(a+i)$ with the proviso $(a)_{(0)}:=1$ \citep[Chapter 2]{Cha(05)}. for any $l=0,1,\ldots,c_{j}$,
\begin{equation}\label{post_dp}
	\text{Pr}[f_{X_{n+1}}=l\,|\,\mathbf{C}_{J}=\mathbf{c},h(X_{n+1})=j] = \frac{\theta}{J}\frac{(c_{j}-l+1)_{(l)}}{\left(\frac{\theta}{J}+c_{j}-l\right)_{(l+1)}}.
\end{equation}
See Appendix~\ref{app:dp_cor} for a proof of how \eqref{post_dp} follows from \Cref{main_pk}, and Appendix~\ref{app:dp_teo} for an alternative proof based on finite-dimensional properties of the DP. It is easy to see that \eqref{post_dp} is a Beta-Binomial distribution, namely a Binomial distribution in which the probability of success at each of the $c_{j}$  trials is a Beta random variable with parameter $(1,\theta/J)$, say $B_{1,\theta/J}$. That is, if $F_{X_{n+1}}$ is a random variable distributed as \eqref{post_dp}, then de Finetti's theorem implies that $c_{j}^{-1}F_{X_{n+1}}$ converges (weakly) to a $B_{1,\theta/J}$ as $c_{j}\rightarrow+\infty$. The posterior distribution in \eqref{post_dp} depends on the sketch $\mathbf{C}_{J}$ only through $C_{h(X_{n+1})}$. 
Thus, under the DP prior, $C_{h(X_{n+1})}$ is a sufficient statistic for estimating $f_{X_{n+1}}$ from $\mathbf{C}_{J}$, making \eqref{post_dp} equivalent to the posterior distribution in \citet[Proposition 1]{Cai(18)}. 

The next theorem, proved in Appendix~\ref{app:char_proof}, characterizes the DP prior as the unique PK prior for which the posterior distribution of $f_{X_{n+1}}$ with respect to $C_{h(X_{n+1})}$ is equivalent to the posterior distribution with respect to $\mathbf{C}_{J}$ and $h(X_{n+1})$.

\begin{thm}\label{char_dp}
The DP prior is the sole PK prior for which \eqref{post_pk} depends on $\mathbf{C}_{J}$ only through $C_{h(X_{n+1})}$.
\end{thm}

In practice, evaluating \eqref{post_dp} requires estimating the unknown parameter $\theta>0$ of the prior from the sketch $\mathbf{C}_{J}$. \citet{Cai(18)} proposed an empirical Bayes approach to estimate $\theta$, which relies on the following finite-dimensional projective property of $P\sim\text{DP}(\theta,G_{0})$:  if $\{B_1,\ldots,B_{k}\}$ is a measurable $k$-partition of $\mathcal{S}$, for $k\geq1$, then $(P(B_1), \ldots, P(B_k))$ follows a Dirichlet distribution with parameter $(\theta G_{0}(B_{1}),\ldots,\theta G_{0}(B_{k}))$ \citep{Fer(73),Reg(01)}. Due to the finite-dimensional projective property and the assumption that $\mathcal{H}_{J}$ is independent of $P$, the sketch $\mathbf{C}_{J}$ is distributed as a Dirichlet-Multinomial distribution, i.e.,
\begin{align}\label{marg}
&\text{Pr}[\mathbf{C}_{J}=\mathbf{c}]=\frac{n!}{(\theta)_{(n)}}\prod_{j=1}^{J}\frac{(\frac{\theta}{J})_{(c_{j})}}{c_{j}!}.
\end{align}
This distribution enables estimating $\theta$ directly, by maximizing the (marginal) likelihood \eqref{marg} over $\theta$  \citep{Cai(18)}. 
The estimated value of $\theta$ can then be plugged into the posterior distribution \eqref{post_dp}.  A fully Bayesian approach is also possible, by assigning a prior distribution to $\theta$ and assessing the resulting posterior distributions through Monte Carlo sampling.

\subsection{Results under the PYP prior}

Consider Definition~\ref{pk} with $\theta=1$, $\mutilde$ being an $\alpha$-Stable CRM, i.e., $\rho(s)=(\alpha/\Gamma(1-\alpha)) s^{-1-\alpha}$ for $\alpha\in(0,1)$ and $g(t)=(\Gamma(\gamma+1)/\Gamma(\gamma/\alpha+1)t^{-\gamma}$ for $\gamma>-\alpha$, where $\Gamma(\cdot)$ denotes the Gamma function, namely $\Gamma(x)=\int_{(0,+\infty)}z^{x-1}\text{e}^{-z}\ddr z$ for $x>0$. Then, $P\sim\text{PK}(\theta, \rho,  g f_T, G_0)$ is a PYP with discount $\alpha$, mass $\gamma$ and base measure  $G_0$ \citep{Pit(03)}; in short, we write $P\sim\text{PYP}(\alpha,\gamma,G_{0})$. 
In this case, a simplified expression for the posterior of $f_{X_{n+1}}$ given $\mathbf{C}_{J}$ and $h(X_{n+1})$ is obtained from Theorem~\ref{main_pk} as follows. For $n\geq0$ and $0\leq k\leq n$,
\begin{displaymath}
	\mathscr{C}(n,k;\alpha)=\frac{1}{k!}\sum_{i=0}^{k}(-1)^{i}{n\choose i}(-i\alpha)_{(n)}
\end{displaymath}
denotes the generalized factorial of $n$ of order $k$, with $\mathscr{C}(0,0;\alpha):=0$, $\mathscr{C}(0,0;\alpha):=0$, and  $\mathscr{C}(n,0;\alpha):=1$ \citep[Chapter 2]{Cha(05)}. Then, from \eqref{post_pk}, for any $l=0,1,\ldots,c_{j}$,
\begin{align}\label{post_pyp}
&\text{Pr}[f_{X_{n+1}}=l\,|\,\mathbf{C}_{J}=\mathbf{c},h(X_{n+1})=j]\\
&\notag \quad= \frac{\gamma}{J} \binom{c_j}{l} (1 - \alpha)_{(l)} \frac{\sum_{\bm i \in \mathcal S(\bm c, j, -l)} \frac{\Gamma\left(\frac{\gamma + \alpha}{ \alpha} + |\bm i|\right)}{J^{|\bm i|}} \prod_{k=1}^J \mathscr{C}(c_k - l \delta_{k,j}, i_k; \alpha)}{\sum_{\bm i \in \mathcal S(\bm c, j, 1)}\frac{ \Gamma\left(\frac{\gamma}{\alpha} + |\bm i|\right)}{J^{|\bm i|}}  \prod_{k=1}^J \mathscr{C}(c_k + \delta_{k,j}, i_k; \alpha)},
\end{align}
where $\mathcal S(\bm c, j, q)$ is the Cartesian product $\times_{1\leq k\leq J}\{0, \ldots, c_k + \delta_{k, j} q\}$, while $\delta_{k,j}$ is the Kronecker delta and $|\bm i| = \sum_{1\leq k\leq J} i_k$. 
See Appendix~\ref{app:post_pyp} for the proof of \eqref{post_pyp}. 
The posterior distribution  \eqref{post_pyp} generalizes \eqref{post_dp}, which is recovered for $\gamma=\theta$ and $\alpha\rightarrow0$; see Appendix~\ref{app:post_pyp_new}.

Unlike that in \eqref{post_dp}, the posterior distribution in \eqref{post_pyp} relies on the entirety of $\mathbf{C}_{J}$. Further, the size of $\mathcal S(\bm c, j, q)$ increases exponentially with $J$ and $n$. For instance, if $J = 10$ and $c_j = 5$ for all $j$, then $|\mathcal S(\bm c, j, q) | \approx 60 \times 10^6$. Consequently, the evaluation of \eqref{post_pyp} is intractable even for moderately large $n$, as it necessitates summations over an exceedingly large number of generalized factorial coefficients, depending on $J$. See  Appendix~\ref{app:computations_pyp} for some approaches to evaluate \eqref{post_pyp}, which still lead to non-trivial computational obstacles. To overcome this challenge, we seek an approximation of \eqref{post_pyp} that depends on $\mathbf{C}_{J}$ only through $C_{h(X_{n+1})}$. The next theorem characterizes the large-sample behaviour of the posterior mean estimator~\eqref{eq:post-mean}.

\begin{thm}\label{char_pyp}
Let $\mathbf{C}_{J}$ be a sketch of $\mathbf{X}_{n}$ under the model in \eqref{eq:exchangeable_model_hash}, with $P \sim \text{PYP}(\alpha,\gamma,G_{0})$ for $\alpha > 0$, and consider an additional (unobservable) $X_{n+1}$. Suppose $h(X_{n+1})=j$, for some fixed $j \in \{1,\ldots,J\}$.
Define $\mathbf{c}_{-j} := (c_1, \ldots, c_{j-1}, c_{j+1}, \ldots, c_J)$. 
Let $\hat{f}_{X_{n+1}}$ denote the posterior mean estimator as defined in~\eqref{eq:post-mean}.
Then,
\begin{align*}
&\lim_{\mathbf{c}_{j}\rightarrow+\infty} \lim_{\mathbf{c}_{-j}\rightarrow+\infty} \frac{\hat{f}_{X_{n+1}}}{c_j} 
 = \frac{\gamma}{\alpha} \cdot \frac{1 - \alpha}{\gamma + J\alpha - \alpha + 1}.
\end{align*}
\end{thm}
See Appendix~\ref{asypyp} for the proof of Theorem~\ref{char_pyp}. 
This result motivates approximating the posterior mean estimator $\hat f_{X_{n+1}}$, in those situations where all the $c_{i}$'s are large, with:
\begin{equation}\label{approx_est}
\tilde f_{X_{n+1}} := c_j \cdot \frac{\gamma}{\alpha} \cdot \frac{1 - \alpha}{\gamma + J\alpha - \alpha + 1}.
\end{equation}

Applying either \eqref{post_pyp} or \eqref{approx_est} requires estimating $(\alpha,\gamma)$ from the sketch $\mathbf{C}_{J}$. Under the PYP prior, the distribution of $\mathbf{C}_{J}$ has no simple closed-form expression, which impedes the use of empirical Bayes or fully Bayesian approaches to estimate $(\alpha,\gamma)$. To address a similar challenge, \cite{Dol(23)} adopted a likelihood-free approach based on the Wasserstein distance \citep{Ber(19)}. That involves sampling independent data sets $\mathbf{X}^{\prime}_n$ from $P\sim\text{PYP}(\alpha^\prime,\gamma^\prime,G_{0})$ and sketching them into $\mathbf{C}^{\prime}_{J}$ using the same hash function $h$. 
Then, $(\alpha,\gamma)$ are estimated by minimizing (a suitable approximation of) the 1-Wasserstein distance between $\mathbf{C}^{\prime}_{J}$ and $\mathbf{C}_{J}$. Since the objective function is not differentiable, \cite{Dol(23)} utilized Bayesian Optimization to estimate the parameters. 

However, their approach is computationally demanding: the cost of evaluating the Wasserstein distance is $\mathcal O(J^3)$, while the cost of simulating $\mathbf{X}^{\prime}_n$ scales super-linearly with $n$, depending on the parameters' values. 
Therefore, we propose an alternative approach: we  consider the first $n^\prime \ll n$ observations, sketch them through $h$, and estimate their frequencies using \eqref{approx_est}. Then, the mean absolute error of the frequency recovery can be easily minimized with respect to the unknown parameters using standard software packages since the loss function is differentiable almost everywhere.

\subsection{Comparison of the posterior under the PYP and DP priors}

We present a numerical illustration of the posterior distributions under the DP and PYP priors; see Section~\ref{sec:numerics} for a more extensive simulation study. We set $n=50$ and $J=10$, considering four sketched datasets: $\mathbf C_{J}^{\text{\tiny{(1)}}}$, $\mathbf C_{J}^{\text{\tiny{(2)}}}$, $\mathbf C_{J}^{\text{\tiny{(3)}}}$, and $\mathbf C_{J}^{\text{\tiny{(4)}}}$, with values reported in Table~\ref{tab:sketch_simu}. Specifically: i) in scenario $\mathbf C_{J}^{\text{\tiny{(1)}}}$, the values $C^{\text{\tiny{(1)}}}_j$'s are constant across $j$; ii) in scenario $\mathbf C_{J}^{\text{\tiny{(2)}}}$, the values $C^{\text{\tiny{(2)}}}_j$'s decay exponentially in $j$; iii) in scenario $\mathbf C_{J}^{\text{\tiny{(3)}}}$, the values $C^{\text{\tiny{(3)}}}_j$'s decay linearly in $j$; iv) in scenario $\mathbf C_{J}^{\text{\tiny{(4)}}}$, the values $C^{\text{\tiny{(4)}}}_j$'s are either nine, five, or one. Additionally, we assume $X_{n+1}$ is mapped into bucket $h^{(i)}(X_{n+1})$ such that $C^{\text{\tiny{(i)}}}_{h^{(i)}(X_{n+1})} = 5$ for $i=1,2,3,4$. We consider a PYP prior with parameter $\gamma=1$ and parameter $\alpha = 0, 0.1, 0.3, 0.5$; recall that the PYP prior with $\alpha=0$ and $\gamma>0$ coincides with the DP prior with parameter $\gamma$. 

\begin{figure}[!htb]
  \centering
  \includegraphics[width=12cm,height=6cm,keepaspectratio]{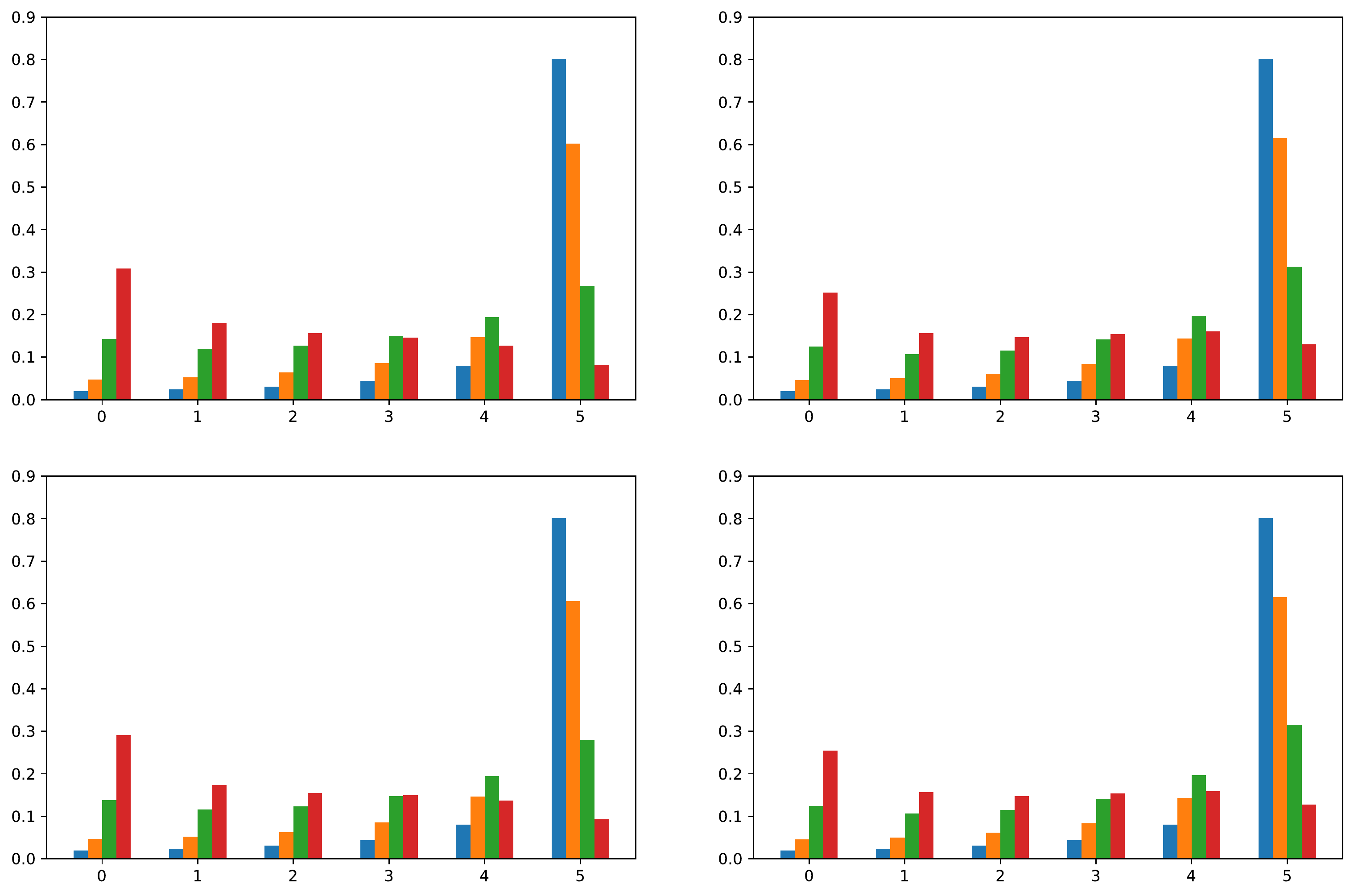}
  \caption{\small{Posterior distribution of $f_{X_{n+1}}$ in the four scenarios of Table~\ref{tab:sketch_simu}. The blue, orange, green, and red bars correspond to a PYP$(\alpha,\gamma)$ prior, with $\alpha$ equal to $0$, $0.1$, $0.3$, and $0.5$ respectively, and $\gamma = 1$. Since the DP posterior ($\alpha=0$) concentrates more of its probability mass at 5, the results suggest the DP prior leads to a more accurate estimate across these 4 scenarios. } }
  \label{fig:py_post}
\end{figure}

Figure~\ref{fig:py_post} summarizes the posterior distribution of $f_{X_{n+1}}$ in each scenario. Given the sufficiency of $C_{h(X_{n+1})}$ under the DP prior, the corresponding posterior distributions remain unchanged across all scenarios. Moreover, the posterior distributions under the DP prior are the most concentrated on larger values. Across all the scenarios, increasing $\alpha$ pushes the posterior towards lower values. Although not clearly evident from the plots, there are slight differences under the PYP priors across scenarios. Particularly, for $\alpha = 0.3$ and $\alpha=0.5$, the posterior mass assigned to 5 is larger in the second scenario compared to other settings.

\begin{table}[!htb]
\centering
  \begin{tabular}{c | c c c c c c c c c c}
     $\mathbf C_{J}^{\text{\tiny{(1)}}}$ & 5 & 5 & 5 & 5 & 5 & 5 & 5 & 5 & 5 & 5  \\
     $\mathbf C_{J}^{\text{\tiny{(2)}}}$ & 14 & 10 & 7 & 5 & 4 &3 & 2 & 2 & 2 & 1   \\
     $\mathbf C_{J}^{\text{\tiny{(3)}}}$ & 10 & 9 & 8 & 7 & 5 & 4 & 3 & 2 & 1 & 1  \\
      $\mathbf C_{J}^{\text{\tiny{(4)}}}$ & 9 & 9 & 9 & 5 & 5 & 5 & 5 & 1 & 1 & 1  \\
  \end{tabular}
  \caption{\label{tab:sketch_simu} \small{Empirical frequencies of the $j$-th bucket (by column) in 4 simulated scenarios.}}
  \end{table}

\section{BNP cardinality recovery}\label{sec21}

\subsection{Setup and overview}

In the ``species'' setting described in Section~\ref{sec2}, we tackle the problem of estimating the cardinality $k_{n}$ of $(x_{1},\ldots,x_{n})$ and the $l$-cardinality $m_{l,n}$ of $(x_{1},\ldots,x_{n})$, for any $l\in[n]$, based on the data sketch $\mathbf{C}_{J}$. Estimating the $m_{l,n}$'s is a refinement of the problem of estimating $k_{n}$, as the estimates of the $m_{l,n}$'s imply an estimate of $k_{n}$, given by
\begin{align*}
  & k_{n}=\sum_{l=1}^{n}m_{l,n},
  & n=\sum_{l=1}^{n}l \cdot m_{l,n}.
\end{align*}
Recovering the $l$-cardinalities implies recovering the partition structure of $(x_{1},\ldots,x_{n})$, which is the partition of $\{1,\ldots,n\}$ induced by the equivalence relation $i\sim j\iff x_{i}=x_{j}$; this is a sufficient statistic for the sample $(x_{1},\ldots,x_{n})$.  Under the model \eqref{eq:exchangeable_model_hash}, for the DP and the PYP priors, we demonstrate that recovering $m_{l,n}$ boils down to computing the posterior distribution of $f_{X_{n+1}}$ given $\mathbf{C}_{J}$. 
Thus, the BNP approach connects frequency and cardinality recovery, enabling the estimation of $f_{X_{n+1}}$ and $k_{n}$ from the same sketch $\mathbf{C}_{J}$. Our approach relies on a distinctive feature of the DP and PYP priors, known as the ``sufficiency'' postulate \citep{Bac(17)}, and does not extend to other PK priors.  Below, we compute estimates of $k_{n}$ and of the $m_{l,n}$'s in terms of posterior expectations, given $\mathbf{C}_{J}$.

\subsection{Background on the PYP prior}

Before stating our main results, we recall the distribution of a random sample $\mathbf{X}_{n}=(X_{1},\ldots,X_{n})$ from $P\sim\text{PYP}(\alpha,\gamma,G_{0})$. Due to the (almost sure) discreteness of $P$, the random sample $\mathbf{X}_{n}$ induces a random partition of $[n]$ into $K_{n}=k\leq n$ blocks, denoted as
$\{S_{1}^{\ast},\ldots,S^{\ast}_{K_{n}}\}$, with frequencies $(N_{1,n},\ldots,N_{K_{n},n})=(n_{1},\ldots,n_{k})$ such that $n_{i}>0$ and $\sum_{1\leq i\leq k}n_{i}=n$. 
Let $M_{l,n}$ be the number of blocks with frequency $l\in[n]$, i.e.,
\begin{displaymath}
M_{l,n}=\sum_{i=1}^{K_{n}}I(N_{i,n}=l),
\end{displaymath}
so that $\sum_{l=1}^{n}M_{l,n}=K_{n}$ and $\sum_{l=1}^{n}lM_{l,n}=n$, and let $\mathbf{M}_{n}=(M_{1,n},\ldots,M_{n,n})$. If we set
\begin{displaymath}
\mathcal{M}_{n}=\left\{(m_{1},\ldots,m_{n})\text{ : }m_{i}\geq0,\,\sum_{l=1}^{n}m_{l}=k,\,\sum_{l=1}^{n}lm_{l}=n\text{ and }k\in[n]\right\},
\end{displaymath}
then, for $\mathbf{m}\in\mathcal{M}_{n}$,
\begin{equation}\label{eq_ewe_py}
\text{Pr}[\mathbf{M}_{n}=\mathbf{m}]=n!\frac{\left(\frac{\gamma}{\alpha}\right)_{(\sum_{l=1}^{n}m_{l})}}{(\gamma)_{(n)}}\prod_{l=1}^{n}\left(\frac{\alpha(1-\alpha)_{(l-1)}}{l!}\right)^{m_{l}}\frac{1}{m_{l}!}.
\end{equation}
The distribution in \eqref{eq_ewe_py} first appeared in \citet[Proposition 9]{Pit(95)}. Assuming $P\sim\text{DP}(\theta,G_{0})$, the distribution of $\mathbf{M}_{n}$ is obtained from \eqref{eq_ewe_py} by setting $\gamma=\theta$ and letting $\alpha\rightarrow0$.

An application of \eqref{eq_ewe_py} yields the conditional distribution of $X_{n+1}$ given $\mathbf{X}_{n}$, known as the predictive distribution or generative scheme of the PYP prior, and also of the DP prior by letting $\alpha \rightarrow 0$ \citep[Chapter 3]{Pit(06)}.
 In particular, for any $l\in[n]$, let
\begin{align*}
  & \mathcal{S}_{0}=\mathbb{S}-\{S_{1}^{\ast},\ldots,S^{\ast}_{K_{n}}\},
  & \mathcal{S}_{l}=\bigcup_{i=1}^{K_{n}}\{S^{\ast}_{i}\in\{S_{1}^{\ast},\ldots,S^{\ast}_{K_{n}}\}\text{ : }N_{i,n}=l\}.
\end{align*}
Above, $\mathcal{S}_{0}$ is the set of ``new'' symbols not observed in $\mathbf{X}_{n}$, while $\mathcal{S}_{r}$ is the set of ``old'' symbols observed in $\mathbf{X}_{n}$ with frequency $r$. From \citet[Proposition 9]{Pit(95)}, 
\begin{equation}\label{bnp_est}
\text{Pr}[X_{n+1}\in\mathcal{S}_{l}\,|\,\mathbf{X}_{n}]=\begin{cases}
 \frac{\gamma+k\alpha}{\gamma+n}, &\mbox{ if } l=0, \\[0.4cm]
\frac{m_{l}(l-\alpha)}{\gamma+n}, &\mbox{ if } l\geq1.
\end{cases}
\end{equation}
The PYP prior is the sole PK prior for which the probability that $X_{n+1}$ is a ``new" symbol depends on $\mathbf{X}_{n}$ only through $K_{n}$, and the probability that $X_{n+1}$ is an ``old" symbol with frequency $l$  depends on $\mathbf{X}_{n}$ only through $M_{l,n}$. 
Further, the DP prior is the sole PK prior for which the probability that $X_{n+1}$ is ``new" depends on $\mathbf{X}_{n}$ only through $n$. These are known as the ``sufficientness'' postulates of the DP and PYP priors \citep{Bac(17)}.

\subsection{Main result}

Our BNP approach integrates the results of Section~\ref{sec2}, under both the DP and PYP priors, with the predictive distribution \eqref{bnp_est}. 
We focus here on the PYP prior, noting that the results for the DP follow by setting $\gamma=\theta$ and letting $\alpha\rightarrow0$. 
Before presenting the main results, we outline the key arguments of our approach. In particular, for any $l\in[n]$,
\begin{equation}\label{card_main_1}
\text{Pr}[f_{X_{n+1}}=l\,|\,\mathbf{C}_{J}=\mathbf{c}]=\sum_{\mathbf{m}\in\mathcal{M}_{n}}\text{Pr}[X_{n+1}\in\mathcal{S}_{l}\,|\,\mathbf{C}_{J}=\mathbf{c},\,\mathbf{M}_{n}=\mathbf{m}]\text{Pr}[\mathbf{M}_{n}=\mathbf{m}\,|\,\mathbf{C}_{J}=\mathbf{c}],
\end{equation}
where
\begin{equation}\label{card_main_2}
\text{Pr}[X_{n+1}\in\mathcal{S}_{l}\,|\,\mathbf{C}_{J}=\mathbf{c},\,\mathbf{M}_{n}=\mathbf{m}]=\text{Pr}[X_{n+1}\in\mathcal{S}_{l}\,|\,\mathbf{M}_{n}=\mathbf{m}]=\frac{m_{l}(l-\alpha)}{\gamma+n},
\end{equation}
with the last identity in \eqref{card_main_2} derived from \eqref{bnp_est}. By combining \eqref{card_main_1} with \eqref{card_main_2}, we obtain:
\begin{equation}\label{card_main_3}
\text{Pr}[f_{X_{n+1}}=l\,|\,\mathbf{C}_{J}=\mathbf{c}]=\frac{l-\alpha}{\gamma+n}\E[M_{l,n}\,|\,\mathbf{C}_{J}=\mathbf{c}], \qquad l\in[n].
\end{equation}
This connects the posterior distribution of $f_{X_{n+1}}$, given $\mathbf{C}_{J}$, to the corresponding conditional expectation of $M_{l,n}$. This is the identity underpinning  our approach to cardinality recovery, effectively linking it to the frequency recovery problem discussed in Section~\ref{sec2}.

From \eqref{card_main_3}, a BNP estimator for $m_{l,n}$ under a squared loss can be immediately derived. This is the conditional expectation of $M_{l,n}$, given $\mathbf{C}_{J}$, which is given for any $l\in[n]$ by
\begin{equation}\label{card_main_4}
\hat{m}_{l,n}=\E[M_{l,n}\,|\,\mathbf{C}_{J}=\mathbf{c}]=\frac{\gamma+n}{l-\alpha}\text{Pr}[f_{X_{n+1}}=l\,|\,\mathbf{C}_{J}=\mathbf{c}].
\end{equation}
Then, since $K_{n}=\sum_{l=1}^n M_{l,n}$, a BNP estimator for $k_{n}$ under a squared loss is of the form
\begin{equation}\label{card_main_5}
\hat{k}_{n}=\E[K_{n}\,|\,\mathbf{C}_{J}=\mathbf{c}]=\sum_{l=1}^{n}\frac{\gamma+n}{l-\alpha}\text{Pr}[f_{X_{n+1}}=l\,|\,\mathbf{C}_{J}=\mathbf{c}].
\end{equation}
Both \eqref{card_main_4} and \eqref{card_main_5} are based on \eqref{card_main_3}, which in turn depends on \eqref{card_main_2} relying only on $M_{l,n}$. Since the PYP prior is the sole PK prior for which \eqref{card_main_2} depends on $\mathbf{X}_{n}$ only through $M_{l,n}$, this approach cannot be extended to other PK priors.  To conclude, we compute,
\begin{equation}\label{card_main_6}
\text{Pr}[f_{X_{n+1}}=l\,|\,\mathbf{C}_{J}=\mathbf{c}]=\sum_{j=1}^{J}\text{Pr}[f_{X_{n+1}}=l\,|\,\mathbf{C}_{J}=\mathbf{c},h(X_{n+1})=j]\text{Pr}[h(X_{n+1})=j\,|\,\mathbf{C}_{j}=\mathbf{c}],
\end{equation}
for any $l\in[n]$, where the distribution of $f_{X_{n+1}}$ conditional on $\mathbf{C}_{j}$ and $h(X_{n+1})$ is given in \eqref{post_pyp}. 
We next provide an expression for \eqref{card_main_6} that simplifies the estimators \eqref{card_main_4} and \eqref{card_main_5}.

\begin{corollary}\label{corol_card}
Let $\mathbf{C}_{J}$ be a sketch of $\mathbf{X}_{n}$ under \eqref{eq:exchangeable_model_hash}, with $P\sim\text{PYP}(\alpha,\gamma,G_{0})$. For $l \in [n]$,
\begin{align}\label{post_card}
\begin{split}
&\text{Pr}[f_{X_{n+1}}=l\,|\,\mathbf{C}_{J}=\mathbf{c}]\\
& \quad=\frac{\frac{\gamma}{J}}{\gamma+n}(1-\alpha)_{(l)}\sum_{j=1}^{J}{c_{j}\choose l}\frac{\sum_{\bm i \in \mathcal S(\bm c, j, -l)} \frac{\Gamma\left(\frac{\gamma + \alpha}{ \alpha} + |\bm i|\right)}{J^{|\bm i|}} \prod_{k=1}^J \mathscr{C}(c_k - l \delta_{k,j}, i_k; \alpha)}{\sum_{\bm{i}\in\mathcal{S}(\mathbf{c},j,0)}\frac{\Gamma\left(\frac{\gamma}{\alpha}+|\bm{i}|\right)}{J^{|\bm{i}|}}\prod_{k=1}^{J}\mathscr{C}(c_{k},i_{k};\alpha)}.
\end{split}
\end{align}
\end{corollary}

See Appendix~\ref{card_proof} for the proof of Corollary~\ref{corol_card}. The BNP estimators for $l$-cardinality and cardinality are derived by combining \eqref{post_card} with \eqref{card_main_4} and \eqref{card_main_5}, respectively.  Note that, under the PYP prior, these estimators face the same computational challenges as those discussed in Section~\ref{sec2} for the frequency estimation problem. 

By contrast, under the DP prior, the estimators simplify significantly.
By setting $\gamma=\theta$ and letting $\alpha\rightarrow0$ in \eqref{post_card}, we prove in Appendix~\ref{card_pyp_pd} that, for any $l\in[n]$,
\begin{equation}\label{post_card_pd}
\text{Pr}[f_{X_{n+1}}=l\,|\,\mathbf{C}_{J}=\mathbf{c}]=\frac{\frac{\theta}{J}}{\theta+n}\sum_{j=1}^{J}\frac{(c_{j}-l+1)_{(l)}}{\left(\frac{\gamma}{J}+c_{j}-l\right)_{(l)}}.
\end{equation}
By combining \eqref{post_card_pd} with  \eqref{card_main_4} and \eqref{card_main_5}, with $\gamma=\theta$ and $\alpha=0$, for any $l\in[n]$:
\begin{align*}
& \hat{m}_{l,n}=\frac{\frac{\theta}{J}}{l}\sum_{j=1}^{J}\frac{(c_{j}-l+1)_{(l)}}{\left(\frac{\gamma}{J}+c_{j}-l\right)_{(l)}}
& \hat{k}_{n}=-\theta\psi\left(1-\frac{\theta}{J}\right)+\frac{\theta}{J}\sum_{j=1}^{J}\psi\left(1-\frac{\theta}{J}-c_{j}\right),
\end{align*}
where $\psi$ is the digamma function, i.e. $\psi(x)=\frac{d}{dx}\log\Gamma(x)$. These estimates for $k_{n}$ and $m_{l,n}$ depend on the prior parameter $\theta$, which can be estimated from $\mathbf{C}_{J}$ as discussed in Section~\ref{sec2}.

\section{BNP frequency recovery in the ``traits" setting} \label{sec3}

The ``traits'' setting for frequency recovery extends the ``species" setting studied in Section~\ref{sec2} by allowing the data points to exhibit nonnegative integer levels of association with multiple symbols. 
We consider $n \geq 1$ data points $(x_1, \ldots, x_n)$ modeled as a random sample $\mathbf{X}_{n} = (X_1, \ldots, X_n)$, where each $X_i$ is represented as $X_i = ((\tilde Y_{i,j}, \tilde A_{i,j}), j = 1, \ldots, K_i)$. Here, $\tilde A_{i,j} \in \mathbb{N}_{0}$ represents the level of association of the $i$-th sample with the trait $\tilde Y_{i,j}$. Sketching $\mathbf{X}_{n}$ is conducted at the trait level using a random hash function $h: \mathbb{S} \rightarrow \{1, \ldots, J\}$ from the hash family $\mathcal{H}_{J}$. This function assigns each $\tilde Y_{i,j}$ to a bucket, incrementing the bucket's counter $C_{h(\tilde Y_{i,j})}$ by $\tilde A_{i,j}$. Consequently, the sketch $\mathbf{C}_{J} = (C_1, \ldots, C_J)$ captures the total association levels for the traits hashed into each bucket. For a new data point $X_{n+1}$, we define the empirical frequency level of a trait $Y_{n+1, r}$ as:
\begin{equation}\label{eq:trait_rec1}
  f_{Y_{n+1, r}} = \sum_{i=1}^n \sum_{j=1}^{K_i} \tilde A_{i, j} I(\tilde Y_{i,j} = Y_{n+1, r}).
\end{equation}
We make use of the sketch $\mathbf{C}_{J}$ to develop a BNP approach to estimate $f_{Y_{n+1, r}}$, assuming that $\mathbf{X}_{n}$ is sampled from a generalized Indian buffet process \citep{Jam(17)}. Under this model, the frequency levels of traits are sufficient statistics, thus paralleling the role of species frequencies in our earlier BNP model for species data in Section~\ref{sec2}.

\begin{figure}[!htb]
\centering
\includegraphics[width=12cm,height=6cm,keepaspectratio]{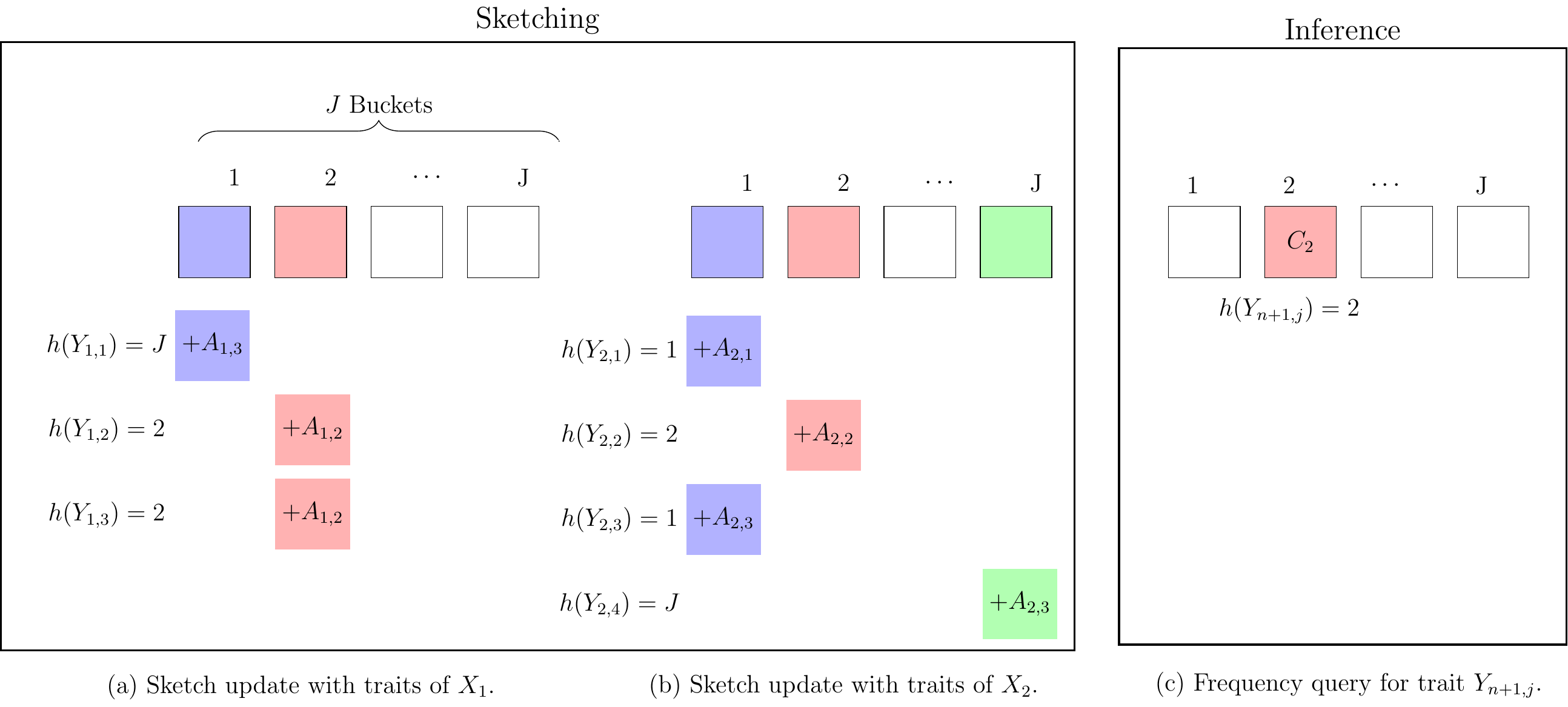}
\caption{\small{Data sketching in the ``traits'' setting. Each $X_i$, for $i \in \{1,\ldots,n\}$, is defined by $K_i$ traits $(Y_{i,1}, \ldots, Y_{i,K_i})$ and exhibits an association level $A_{i,j}$ with each $Y_{i,j}$. A hash function maps each trait into one of $J$ possible ``buckets'' (shown in different colors), and the corresponding counters are incremented by one. The recovery problem is to estimate the empirical frequency of a ``trait'' $X_{Y+1,j}$ based on the information contained in the sketch and the bucket assignments for $Y_{n+1,j}$.}} \label{fig:diagram-traits}
\end{figure}

Topic modeling is a prominent application of ``traits'' allocations, particularly through the use of a multinomial naive Bayes classifier to categorize documents into topics. Suppose each data point $(X_i, T_i)$ consists of a document and a topic label $T_i \in \{1, \ldots, M\}$. If $n_j$ is the number of documents in each topic, the naive Bayes classification rule is given by:
\begin{equation}\label{eq:multnb}
  \prob[T_i = m \mid X_i = \{(y_{i,j}, a_{i,j})\}_{j=1}^{K_i}] \propto \prob[T_i = m]\prod_{j=1}^{K_i} (\pi^m_{y_{i,j}})^{a_{i, j}},
\end{equation}
where $\pi^m_{y_{i,j}}$ is the relative frequency of the trait $y_{i,j}$ in documents with topic $m$. For sketched data, $\pi^m_{y_{i,j}}$ can be replaced by the Bayesian estimator $\hat{f}^m_{y_{i,j}} / n_j$. Another application is in feature engineering, computing sketched ``tf-idfs''. In information retrieval, tf-idf is a preprocessing to adjust the $\tilde A_{i,k}$'s based on the frequency of documents with trait $\tilde Y_{i,k}$:
\begin{equation}\label{tf-idfs}
    \frac{\tilde A_{i,k}}{\sum_{j=1}^{K_i} \tilde A_{i, j}} \log \frac{n}{\sum_{i=1}^n I(\tilde Y_{i, k} \in X_i)}.
\end{equation}
It is known that tf-idfs can replace raw frequencies in \eqref{eq:multnb} \citep{Ren(03)}. One may consider a sketched version of \eqref{tf-idfs}, where the counter is incremented by one instead of $\tilde A_{i, k}$, so that $f_{Y_{n+1}, r}$ is the number of documents containing the word $Y_{n+1}$. Lastly, \cite{Zho(16)} combines elements from the naive Bayes approach with a BNP model for ``traits'' allocations, further illustrating the adaptability and potential of these methods.

\subsection{BNP model and main result}

We recall the definition of the generalized Indian buffet process; see also \citet{Bro(15)}, \citet{Bro(18)}.
\begin{defi}
  Let $\mutilde \sim \text{CRM}(\theta, \rho, G_0)$, that is $\mutilde= \sum_{k \geq 1} J_k \delta_{W_k}$, and let $G_A$ be a probability mass function over $\N_0$. We say that a random variable $X$ given $\mutilde$ is distributed as a generalized Indian buffet process with parameter $G_A$ if
$X = \sum_{k \geq 1} A_{k} \delta_{W_k}$,
where $(A_{k})_{k\geq1}$ is independent of $(W_{k})_{k\geq1}$ and such that $A_{k} \mid J_k \sim G_A(J_k)$ for $k\geq1$.
\end{defi}

We write $X \mid \mutilde \sim \mbox{IBP}(G_A \mid \mutilde)$ to denote that $X$ is distributed according to a generalized Indian buffet process with parameter $G_A$. While $X$ represents an (infinite) random measure, it is assumed that only a finite number of the $A_{i,k}$'s are non-zero. The representation of $X$ as a collection of displayed traits with their corresponding levels of association follows by letting $\{(\tilde Y_{ j}, \tilde A_{j})\}_j =  \{(W_k \, : A_{k} > 0, A_{k})\}_k$. Accordingly, we write the BNP model as:
\begin{align}\label{eq:trait_model_hash}
  \begin{split}
    C_{j} &\,=\, \sum_{i=1}^{n} \sum_{k \geq 1} A_{i,k} I( h(W_k) = j), \qquad j\in[J], \\
    h&\,\sim\,\mathcal{H}_{J}, \\
    X_1,\ldots,X_{n}\,|\,\mutilde &\,\simiid\, \mbox{IBP}(G_A \mid \mutilde), \\
    \mutilde&\,\sim\,\mbox{CRM}(\theta, \rho, G_0).
  \end{split}
\end{align}
Under this model \eqref{eq:trait_model_hash},
\begin{displaymath}
  f_{Y_{n+1, r}} = \sum_{i=1}^n X_i(Y_{n+1, r})
\end{displaymath}
is the total weighted frequency from equation \eqref{eq:trait_rec1}. The sketch $\mathbf{C}_{J}$ is derived through random hashing of $\mathbf{X}_{n}$, with each bucket $j$ containing the hashed traits of $\mathbb S$ where $h(\omega) = j$, denoted as $D_j = h^{-1}({j}) = \{\omega \in \mathbb S \, : h(\omega) = j\}$. 
For a new data point $X_{n+1}$, let $Y_{n+1, r}$ be the trait belonging to $X_{n+1}$ whose frequency we are interested in, and let $B_j = X_{n+1}(D_j)$ be the increment to the $j$-th bucket of the sketch given by $X_{n+1}$, for $j\in[J]$. Let $\mathbf B = (B_1, \ldots, B_J)$.

In the ``species'' setting, the posterior distribution \eqref{post_pk} represented the probability that $X_{n+1}$ appeared $l+1$ times conditional on the augmented sketch $(C_1, \ldots, C_j + 1, \ldots, C_J)$ and $h(X_{n+1}) = j$. This has a natural counterpart in the ``traits'' setting, namely
\begin{equation}\label{eq:trait_rec_all}
  \prob\left[f_{Y_{n+1,r}} = l, X_{n+1}(Y_{n+1,r}) = a \mid h(Y_{n+1,r}) = j, \mathbf C_J = \mathbf c_J, \mathbf B = \mathbf b \right],
\end{equation}
from which the posterior distribution of $f_{Y_{n+1,r}}$, given $\mathbf{C}_{J}$, $\mathbf{B}$, and $h(Y_{n+1,r})$, follows by Bayes' theorem. The next proposition establishes that the pair $(C_{h(Y_{n+1, r})}, B_{h(Y_{n+1, r})})$ is a sufficient statistic for the estimation of $f_{Y_{n+1,r}}$, with respect to $(h(Y_{n+1,r}), \mathbf C_J, \mathbf B)$.
\begin{prp}\label{prop:rec_trait_all}
Let $\mathbf{C}_{J}$ be a sketch of $\mathbf{X}_{n}$ under the model \eqref{eq:trait_model_hash}. Assuming $X_{n+1}$ is an additional, unobserved random sample, for $l = 0, \ldots, c$, the following holds:
  \begin{align}\label{eq:trait_rec_single}
    \begin{split}
   &\prob\left[f_{Y_{n+1,r}} = l, X_{n+1}(Y_{n+1,r}) = a \mid h(Y_{n+1,r}) = j, \mathbf C_J = \mathbf c_J, \mathbf B = \mathbf b \right]\\
      & \quad=\prob\left[f_{Y_{n+1,r}} = l, X_{n+1}(Y_{n+1,r}) = a \mid h(Y_{n+1,r}) = j,  C_{j} = c, B_{j} = b \right].
    \end{split}
  \end{align}
\end{prp}

See Appendix~\ref{app:prof_trait_all} for the proof of Proposition~\ref{prop:rec_trait_all}. The role of Proposition~\ref{prop:rec_trait_all} is to highlight an interesting distinction between the ``species'' and ``traits'' settings in frequency recovery.  In the ``species" setting, Theorem~\ref{char_dp} establishes that the DP prior is the sole prior for which the posterior distribution of $f_{X_{n+1}}$, given $\mathbf{C}_{J}$ and $h(X_{n+1})$, matches the distribution given just $C_{h(X_{n+1})}$. Conversely, in the ``traits" setting, Proposition~\ref{prop:rec_trait_all} demonstrates that for any CRM prior, the posterior distribution of $f_{Y_{n+1,r}}$, given $\mathbf{C}_{J}$, $\mathbf{B}$, and $h(Y_{n+1,r})$, aligns with that given $C_{h(Y_{n+1,r})}$ and $B_{h(Y_{n+1,r})}$. This interesting phenomenon indicates that all CRM priors in the ``traits" setting satisfy a sufficiency property analogous to the DP prior in the ``species" setting. The following theorem provides an expression for \eqref{eq:trait_rec_single}.
\begin{thm}\label{teo:trait_general}
Let $\mathbf{C}_{J}$ represent a sketch of $\mathbf{X}_{n}$ under the model \eqref{eq:trait_model_hash} with $X_{n+1}$ being an additional (unobservable) random sample. 
For any $l \in \{0, \ldots, c\}$,
\begin{align}\label{eq:post_trait}
    \begin{split}
  & \prob\left[f_{Y_{n+1,r}} = l, X_{n+1}(Y_{n+1,r}) = a \mid h(Y_{n+1,r}) = j,  C_{j} = c, B_{j} = b \right]\\
  & \quad= \frac{\theta}{J}\frac{\prob\left[ \sum_{k \geq 1} \sum_{i=1}^n A^\prime_{i, k} = c - l, \sum_{k \geq 1}  A^\prime_{n+1, k} = b - a\right]}{\prob\left[ \sum_{k \geq 1} \sum_{i=1}^n A^\prime_{i, k} = c, \sum_{k \geq 1}  A^\prime_{n+1, k} = b\right]}  \\
  & \quad\quad \times  \int_{\R_+} \prob \left[ \sum_{i=1}^n \tilde A_{i} = l, \tilde A_{n+1} = a \mid s \right] \rho(s) \dd s,
\end{split}
\end{align}
where $\tilde A_{1} \ldots, \tilde A_{n+1} \mid s \simiid G_A(\cdot \mid s)$ and the $A^\prime_{ik}$'s are the projection on the second coordinate of the points of $N^\prime = \{(J^\prime_k, (A^\prime_{i,k})_{i=1}^{n+1})\}_{k \geq 1}$, which is a Poisson process on $\R_+ \times \N_0^{n+1}$ with L\'evy intensity
$ (\theta/J) [\prod_{1 \leq 1 \leq n+1} G_A(\dd a_i \mid s) ] \rho(s) \dd s$.
\end{thm}

See Appendix~\ref{app:proof_trait_general} for a proof of Theorem~\ref{teo:trait_general}. The application of Theorem~\ref{teo:trait_general} necessitates defining the CRM prior via its L\'evy intensity and choosing a distribution $G_{A}$ for the trait association levels. Below, we apply Theorem~\ref{teo:trait_general} with $G_{A}$ modeled as a Poisson distribution, which simplifies \eqref{eq:post_trait}, and then we consider modeling  $G_{A}$ as a Bernoulli distribution.

\subsection{Results under $G_{A}$ Poisson}\label{sec:poisson_traits}
Let $G_A(\cdot \mid J_k)$ be the Poisson distribution with parameters $\lambda J_k$, for fixed $\lambda > 0$.  Leveraging the closeness under convolution of Poisson distributions, we can simplify \eqref{eq:post_trait} considerably.

\begin{prp}\label{prop:poisson_traits}
Let $\mathbf{C}_{J}$ be a sketch of $\mathbf{X}_{n}$ under the model \eqref{eq:trait_model_hash} with $G_A(J_k)$ being the Poisson distribution with parameter $rJ_k$ for a fixed $r > 0$, and let $X_{n+1}$ be an additional (unobservable) random sample. Further, let $\psi$ be the function defined in \eqref{eq_cond}, and let $\phi^{(n)}(u) = (-1)^{n} \frac{\dd^n}{\dd u^n} e^{-\theta/J \psi(u)}$ and  $\kappa(u, n) = \int_{\R_+} e^{-un} s^n \rho(s) \dd s$. Then, for $l=0, \ldots, c$,
  \begin{align}\label{post_poiss}
    \begin{split}
    &\prob\left[f_{Y_{n+1,r}} = l, X_{n+1}(Y_{n+1,r}) = a \mid h(Y_{n+1,r}) = j,  C_{j} = c, B_{j} = b \right] \\
    & \quad=\frac{\theta}{J} \binom{c}{l} \binom{b}{a}  \frac{\phi^{(c-l+b-a)}((n+1) \lambda) }{\phi^{(c+b)}((n+1) \lambda)} \kappa(l+a, (n+1) \lambda).
  \end{split}
  \end{align}
\end{prp}

See Appendix~\ref{app:proof_poisson} for a proof of Proposition~\ref{prop:poisson_traits}. 
We further specialize Proposition~\ref{prop:poisson_traits} by assuming that $\mutilde$ is a Gamma CRM, i.e., $\rho(s) = s^{-1} \exp\{-s\}$. This simplifies \eqref{post_poiss} to:
\begin{align}\label{eq:pois_gamma}
\begin{split}
&\prob\left[f_{Y_{n+1,r}} = l, X_{n+1}(Y_{n+1,r}) = a \mid h(Y_{n+1,r}) = j,  C_{j} = c, B_{j} = b \right] \\
  & \quad=\frac{\theta}{J}\binom{c}{l} \binom{b}{a} (l + a - 1)! \frac{\Gamma(\theta/J + c + b - l - a)}{\Gamma(\theta/J + c + b)}.
\end{split}
\end{align}
See Appendix~\ref{app:poisson} for the proof of~\eqref{eq:pois_gamma}.  As a generalization of the Gamma CRM, we consider the generalized Gamma CRM, i.e., $\rho(s) = \alpha \Gamma(1-\alpha)^{-1} s^{- \alpha - 1} e^{-\tau s}$ for $\alpha\in[0,1)$ and $\tau>0$ \citep{Bri(99)}; see also \citet{Pit(03)} for details. 
This simplifies \eqref{post_poiss} to:
\begin{align}\label{eq:pois_gg}
\begin{split}
&\prob\left[f_{Y_{n+1,r}} = l, X_{n+1}(Y_{n+1,r}) = a \mid h(Y_{n+1,r}) = j,  C_{j} = c, B_{j} = b \right] \\
 & \quad= \frac{\theta}{J} \binom{c}{l} \binom{b}{a}\frac{\alpha (1 - \alpha)_{(l + a -1)}}{(\tau + (n+1)r)^{-\alpha + l + a}} \frac{
        \sum_{i=1}^{c-l+b-a} \left(\frac{\theta}{J}\right)^i \frac{\calC(c - l + b -a, i; \alpha)}{(\tau + (n+1)\lambda)^{c - l + b -a - \alpha i}}
    }{
    \sum_{i=1}^{c+b} \left(\frac{\theta}{J}\right)^i \frac{\calC(c + b , i; \alpha)}{(\tau + (n+1)\lambda)^{c + b - \alpha i}}
    }.
  \end{split}
\end{align}
See Appendix~\ref{app:poisson} for the proof of~\eqref{eq:pois_gg}. Both \eqref{eq:pois_gamma} and \eqref{eq:pois_gg} are easy to evaluate. 

As explained in Section~\ref{sec2}, applying \eqref{eq:pois_gamma} or \eqref{eq:pois_gg} necessitates estimating the unknown parameters $\lambda$ and the prior parameters from the sketch $\mathbf{C}_{J}$. Given the convolution properties of the Poisson distribution, the distribution of $\mathbf C_j$ under the model \eqref{eq:trait_model_hash}, where $G_A(J_k)$ is a Poisson$(\lambda J_k)$, corresponds to the distribution of the following hierarchical model:
\begin{align} \label{eq:cj_model}
\begin{split}
	C_j \mid T^\prime_j & \quad \simind \quad \mbox{Poi}(n \lambda T^\prime_j), \\
	T^\prime_j & \quad \simind \quad f_{D_j} (\theta, \rho, G_0),
      \end{split}
\end{align}
where $f_{D_j}(\theta, \rho, G_0)$ is the distribution of $\mutilde(D_j)$. If $\mutilde$ is modeled according to a Gamma CRM, then $f_{D_j}$ follows a Gamma distribution with parameter $(\theta/J, 1)$. Consequently,
\begin{equation}\label{like_ibp}
	\prob[\mathbf  C_J = \mathbf c_J] = \frac{(n \lambda)^n}{(1 + n\lambda)^{\theta + n}} \prod_{j=1}^J \frac{\left(\frac{\theta}{J}\right)_{(c_j)}}{c_j!},
\end{equation}
which allows estimating $\theta$ and $\lambda$ by maximizing the marginal likelihood \eqref{like_ibp}. If $\mutilde$ is a generalized Gamma CRM, the distribution of $f_{D_j}$ is not conjugate to the Poisson distribution. This prevents obtaining a closed-form distribution for the sketch $\mathbf{C}_{J}$. However, we can apply a likelihood-free method similar to that discussed for the PYP prior in Section~\ref{sec2}.

\subsubsection{Comparisons under the Gamma and generalized Gamma process priors}

\begin{figure}[!htb]
  \centering
  \includegraphics[width=22cm,height=14cm,keepaspectratio]{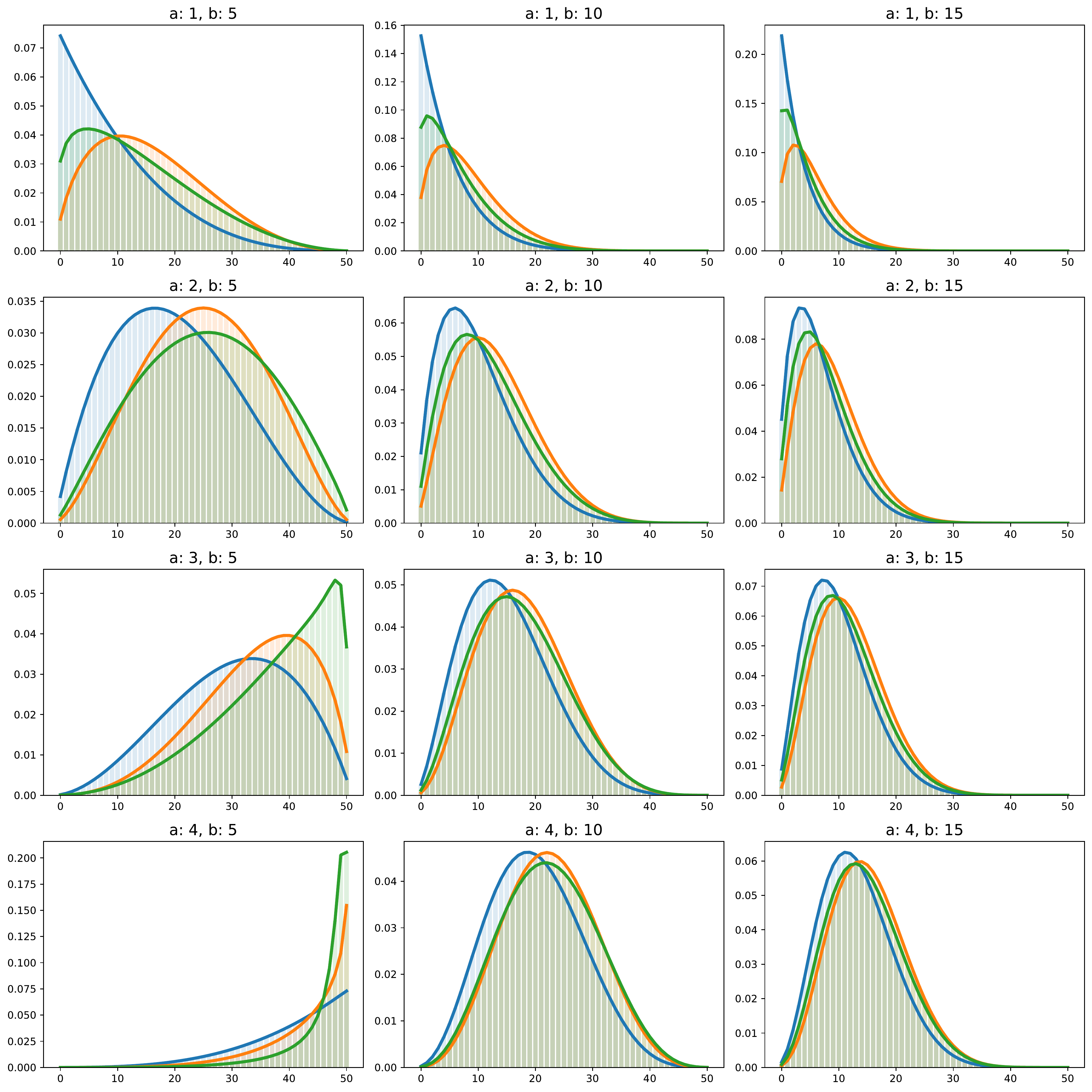}
  \caption{\small{Posterior distribution of the empirical frequency level $f_{Y_{n+1,r}}$, assuming the Poisson setting with a Gamma CRM prior (blue line) and with a generalized Gamma CRM prior (orange and green lines). In all the panels, $a$ and $b$ vary, while we fix $c=50$, $m=1000$, $J=50$, $\theta=0.3$, $\tau = 1$, $\lambda=1$, and $\alpha = 0.25, 0.75$ for the orange and green lines respectively.}}
  \label{fig:poi_post}
\end{figure}

We compare the posterior distributions of $f_{Y_{n+1, r}}$ under the Gamma CRM and generalized Gamma CRM priors, with $G_A(\cdot \mid J_k)$ modeled as a Poisson distribution with parameter $\lambda J_k$, for fixed $\lambda>0$. Notably, when $a = b = 1$, the posterior distribution under the Gamma CRM prior coincides with the posterior distribution obtained under the DP prior, which is displayed in \eqref{post_dp}. In other cases, the posterior distribution of $f_{Y_{n+1, r}}$ incorporates additional information from $X_{n+1}(Y_r) = a$ and $X_{n+1}(D_{h(Y_{n+1,r})}) = b$, providing details on the relative frequency of trait $Y_{n+1,r}$ within its corresponding bucket. Figure \ref{fig:poi_post} contrasts the posterior distributions under Gamma CRM and generalized Gamma CRM priors, illustrating significant differences influenced by the specific prior settings, especially for small values of $X_{n+1}(Y_{n+1, r})$. Moreover, for a given $X_{n+1}(D_{h(Y_{n+1,r})}) = b$, values of $X_{n+1}(Y_r)$ close to $b$ push the posterior distribution to large values, reflecting the prevalence of the trait $Y_{n+1,r}$ within that bucket.  Conversely, values of $X_{n+1}(Y_r)$ considerably smaller than $b$ indicate the rarity of the trait within the bucket, driving the posterior of $f_{Y_{n+1, r}}$ towards lower values.

\subsection{Results under $G_{A}$ Bernoulli}

Let $G_A(\cdot \mid J_k)$ be the Bernoulli distribution with parameter $J_k$. Under this assumption, the levels of associations of traits simply refer to their presence or absence. For this reason, traits are typically referred to as features. The conditional distribution of the random variable $S_n = \sum_{k \geq 1} \sum_{i=1}^k A^\prime_{i, k}$ in \eqref{eq:post_trait}, given $\mutilde$, is the distribution of the sum of independent Bernoulli random variables with parameters $J^\prime_1, \ldots, J^\prime_1, J^\prime_2, \ldots, J^\prime_2, \ldots$, where each $J^\prime_k$ appears exactly $n$ times; this is the Poisson-Binomial distribution \citep{Che(97)}. Similarly, $Z = \sum_{k \geq 1} A^\prime_{n+1, k}$ follows a Poisson-Binomial distribution with parameter $J^\prime_1, J^\prime_2, \ldots$. While the Poisson-Binomial distribution has a complicated expression, one may combine Le Cam's theorem \citep{LeCam(60)} with the results in Section~\ref{sec:poisson_traits} to approximate \eqref{eq:post_trait} in the Bernoulli setting. For notation's sake, considering $S_n$ and $Z$ as previously defined, we introduce $\tilde S_n$ and $\tilde Z$ such that, given $T^\prime = \sum_{k \geq 1} J^\prime_k$, $\tilde S_n \mid T^\prime \sim \mbox{Poi}(nT^\prime)$ and $\tilde Z \mid T^\prime \sim \mbox{Poi}(T^\prime)$. We first observe that the posterior distribution of $f_{Y_{n+1, r}}$ in \eqref{eq:post_trait} is proportional to
\begin{align*}
 & \prob\left[f_{Y_{n+1,r}} = l, X_{n+1}(Y_{n+1,r}) = 1 \mid h(Y_{n+1,r}) = j,  C_{j} = c, B_{j} = b \right]\\
  &\qquad\propto\prob\left[f_{Y_{n+1,r}} = l, X_{n+1}(Y_{n+1,r}), S_n =  c-l, Z = b-1\right].
\end{align*}
In the next theorem, we provide an approximation of the distribution of $f_{Y_{n+1,r}}$ by replacing $S_n$ and $Z$ with $\tilde S_n$ and $\tilde Z$ respectively. We also estimate the approximation error.

\begin{thm}\label{prop:bern_traits}
Let $\mathbf{C}_{J}$ be a sketch of $\mathbf{X}_{n}$ under the model \eqref{eq:trait_model_hash} with $G_A(J_k)$ being the Bernoulli distribution with parameter $J_k$, and let $X_{n+1}$ be an additional (unobservable) random sample. Furthermore, let $\psi$ be the function defined in \eqref{eq_cond}, and let $\phi^{(n)}(u) = (-1)^{n} \frac{\dd}{\dd u} e^{-\theta/J \psi(u)}$ and  $\kappa(u, n) = \int_{\R_+} e^{-un} s^n \rho(s) \dd s$. Then, the posterior distribution of $f_{Y_{n+1, r}}$ can be approximated by the distribution of the random variable $\tilde f_{Y_{n+1, r}}$ such that, for $l=0, \ldots, c$:
  \begin{align}\label{eq:post_features}
    \begin{split}
    & \prob\left[\tilde f_{Y_{n+1,r}} = l, X_{n+1}(Y_{n+1,r}) = 1\right]\\
     & \quad= \prob\left[f_{Y_{n+1,r}} = l, X_{n+1}(Y_{n+1,r}) = 1 \mid \tilde S_n =  c-l, \tilde Z = b-1\right] \\
     & \quad\propto
    (c - l + 1)_{(l)} \binom{n}{l} \phi^{(c+b - l - 1)}(n+1) \int (s^{l+1} - s^{n+1}) \rho(s) \dd s.
  \end{split}
  \end{align}
  Moreover, the total variation distance between the random vectors $[f_{Y_{n+1, r}}, \allowbreak X_{n+1}(Y_{n+1, r}),S_n, Z]$ and $[f_{Y_{n+1, r}}, X_{n+1}(Y_{n+1, r}), \tilde S_n, \tilde Z]$
 is upper bounded by $(2\theta/J) \int_{\R_+} e^{-\psi(u)} \kappa(u, 2) \dd u$.
\end{thm}


\section{Numerical illustrations}\label{sec:numerics}

\subsection{Frequency recovery}

\subsubsection{Limitations of the DP with heavy-tailed data distributions}

To illustrate the limitations of the DP prior, we conducted two simulations with $n = 500,000$ data points, simulated either from a DP with parameters $\theta = 5, 10, 20, 100$ or from a Zipf distribution with tail parameters $c = 1.18, 1.54, 1.82, 2.22$. The Zipf distribution, applicable to infinitely many items, assigns the probability of the $k$-th item as $k^{-c} / \zeta(c)$,
where $\zeta(\cdot)$ denotes the Riemann's zeta function, defined by $\zeta(c) = \sum_{k \geq 1} k^{-c}$.

The DP is known to produce distributions with geometric (light) tails \citep{Teh(10)}, which may not adequately capture data with heavier tails behaviours. In contrast, the Zipf distribution exhibits power-law tails, where the parameter $c$ directly influences the decay rate of the tail: lower values of $c$ indicate heavier tails, representing larger fractions of low-frequency items.

We assume model \eqref{eq:exchangeable_model_hash} in combination with DP prior for $P$, whose total mass parameter is estimated from the sketch as discussed in Section~\ref{sec:dp_freq}.
We let $J$ vary between $100$ and $5,000$.
As evaluation metric, we follow \cite{Dol(23)} and consider the mean absolute error (MAE) stratified by the true frequency of the tokens. That is, for all distinct symbols $s$ appearing in the original data set, we compute
\begin{displaymath}
    \text{MAE}_m = \frac{1}{\sum_{s \in \mathbb S} I(f_s \in (l_m, u_m])} \sum_{s \in \mathbb S} |f_s - \hat f_s | \, I\left(f_s \in (l_m, u_m]\right),
\end{displaymath}
where $f_s = \sum_{i=1}^n I(X_i = s)$ is empirical frequency, $\hat f_s$ its  estimate, and $(l_m, u_m]_{m \geq 1}$ are non-overlapping frequency bins.

\begin{figure}[!htb]
  \centering
  \includegraphics[width=12cm,height=6cm,keepaspectratio]{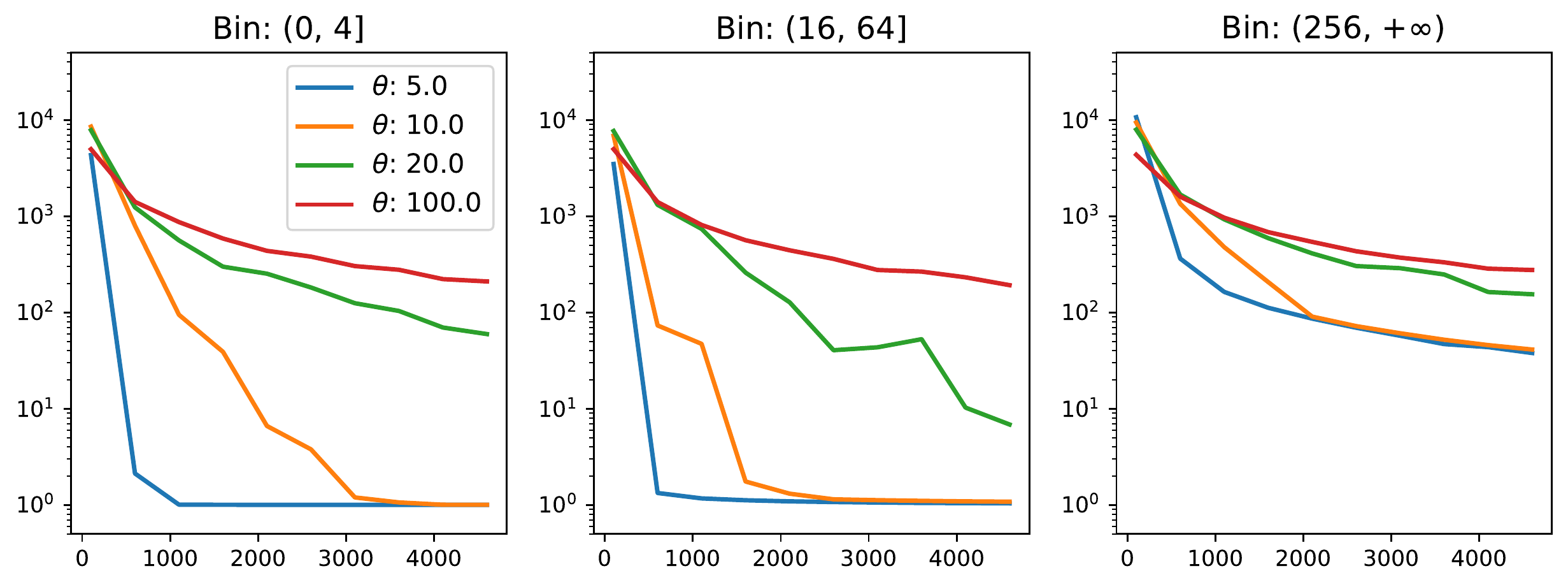}
  \includegraphics[width=12cm,height=6cm,keepaspectratio]{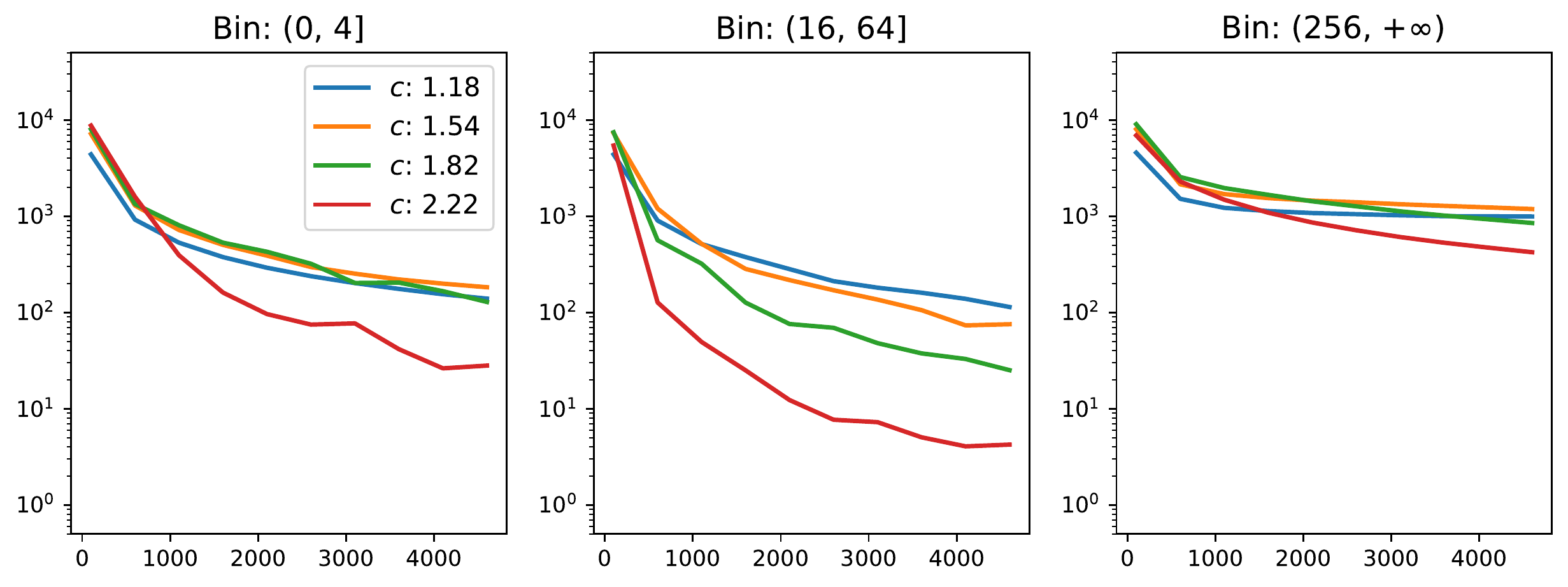}
  \caption{\small{MAEs as a function of $J$ for the frequency recovery problem under a DP prior. Top: data simulated from a DP with parameter $\theta$. Bottom:  data from a Zipf distribution with parameter $c$.}}
  \label{fig:dp_simu}
\end{figure}

The top row of Figure~\ref{fig:dp_simu} shows the MAEs when data are generated from a DP with parameter $\theta \in \{5, 10, 20, 100\}$. In this case, the MAEs for low and mid-frequency tokens decrease rapidly with $J$, especially if $\theta$ is small. 
This is expected since lower values of $\theta$ correspond to fewer distinct symbols in the sample, reducing the likelihood of hash collisions.

The results for Zipf-distributed data are displayed in the bottom row of Figure~\ref{fig:dp_simu}.
In this case, the MAEs are considerably larger, especially when $c$ is small. This is unsurprising since, as discussed above, smaller values of $c$ correspond to a greater number of distinct symbols in the sample, increasing the likelihood of hash collisions.
Moreover, the DP prior is not suited to modeling heavy-tailed data, which is clearly reflected in the very high MAEs associated with the low-frequency tokens.

\subsubsection{Frequency recovery with the approximate posterior under the PYP}

\begin{figure}[!htb]
  \centering
  \includegraphics[width=12cm,height=6cm,keepaspectratio]{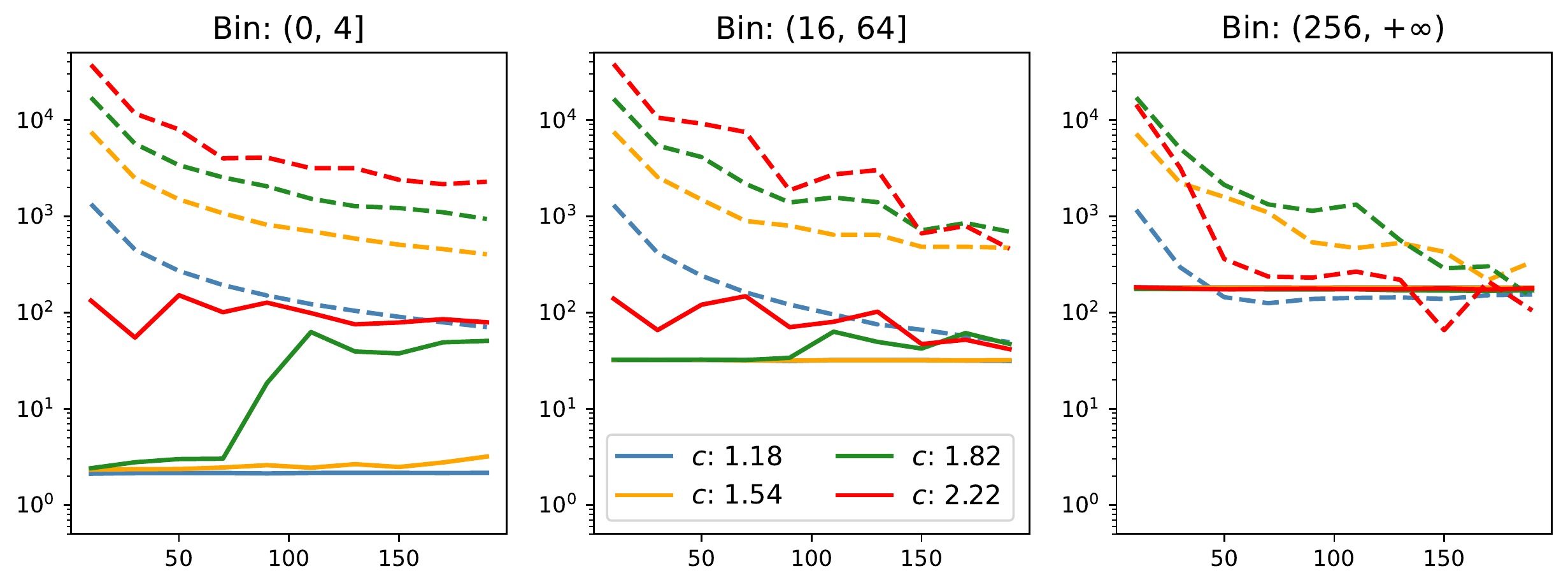}
  \caption{\small{MAEs as a function of $J$ for the frequency recovery problem of Zipf($c$) data under a PYP prior. Solid lines correspond to prior's parameters $(\alpha,\theta)$ estimated with our method; dashed lines correspond to estimates of $(\alpha,\theta)$ obtained as in \cite{Dol(23)}.}}
  \label{fig:py_asyn_simu}
\end{figure}

We consider now the same synthetic datasets of \cite{Dol(23)}, which consist of $n=500,000$ samples from the Zipf  distribution with parameter $c=1.18, 1.54, 1.82, 2.22$. We make use of \eqref{approx_est} to estimate the frequencies via the posterior expectation, and then we compare estimating the parameters via the likelihood-free approach in \cite{Dol(23)} (specifically, we refer to their estimated values) and our method based on minimizing the recovery error on the first 10,000 samples. The results are displayed in Figure~\ref{fig:py_asyn_simu}. 

Note that, to satisfy the asymptotic regime of Theorem~\ref{char_pyp}, we set $J$ to smaller values compared to the previous simulation. It is clear that the proposed estimation method leads to a significant enhancement for the frequency recovery problem. Additionally, employing a PYP prior yields markedly lower MAEs compared to a DP prior, particularly for low frequency tokens. Interestingly, while larger values of the parameter $c$ correspond to better performance under the DP prior, the opposite holds true under the PYP.

This can be explained as follows: when the tail parameter of the Zipf distribution, $c$, is large, the data generating process exhibits lighter tails and fewer data points contribute to most of the total mass, commonly referred to as ``heavy hitters''. 
When predicting the frequency of these heavy hitters, the asymptotic regime under which \eqref{approx_est} is derived may not hold, as $c_j$ may be comparable to $n$ for some $j$ while $c_k \approx 0$ for other $k \neq j$.

A similar scenario might arise if $J$ is large. From \eqref{eq:post-mean}, it is clear that our estimator applies linear shrinkage to $c_j$ by a term inversely proportional to $J$. Therefore, if $J$ is large and $c_j$ is not sufficiently large, the estimator over-shrinks the frequency to zero. 
This becomes evident in our example when $c = 1.82$: in this case, the MAEs associated with low and mid-frequency tokens show an increase as $J$ increases, which can be attributed to the invalidity of the assumptions underlying Theorem~\ref{char_pyp} for larger values of $J$.

\subsubsection{Application to the Gutenberg corpus}\label{sec:gut_freq}

\begin{figure}[!htb]
\centering
\includegraphics[width=.28\linewidth,valign=c]{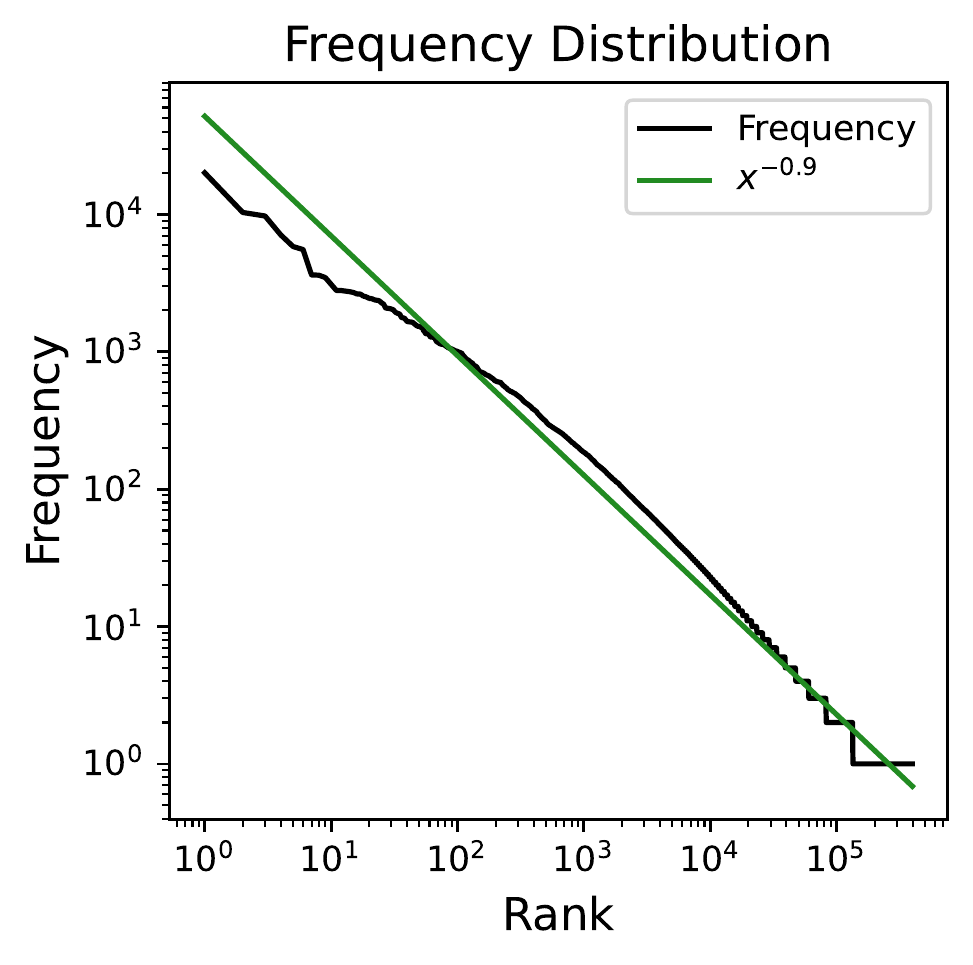} %
\includegraphics[width=0.6\textwidth,valign=c]{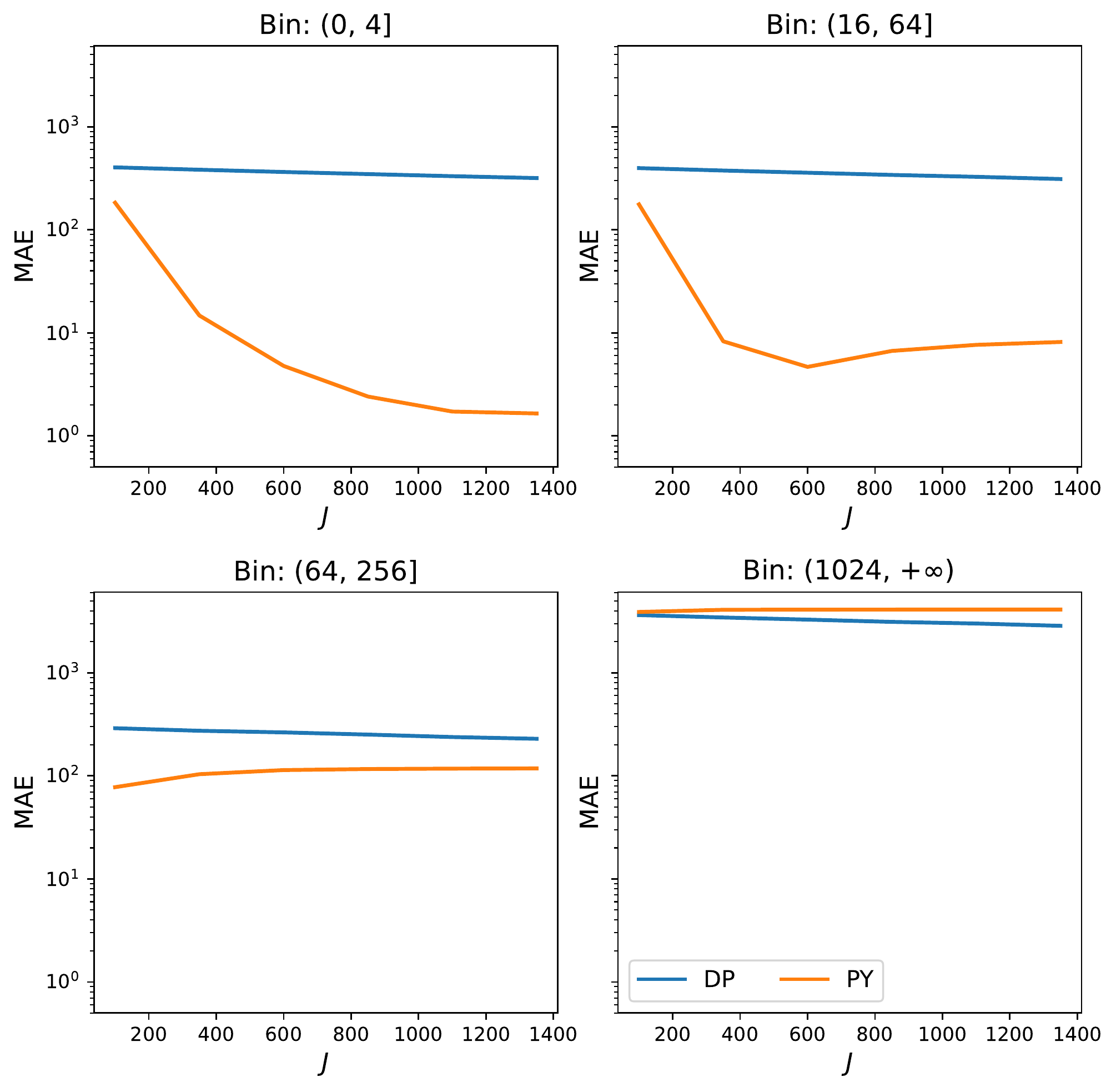}
\caption{\small{Left: frequency distribution for the Gutenberg corupus' bigrams. Right: MAEs as a function of $J$ for the frequency recovery problem.}}
\label{fig:bigrams_freq}
\end{figure}

The second data set comprises 18 open-domain classic pieces of English literature from the Gutenberg Corpus~\citep{Gutenberg}.
These data are pre-processed with the same approach of \citet{Ses(22)}: 
punctuation and unusual words are removed, retaining only words found in an English dictionary of size 25,487. Subsequently, 1,700,000 consecutive word pairs, or \textit{2-grams}, are extracted. These 2-grams are then sketched using a random hash function, as usual.
As shown in Figure~\ref{fig:bigrams_freq} (left), these data exhibit a clear power law distribution characterized by numerous bigrams with low frequency.

We compare the frequency estimates obtained under the DP and PYP priors, with the latter utilizing the large-$n$ approximation derived in Theorem~\ref{char_pyp}. 
Given the large sample size, it is reasonable to expect that such an approximation should be quite accurate  in this case. 
Figure~\ref{fig:bigrams_freq} (right) displays the results, averaged over 20 independent hash functions. As anticipated, the PYP outperforms the DP significantly for low and mid-frequency tokens, while their performance is essentially equivalent for very high frequency ones (note that the number of tokens appearing more than 1024 times is less than 100).

\subsection{Cardinality recovery}

\subsubsection{Experiments with synthetic data}

We investigate here the empirical performance on simulated data of the BNP cardinality estimator derived in Section~\ref{sec21}, for the special cases of the DP prior of the general PYP prior.
Synthetic data are generated from the BNP model in~\eqref{eq:exchangeable_model_hash} using different values of the PYP parameters $(\alpha,\theta)$, and then they are sketched using a hash function of width 128.
Our BNP estimates are compared to the ground truth, which is available in these experiments because we know the prior parameters of the data-generating model and have access to the non-sketched data.
All experiments are repeated 20 times and the results averaged, utilizing independent data sets and independent hash functions.

\begin{figure}[!htb]
\centering
\includegraphics[width=0.8\linewidth]{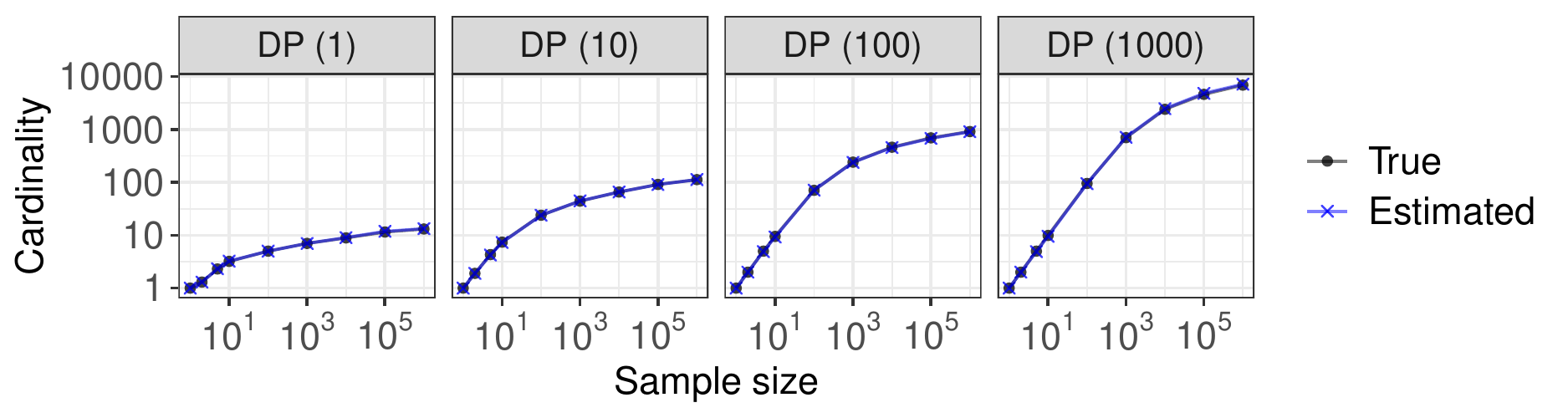}
\caption{\small{True and estimated cardinality in synthetic data from a DP prior model with different parameters, as a function of the sample size.
The true and estimated cardinalities are almost indistinguishable.}}
\label{fig:exp-dp-dp-distinct}
\end{figure}

Figure~\ref{fig:exp-dp-dp-distinct} compares the true and estimated cardinality as a function of the sample size $n$, separately for data generated from DP prior models with $\alpha=0$ and different values of $\theta$ (Figure~\ref{fig:exp-dp-dp-distinct}a), and for PYP models with $\theta=100$ and different values of $\alpha$ (Figure~\ref{fig:exp-dp-dp-distinct}b).
To make the computations practical, all estimates are computed under the (possibly mis-specified) assumption that $\alpha=0$, and estimating $\theta$ empirically via maximum marginal likelihood.
As predicted by the theory, the results confirm our estimates are accurate when the DP prior is well-specified, while otherwise they tend to underestimate the number of distinct species, especially if $n$ is large.
Similar results are shown by Figure~\ref*{fig:exp-dp-zipf-distinct} in Appendix~\ref{app:numerical}, which reports on analogous experiments based on synthetic data generated from a Zipf distribution.

\subsubsection{Experiments with real data}

In addition to the Gutenberg corpus' bigrams dataset discussed in Section~\ref{sec:gut_freq}, we consider here two additional datasets.
The first one was made publicly available by the National Center for Biotechnology Information~\citep{hatcher2017virus} and contains 43,196 sequences of approximately 30,000 nucleotides each, from SARS-CoV-2 viruses.
For each sequence, we extract a list of all contiguous DNA sub-sequences of length $16$ (i.e., {\em 16-mers}), and then we sketch the resulting data set with the random hash function.
The last data set is discussed in \citet{rojas2018personalized} and contains a list of 3,577,296 IP addresses, which we sketch directly without pre-processing; these data were made publicly available through the Kaggle machine-learning competition website.
For all data sets, in these experiments we separately consider sketches obtained with two distinct hash functions, one with a width of 128 and another with a width of 4096.

\begin{figure}[!htb]
\centering
\includegraphics[width=0.8\linewidth]{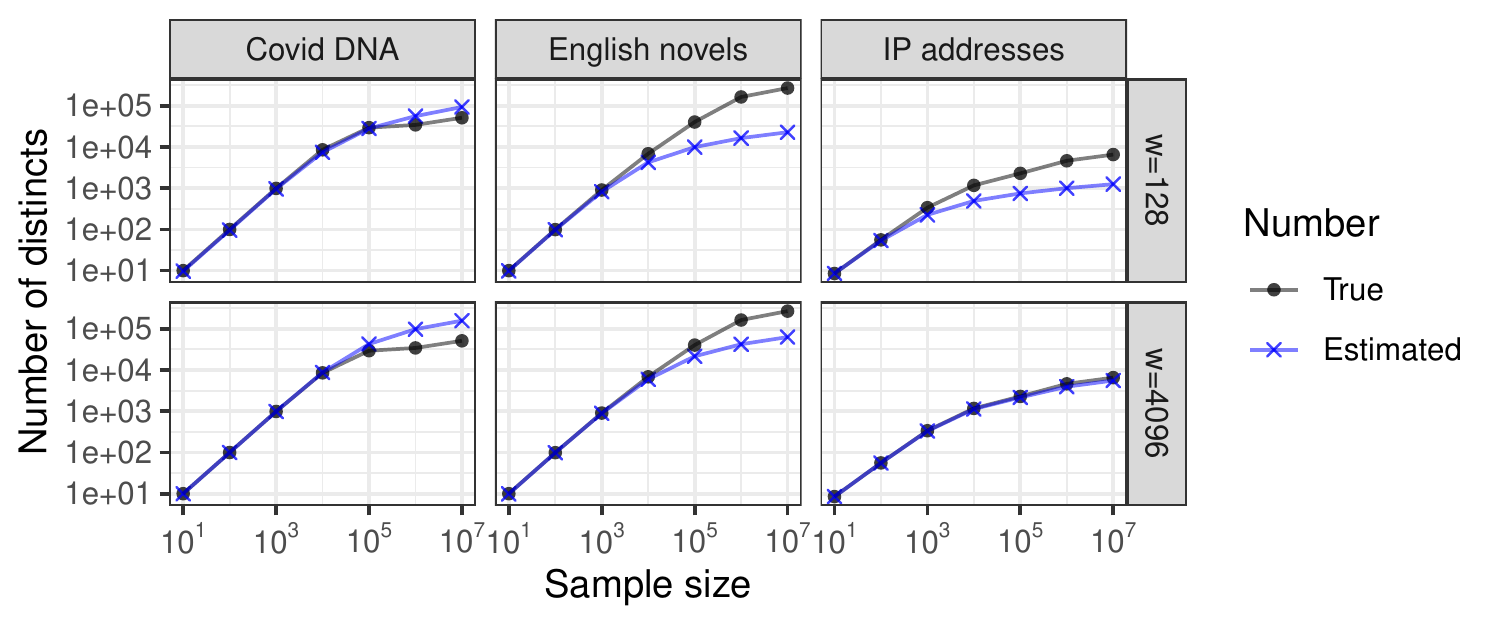}
\caption{\small{True and estimated cardinality in real data sets sketched with random hash functions of different width, as a function of the sample size.}}
\label{fig:exp-dp-data-distinct}
\end{figure}

Figure~\ref{fig:exp-dp-data-distinct} compares the true and estimated cardinality for random subsets of the three aforementioned data sets, as a function of the sample size and for two different values of the hash width.
Our BNP cardinality estimators are computed assuming $\alpha=0$, and estimating $\theta$ empirically via maximum marginal likelihood. 
All results are averaged over 20 independent experiments with different hash functions.
The results show the cardinality estimates obtained under the DP prior are relatively accurate for the DNA data set, which does not exhibit power-law tail behaviour \citep{Ses(22)}, but tend to underestimate the true missing mass in the other cases, especially if the hash function width is small.

\section{Discussion}\label{sec4}


In the ``species'' setting, while \citet{Cai(18)} and \citet{Dol(21),Dol(23)} employed BNPs to create learning-augmented versions of the CMS using DP and PYP priors, our approach diverges from the CMS framework as we condition on the information contained in the entire sketch. 
Moreover, our approach encompasses the broader class of PK priors and also facilitates cardinality recovery using the same data sketch. 
PK priors have been widely applied in BNPs, assuming access to true data, largely for their mathematical tractability, which leads to posterior inferences that are straightforward, computationally efficient, and scalable to large datasets \citep{Lij(10)}. However, we have shown these advantages of PK priors do not necessarily extend to BNP inferences from sketched data, presenting unique challenges and constraints in selecting appropriate prior distributions.

This challenge creates opportunities for: (i) further exploration of both exact and approximate algorithms for the efficient numerical evaluation of the posterior under the PYP prior; (ii) continued investigation into large-sample approximations of the posterior under the PYP prior, potentially with reliable error bounds; and (iii) the consideration of alternative priors aimed at simplifying posterior distributions while providing more flexible tail behaviors than those offered by the DP.

The ``traits'' setting of the frequency recovery problem introduces another novel aspect that goes beyond the works of \citet{Cai(18)} and \citet{Dol(21),Dol(23)}, highlighting the adaptability of the BNP framework to diverse data structures. 
Unlike the ``species'' setting, the ``traits'' setting demonstrates greater flexibility in the choice of prior distributions, enabling the tractable evaluation of posterior distributions for large sample sizes, especially under the assumption of a Poisson distribution for the levels of trait associations. 
Such a desirable property stems from the Poisson process formulation of CRMs, which determines a posterior distribution that depends on the sketch $\mathbf{C}_{J}$ only through $C_{h(X_{n+1})}$. 

In general, for both the ``species" and ``traits" settings, we argue that the BNP approach is also flexible with respect to the object of interest, due to the use of the posterior distribution as the main tool to obtain estimates. In the future, our BNP approach may also be extended to recover other quantities beyond those studied in this paper. 
Of notable interest in the CMS literature is the problem of recovering the (cumulative) empirical frequency of $s\geq1$ new data points \citep[Chapter 3]{Cor(20)}. We refer to \citet{Dol(23)} for a discussion in the context of learning-augmented versions of the CMS.


\section*{Acknowledgements}
M.~B.~and S.~F.~were funded by the European Research Council under the Horizon 2020 research and innovation programme, grant 817257. M.~B.~and S.~F.~also gratefully acknowledge support from the Italian Ministry of Education, University and Research, ``Dipartimenti di Eccellenza" grant 2023-2027.


\appendix


\section{Proofs}

\subsection{Proof of \Cref{main_pk}}\label{app:proof_main}

By the definition of conditional probability, we have that
\begin{equation}\label{eq:sketch_post_fract}
  \prob\left[ f_{X_{n+1}} = l \mid \bm C_{J} = \bm c, h(X_{n+1}) = j\right] = \frac{\prob\left[ f_{X_{n+1}} = l, \bm C_{J} = \bm c, h(X_{n+1}) = j\right]}{\prob\left[ \bm C_{J} = \bm c, h(X_{n+1}) = j\right]}.
\end{equation}
Let $g(t) = Z_T t^{-\gamma} e^{-\beta t}$, where $Z_T$ is a normalizing constant.

\paragraph{Denominator.}
We have
\begin{align}
  & \prob[\bm C_{J} = \bm c, h(X_{n+1}) = j] \nonumber \\
  & \qquad = \binom{n}{c_1, \ldots, c_J} \prob\Big[X_1 \in D_j, \ldots, X_{c_j} \in D_j, X_{n+1} \in D_j,  \nonumber \\
  & \hspace{4cm} X_{i} \in D_1, i \in [c_j: c_j + c_1], \ldots, X_{i} \in D_j, i \in [\sum_{k=1}^{J-1} c_k, n] \Big] \nonumber \\
  & \qquad = \binom{n}{c_1, \ldots, c_J} \E\left[ P(D_j)^{c_j + 1} \prod_{k \neq j} P(D_k)^{c_k} \right]. \label{eq:denom_probs}
\end{align}
From \eqref{eq:denom_probs}, writing $P(\cdot) := \mutilde(\cdot) / \mutilde(\Sd)$ and the expression of $g(t) = g(\mutilde(\Sd))$
\begin{align*}
   & \prob[\bm C_{J} = \bm c, h(X_{n+1}) = j] = \binom{n}{c_1, \ldots, c_j} \E\left [P(D_j)^{c_j + 1} \prod_{k \neq j} P(D_k)^{c_k}\right] \\
   & \qquad = Z_T  \binom{n}{c_1, \ldots, c_j} \E\left[\mutilde(\Sd)^{-n  - 1- \gamma} e^{-\beta \mutilde(\Sd)}  \mutilde(D_j)^{c_j + 1} \prod_{k \neq j} \mutilde(D_k)^{c_k}\right] \\
   & \qquad = Z_T  \binom{n}{c_1, \ldots, c_j} \int_{\R_+ } \frac{u^{n+\gamma}}{\Gamma(n + \gamma + 1)} \E\left[e^{-(u + \beta) \mutilde(D_j)} \mutilde(D_j)^{c_j + 1} \right] \\
   & \qquad \qquad \qquad \times \prod_{k \neq j} \E\left[e^{-(u + \beta) \mutilde(D_k)} \mutilde(D_k)^{c_k} \right] \\
   & \qquad = Z_T \binom{n}{c_1, \ldots, c_j} \int_{\R_+ } \frac{u^{n+\gamma}}{\Gamma(n + \gamma + 1)} (-1)^{c+1} \frac{\dd^{c_j+1}}{\dd z^{c_j+1}} e^{- \theta \psi(z)/J}|_{(u + \beta)}   \\
   & \hspace{4cm}\times \prod_{k\neq j} (-1)^{c_k}\frac{\dd^{c_k}}{\dd z^{c_k}} e^{- \theta \psi(z)/J}|_{(u + \beta)} \dd u,
\end{align*}
where the second equality follows from the definition of PK model and the third from the Gamma identity $\mutilde(\Sd) ^{-n -1 - \gamma} = \int_{\R_+} u^{n+\gamma} / \Gamma(n + \gamma + 1) e^{-u \mutilde(\Sd)} \dd x$, an application of Fubini's theorem, and the independence property of CRMs. The last equality follows from the properties of the exponential function and the definition of Laplace exponent of the CRM.

\paragraph{Numerator.}
Let $B_{\omega}$ denote a ball of radius $\varepsilon$ around $\omega \in \Sd$.
We have
\begin{align*}
    &\prob\left[ \sum_{i=1}^n \indicator_{(X_{n+1})}(X_i) = l, \bm C_{J} = \bm c, h(X_{n+1}) = j\right] = \\
    & \qquad \lim_{\varepsilon \rightarrow 0} \int_{D_j} \prob\left[\sum_{i=1}^n \indicator_{B_{\omega^*}}(X_i) = l, \bm C_{J} = \bm c, X_{n+1} \in B_{\omega^*} \right] \dd \omega^*
\end{align*}
We consider the integrand,
\begin{align}
    & \prob\left[\sum_{i=1}^n \indicator_{B_{\omega^*}}(X_i) = l, \sum_{i=1}^n \indicator_{h(\omega^*)}(h(X_i)) = c, X_{n+1} \in B_{\omega^*} \right] \nonumber \\
    & \qquad = \binom{n}{l}\binom{n-l}{c_1, \ldots, c_j - l, \ldots c_J} \prob\Big[X_1 \in B_{\omega^*} \ldots X_l \in B_{\omega^*}, X_{n+1} \in B_{\omega^*}, \nonumber \\
    & \hspace{5cm} X_{l+1} \in D_j \setminus B_{\omega^*} \ldots, X_{c_j} \in D_j \setminus B_{\omega^*}, \nonumber \\
    &  \hspace{5cm} X_{i} \in D_1, i \in [c_j: c_j + c_1], \ldots, X_{i} \in D_j, i \in [\sum_{l=1}^{J-1} c_l, n] \Big] \nonumber \\
     & \qquad = \binom{n}{l, c_1, \ldots, c_j - l, \ldots c_J} \E\left[P(B_{\omega^*})^{l+1} P(D_j \setminus B_{\omega^*})^{c_j-l} \prod_{k \neq j} P(D_k)^{c_k} \right]. \label{eq:num_probs}
\end{align}
To evaluate the expected value, we proceed as in the case of the denominator and write
\begin{align*}
    &  \E\left[P(B_{\omega^*})^{l+1} P(D_j \setminus B_{\omega^*})^{c_j-l} \prod_{k \neq j} P(D_k)^{c_k} \right] \\
    & \qquad  = Z_T \int_{\R_+ } \frac{u^{n + \gamma}}{\Gamma(n + \gamma + 1)} \E \left[ e^{- (u + \beta) \mutilde(B_{\omega^*})} \mutilde(B_{\omega^*})^{l+1} \right]  \\
    & \qquad \qquad \times
    \E \left[ e^{- (u + \beta) \mutilde(D_j \setminus B_{\omega^*})}  \mutilde(D_j \setminus B_{\omega^*})^{c-l}  \right]
    \prod_{k \neq j} \E\left[e^{-(u + \beta) \mutilde(D_k)} \mutilde(D_k)^{c_k} \right]
    \dd u.
\end{align*}
The two latter expected values above can be computed as in the denominator case.
For the first expectation instead, letting $\varepsilon \rightarrow 0$ we have
\[
    \E \left[ e^{- (u + \beta) \mutilde(B_{\omega^*})} \mutilde(B_{\omega^*})^{l+1} \right] \rightarrow \kappa(u + \beta, l+1) \theta G_0(\dd \omega^*),
\]
where $\kappa(u, l) = \int_{\R_+} s^{l} e^{-us} \rho(s) \dd s$. This can be verified using, for instance, Lemma 1 in \cite{Cam(19)}.
Hence, the numerator equals
\begin{align*}
    & \prob\left[ \sum_{i=1}^n \indicator_{(X_{n+1})}(X_i) = l, \bm C_{J} = \bm c, h(X_{n+1}) = j\right] \\
    & \qquad \frac{1}{l!}\binom{n}{c_1, \ldots, c_j - l, \ldots c_J}\frac{Z_T}{\Gamma(n + \gamma + 1)} \\
    & \hspace{2cm} \times\int_{D_j} \int_{\R_+} u^{n + \gamma} (-1)^{c_j - l} \left(\frac{\dd^{c_j-l}}{\dd (u+\beta)^{c_j-l}} e^{- \theta / J \psi(u+\beta)}\right) \times \\
    & \hspace{3cm}\times \prod_{k\neq j} (-1)^{c_k}\frac{\dd^{c_k}}{\dd z^{c_k}} e^{- \theta \psi(z)/J}|_{(u + \beta)} \kappa(u + \beta, l+1) \theta \dd u G_0(\dd \omega^*),
\end{align*}
where we can further integrate with respect to $\dd \omega^*$ and observe $\int_{D_j} G_0(\dd \omega^*) = 1 / J$ thanks to the universality assumption on the hash function $h$.
Combining numerator and denominator, and using the definition of $\phi^{(n)}(u)$, yields the result.

\subsection{Proof of \Cref{char_dp}}\label{app:char_proof}

Recalling the definition of $\psi^{(n)}(u)$, the claim of the theorem entails that
\[
    \frac{\int_{\R_+} u^{n + \gamma} \frac{\dd^{c_j-l}}{\dd (z)^{c_j-l}} e^{- \theta / J \psi(z)}|_{u+\beta}
    \prod_{k\neq j} \frac{\dd^{c_k}}{\dd z^{c_k}} e^{-\theta \psi(z)/J}|_{(u + \beta)} \kappa(u + \beta, l+1) \dd u }{\int_{\R_+ } u^{n+\gamma} \frac{\dd^{c_j+1}}{\dd z^{c_j+1}} e^{- \theta \psi(z)/J}|_{(u + \beta)}  \prod_{k\neq j} \frac{\dd^{c_k}}{\dd z^{c_k}} e^{- \theta \psi(z)/J}|_{(u + \beta)} \dd u} = f(n, c_j, l),
\]
where $f$ is an unknown function which cannot depend on $c_k, k \neq j$.
Then
\begin{multline*}
    \int_{\R_+} u^{n+ \gamma} \frac{\dd^{c_j-l}}{\dd (z)^{c_j-l}} e^{- \theta / J \psi(z)}|_{u+\beta}
    \prod_{k\neq j} \frac{\dd^{c_k}}{\dd z^{c_k}} e^{-\theta \psi(z)/J}|_{(u + \beta)} \kappa(u + \beta, l+1) \dd u = \\
    f(n, c_j, l) \int_{\R_+ } u^{n+\gamma} \frac{\dd^{c_j+1}}{\dd z^{c_j+1}} e^{- \theta \psi(z)/J}|_{(u + \beta)}  \prod_{k\neq j} \frac{\dd^{c_k}}{\dd z^{c_k}} e^{- \theta \psi(z)/J}|_{(u + \beta)} \dd u,
\end{multline*}
or, equivalently,
\begin{multline*}
     \int_{\R_+} u^{n + \gamma} \prod_{k\neq j} \frac{\dd^{c_k}}{\dd z^{c_k}} e^{- \theta \psi(z)/J}|_{(u + \beta)} \times \\
     \left( \frac{\dd^{c_j-l}}{\dd (z)^{c_j-l}} e^{- \theta / J \psi(z)}|_{u+\beta} \kappa(u + \beta, l+1) - f(n, c_j, l) \frac{\dd^{c_j+1}}{\dd z^{c_j+1}} e^{- \theta \psi(z)/J}|_{(u + \beta)}  \right) \dd u = 0.
\end{multline*}
Since the above equality must hold true for all values of $c_k, k \neq j$, it follows that, $\forall u \geq 0$,
\begin{equation}\label{eq:diffeq}
     \frac{\dd^{c_j-l}}{\dd (z)^{c_j-l}} e^{- \theta / J \psi(z)}|_{u+\beta} \kappa(u + \beta, l+1) - f(n, c_j, l) \frac{\dd^{c_j+1}}{\dd z^{c_j+1}} e^{- \theta \psi(z)/J}|_{(u + \beta)} \equiv 0.
\end{equation}

The simplest nontrivial case is when $c_j = 1, l=0$ (indeed note that if $c_j = 0$, $l=0$ by definition and we get $f(n, c_j, l) = 1$).
Now, note that
\begin{align*}
    \frac{\dd}{\dd z} e^{- \theta / J \psi(z)} &= e^{- \theta / J \psi(z)} \left(- \frac{\theta}{J}\right) \frac{\dd}{\dd z} \psi(z), \\
    \frac{\dd^2}{\dd z^2} e^{- \theta / J \psi(z)} &= e^{- \theta / J \psi(z)} \left[\left(- \frac{\theta}{J}\right) \frac{\dd}{\dd z} \psi(z)\right]^2 + e^{- \theta / J \psi(z)} \left(- \frac{\theta}{J}\right) \frac{\dd^2}{\dd z^2} \psi(z).
\end{align*}
Plugging these into \eqref{eq:diffeq} we get
\begin{multline*}
    e^{- \theta / J \psi(u+\beta)} \Bigg\{  - \frac{\theta}{J} \frac{\dd}{\dd z} \psi(z)|_{ u+\beta} \kappa(u+\beta, 1) + \\
     - f(n, c_j, l) \left( \left[\left(- \frac{\theta}{J}\right) \frac{\dd}{\dd z} \psi(z)\right]^2_{u+\beta} \right) - \frac{\theta}{J}  \frac{\dd^2}{\dd z^2} \psi(z)|_{u+\beta}  \Bigg\} = 0.
\end{multline*}
Given the positivity of the exponential function, we can set the term in the curly brackets equal to zero.
Since $\dd / \dd z \, \psi(z) = - \kappa(z, 1)$, the term in the curly brackets above reduces to
\[
     \frac{\theta}{J} \kappa(u + \beta, 1)^2 - f(n, c_j, l) \left( \frac{\theta^2}{J^2} \kappa(u + \beta, 1)^2 + \frac{\theta}{J} \frac{\dd}{\dd z} \kappa(z, 1)|_{u + \beta}\right) = 0,
\]
which entails
\[
    \frac{\dd}{\dd z} \kappa(z, 1) = - \left(\frac{\theta}{J}  - \frac{1}{f(m, c_j, l)} \right) \kappa(z, 1)^2.
\]
Letting $w : = \left(\frac{\theta}{J}  - \frac{1}{f(n, c_j, l)}\right)^{-1}$, the differential equation above can be seen to have solution
\[
    \kappa(z, 1) = \frac{1}{c + \frac{z}{w}} = \frac{w}{\tau + z},
\]
where $c$ is an arbitrary constant and $\tau = c w$.

Let now $c_j = l = 0$ and recall that in this case $f \equiv 1$. Plugging these into \eqref{eq:diffeq} we get
\begin{align*}
    e^{- \theta / J \psi(u+\beta)} \kappa(u+\beta, 1) - \frac{\dd}{\dd z} e^{- \theta / J \psi(z)} &= 0,
\end{align*}
which leads to
\[
    e^{- \theta / J \psi(u+\beta)} \left(\frac{w}{\tau + z} + \frac{\theta}{J} \frac{\dd }{\dd z} \psi(z) \right)_{| u + \beta} = 0.
\]
Setting the term in the parentheses equal to zero, we obtain that
\begin{equation}\label{eq:lap_exp}
    \psi(z) = K \log(\tau + z).
\end{equation}
Hence, we have shown that if $\prob\left[ f_{X_{n+1}} = l \mid \bm C_{J} = \bm c, h(X_{n+1}) = j\right]$ does not depend on $c_k$, $k \neq j$, the CRM $\mutilde$ must have L\'evy exponent \eqref{eq:lap_exp}.

Let now $\mutilde^\prime$ be a CRM with L\'evy exponent $\psi$ as above. We first note that, without loss of generality, we can set $K=1$ since setting $K \neq 1$ the Laplace transform
\[
    \E e^{- z \mutilde^\prime(A)} = e^{- K \log(\tau + z) \theta G_0(A)} = (\tau + z)^{-K\theta G_0(A)} =: F^\prime(z)
\]
simply amounts to rescaling the total mass parameter.

We now show that $\tau$ must necessarily be equal to one.
Note that if $\mutilde_G$ is a Gamma process, then its L\'evy exponent is $\log(1 + z)$, which is \eqref{eq:lap_exp} but shifted of a term $(1 - \tau)$:
\begin{equation}\label{eq:laplace_shift}
     \E e^{- z \mutilde_G(A)} = (1 + z)^{-K\theta G_0(A)} =  F^\prime(z + (1 - \tau)).
\end{equation}
Let $f_A$ be the probability density function of $\mutilde_G(A)$, and $f^\prime_A$ be the probability density function of $\mutilde^\prime(A)$. By the properties of the Laplace transform, \eqref{eq:laplace_shift} is equivalent to
\[
    f_A(t) = e^{(\tau - 1)t} f^\prime_A(t), \quad \text{for all } t,
\]
which is clearly impossible if $\tau \neq 1$ since we must have that both $f_A$ and $f^\prime_A$ must integrate to one since they are probability density functions.
Hence, we have shown that if $\prob\left[ f_{X_{n+1}} = l \mid \bm C_{J} = \bm c, h(X_{n+1}) = j\right]$ does not depend on $c_k, k \neq j$, the underlying CRM must be a Gamma process.

We are left with two degrees of freedom: namely the parameters $\gamma$ and $\beta$ defining the change of measure in the Poisson-Kingman model.
However, note that if $\tilde \mu$ is a Gamma process, the resulting PK model with $g(t) \propto t^{-\gamma} e^{-\beta t}$ is still a DP. This can be checked, for instance, starting from Eq. (177) in \cite{Jam(02)}. This concludes the proof.

\subsection{Proof of Equation \eqref{post_dp} as a special case of Theorem~\ref{main_pk} }\label{app:dp_cor}

The DP is obtained by normalizing a Gamma process, whose L\'evy intensity is $\theta s^{-1} e^{-s} \dd s G_0(\dd x)$.
Hence,
\begin{align*}
    \psi(u) &= \int_{\R_+}(1 - e ^{-us}) \rho(s) \dd s = \int_{\R_+}(1 - e ^{-us}) s^{-1} e^{-s}\dd s \\
    &= \int_{\R_+} \int_{\R_+} (1 - e ^{-us}) e^{-ts} e^{-s} \dd t \, \dd s \\
    &= \int_{\R_+} \frac{u}{(t+1)(t+u+1)} \dd t = \log(1 + u),
\end{align*}
where the third equality follows from the identity $s^{-1} = \int_{\R_+} e^{-ts} \dd t$ and the fourth one by an application of Fubini's theorem.
Then it follows
\[
  (-1)^n \frac{\dd^u}{\dd u^n} e^{-\theta \psi(u)} = \frac{\Gamma(\theta + n)}{\Gamma(\theta)} (1 + u)^{- \theta -n},
\]
Moreover, $k(l, u) = \Gamma(l) / (u+1)^{l}$.
Hence, the integral at the numerator of \eqref{post_pk} can be evaluated as
\begin{align*}
   &\frac{\Gamma(\theta/J + c_j - l)}{\Gamma(\theta/J)} \prod_{k \neq j} \frac{\Gamma(\theta/J + c_k)}{\Gamma(\theta/J)} \Gamma(l+1)  \int_{\R_+} u^{n} + (1 + u)^{-\theta - n - 1} \\
   & \qquad = \frac{\Gamma(\theta/J + c_j - l)}{\Gamma(\theta/J)} \prod_{k \neq j} \frac{\Gamma(\theta/J + c_k)}{\Gamma(\theta/J)} \Gamma(l+1) \frac{\Gamma(n+1)}{\Gamma(\theta + n + 1)}.
\end{align*}
Similarly, the integral at the denominator of \eqref{post_pk}  equals
\begin{align*}
   &\frac{\Gamma(\theta/J + c_j + 1)}{\Gamma(\theta/J)} \prod_{k \neq j} \frac{\Gamma(\theta/J + c_k)}{\Gamma(\theta/J)}  \int_{\R_+} u^{n} + (1 + u)^{-\theta - n - 1} \\
   & \qquad = \frac{\Gamma(\theta/J + c_j + 1)}{\Gamma(\theta/J)} \prod_{k \neq j} \frac{\Gamma(\theta/J + c_k)}{\Gamma(\theta/J)} \frac{\Gamma(n+1)}{\Gamma(\theta + n + 1)}.
\end{align*}
Combining these expressions together, we have
\[
    \prob\left[ f_{X_{n+1}} = l \mid \bm C_{J} = \bm c, h(X_{n+1}) = j\right] = \frac{\theta}{J} \binom{c_j}{l} l! \frac{\Gamma(\theta/J + c_j - l)}{\Gamma(\theta / J + c_j + 1)}.
\]

\subsection{Proof of Equation \eqref{post_dp} from the finite dimensional laws of the DP prior}\label{app:dp_teo}

In this proof, we exploit only the original characterization of the DP in terms of its finite-dimensional distributions. That is, given $\theta >0$ and $G_0$ a probability measure on $(\Sd)$, $P$ is a DP with mean measure $\theta G_0$ if and only if, for any $n > 0$ and any $n$ measurable partition $A_1, \ldots, A_n$ of $\X$,
\begin{equation}\label{eq:dp_def}
    (P(A_1), \ldots, P(A_n)) \sim \mathrm{Dir}_n (\theta G_0(A_1), \ldots, \theta G_0(A_n)),
\end{equation}
where $\mathrm{Dir}_n$ denotes the $n-1$ dimensional Dirichlet distribution.

We argue as in \Cref{app:proof_main}.
To compute the denominator in \eqref{eq:sketch_post_fract}, from \eqref{eq:denom_probs}, \eqref{eq:dp_def}
\begin{align*}
   & \prob[\bm C_{J} = \bm c, h(X_{n+1}) = j] \nonumber \\
   & \qquad = \binom{n}{c_1, \ldots, c_J}  \frac{\Gamma(\theta)}{\Gamma(\theta + n + 1)} \frac{\Gamma(\theta / J + c_j + 1)}{\Gamma(\theta / J)} \prod_{k \neq j} \frac{\Gamma(\theta / J + c_k)}{\Gamma(\theta / J)}, \nonumber
\end{align*}
where we also exploited the e uniformity of the hash function which ensures that $G_0(D_\ell) = G_{0}(h^{-1}(\{\ell\})) = 1 / J$ for any $\ell = 1, \ldots, J$.

Similarly, to compute the numerator, consider \eqref{eq:num_probs}. Since $P$ is a DP,
\begin{align*}
    & \E\left[P(B_{\omega^*})^{l+1} P(D_j \setminus B_{\omega^*})^{c_j-l} \prod_{k \neq j} P(D_k)^{c_k} \right] \\
    & \qquad = \frac{1}{l!}\binom{n}{c_1, \ldots, c_j - l, \ldots c_J} \frac{\Gamma(\theta)}{\Gamma(\theta + n + 1)} \frac{\Gamma(\alpha G_0(B_{\omega^*}) + l + 1)}{\Gamma(\theta G_0(B_{\omega^*}))} \\
    & \qquad \qquad \times \frac{\Gamma(\theta G_0(D_j \setminus B_{\omega^*}) + c_j - l)}{\Gamma(\theta G_0(D_j \setminus B_{\omega^*}))} \prod_{k \neq j} \frac{\Gamma(\theta / J + c_k)}{\Gamma(\theta / J)}.
\end{align*}
We now let $\varepsilon \rightarrow 0$. First note that, of course
\[
    \frac{\Gamma(\theta G_0(D_j \setminus B_{\omega^*}) + c_j - l)}{\Gamma(\theta G_0(D_j \setminus B_{\omega^*}))} \rightarrow \frac{\Gamma(\theta G_0(D_j) + c_j - l)}{\Gamma(\theta G_0(D_j))} = \frac{\Gamma(\theta / J + c_j - l)}{\Gamma(\theta / J)}.
\]
To evaluate the limit of  $\Gamma(\theta G_0(B_{\omega^*}) + l + 1) / \Gamma(\theta G_0(B_{\omega^*}))$, we first unroll the numerator using the recurrence relation $\Gamma(z + 1) = z \Gamma(z)$ $l$ times so that
\[
    \frac{\Gamma(\theta G_0(B_{\omega^*}) + l + 1)}{\Gamma(\theta G_0(B_{\omega^*}))} = (\theta G_0(B_{\omega^*}) + l) \cdots (\theta G_0(B_{\omega^*})) = \theta \Gamma(l+1) G_0(B_{\omega^*}) + o(G_0(B_{\omega^*})).
\]
Letting now $\varepsilon \rightarrow 0$, we can ignore higher order infinitesimals and get that
\[
    \frac{\Gamma(\theta G_0(B_{\omega^*}) + l + 1)}{\Gamma(\theta G_0(B_{\omega^*}))} \rightarrow \theta \Gamma(l+1) G_0(\dd \omega^*),
\]
which leads to
\begin{align*}
    & \prob\left[ \sum_{i=1}^n \indicator_{(X_{n+1})}(X_i) = l, \bm C_{J} = \bm c, h(X_{n+1}) = j\right] \\
    & \qquad = \int_{D_j} \prob\left[\sum_{i=1}^n \indicator_{X_{n+1}}(X_i) = l, \bm C_{J} = \bm c, X_{n+1} \in \dd \omega^* \right] \\
    & \qquad = \frac{1}{l!}\binom{n}{c_1, \ldots, c_j - l, \ldots c_J} \frac{\Gamma(\theta)}{\Gamma(\theta + n + 1)} \theta \Gamma(l+1) G_0(B_{\omega^*}) \\
    & \qquad  \qquad \times \frac{\Gamma(\theta/J + c_j - l)}{\Gamma(\theta / J)} \prod_{k \neq j} \frac{\Gamma(\theta / J + c_k)}{\Gamma(\theta / J)} \int_{D_j} G_0(\dd \omega^*).
\end{align*}
Integration with respect to $\dd \omega^*$ is now straightforward and $\int_{D_j} G_0(\dd \omega^*) = 1 / J$.

Combining numerator and denominator yields the proof.

\subsection{Proof of \eqref{post_pyp} as a special case of Theorem~\ref{main_pk}}\label{app:post_pyp}

In the case of a PYP, we have $\beta = 0$, $\psi(u) = u^{\alpha}$ and
\[
    \kappa(u, l) = \alpha u^{\alpha - l} \frac{\Gamma(l - \alpha)}{\Gamma(1 - \alpha)} = \alpha u^{\alpha - l} (1 - \alpha)_{(l-1)}.
\]
Using the Faa di Bruno formula, we have \citep[see, e.g., Lemma 1 in][]{Cam(19)}:
\begin{align*}
    (-1)^c \frac{\dd^c}{\dd u^c} e^{-u^{\alpha} / J} &= e^{-u^{\alpha} / J} \sum_{i=0}^c \left(\frac{1}{J}\right)^i \sum_{(*)} \frac{1}{i!} \binom{c}{k_1 \cdots k_i} \prod_{j=1}^i \kappa(u, k_j) \\
    &= e^{-u^\alpha /J} \sum_{i=0}^c \left(\frac{1}{J}\right)^i \frac{\alpha^i}{u^{c - \alpha i} i!} \sum_{(*) } \binom{c}{k_1 \cdots k_i} \prod_{j=1}^i (1 - \alpha)_{k_j - 1} \\
    &= e^{-u^\alpha /J} \sum_{i=0}^c \left(\frac{1}{J}\right)^i \frac{1}{u^{c - \alpha i}} \mathscr{C}(c, i; \alpha),
\end{align*}
where $\mathscr{C}(c, i; \alpha)$ is the generalized factorial coefficient and $(*)$ denotes the summation over positive integers $(k_1, \ldots, k_i)$ such that $\sum_{j=1}^i k_j = c$.

Recall the definition of the multi-index set $S(\bm c, j, q)$ as in the main text below \eqref{post_pyp}.
Then, the integral at the numerator of \eqref{post_pk} equals
\begin{align*}
    & \int_{\R_+} u^{b+\gamma} e^{-u^{\alpha} / J} \sum_{i_j = 0}^{c_j - l} J^{-i_j} \frac{\mathscr{C}(c_j - l, i_j; \alpha)}{u^{c_j - l - \alpha i_j}} \\
    & \hspace{2cm} \times \Bigg\{ \prod_{k \neq j}  e^{-u^{\alpha} / J} \sum_{i_k = 0}^{c_k} J^{-i_k} \frac{\mathscr{C}(c_k, i_k; \alpha)}{u^{c_k - \alpha i_k}} \Bigg\} \alpha u^{\alpha - l - 1} (1 - \alpha)_{(l)} \dd u\\
    & = \sum_{i_1 = 0}^{c_1} \cdots \sum_{i_j = 0}^{c_j - l} \cdots \sum_{i_J = 0}^{c_J} J^{- \sum_k i_k} \mathscr{C}(c_j - l, i_j; \alpha) \\
    & \hspace{2cm} \times \prod_{k \neq j} \mathscr{C}(c_k, i_k; \alpha)  \alpha (1- \alpha)_{(l)} \int_{\R_+} e^{-u^\alpha} u^{\gamma - 1 + \alpha \sum_k i_k + \alpha} \dd u\\
    & = \sum_{i_1 = 0}^{c_1} \cdots \sum_{i_j = 0}^{c_j - l} \cdots \sum_{i_J = 0}^{c_J} J^{- \sum_k i_k} \mathscr{C}(c_j - l, i_j; \alpha) \\
    & \hspace{2cm} \times \prod_{k \neq j} \mathscr{C}(c_k, i_k; \alpha) (1 - \alpha)_{(l)} \Gamma\left(\frac{\gamma + \alpha}{\alpha} + \sum_k i_k \right) \\
    &= (1 - \alpha)_{(l)} \sum_{\bm i \in \mathcal S(\bm c, j, -l)} \Gamma\left(\frac{\gamma + \alpha}{ \alpha} + |\bm i|\right) J^{-|\bm i|} \prod_{k=1}^J \mathscr{C}(c_k - l \delta_{k,j}, i_k; \alpha).
\end{align*}
Similarly, the integral at the denominator  of \eqref{post_pk} equals
\begin{align*}
   & \sum_{i_1 = 0}^{c_1} \cdots \sum_{i_j = 0}^{c_j+1} \cdots \sum_{i_J = 0}^{c_J} J^{- \sum_k i_k} \mathscr{C}(c_j +1, i_j; \alpha) \prod_{k \neq j} \mathscr{C}(c_k, i_k; \alpha)  \Gamma\left(\frac{\gamma}{\alpha} + \sum_k i_k \right) \\
    &= \sum_{\bm i \in \mathcal S(\bm c, j, 1)} \Gamma\left(\gamma / \alpha + |\bm i|\right) J^{-|\bm i|} \prod_{k=1}^J \mathscr{C}(c_k + \delta_{k,j}, i_k; \alpha).
\end{align*}
Combining these gives \eqref{post_pyp}.

\subsection{Proof of Equation \eqref{post_dp} from Equation \eqref{post_pyp}}\label{app:post_pyp_new}

The proof relies on \citet[Theorem 2.16]{Cha(05)}, which characterizes the behaviour of generalized factorial coefficients as $\alpha\rightarrow0$. For $n\geq0$ and $0\leq k\leq n$ it holds that
\begin{equation}\label{gen_fact_stir}
\lim_{\alpha\rightarrow+\infty}\frac{\mathscr{C}(n,k;\alpha)}{\alpha^{k}}=|s(n,k)|,
\end{equation}
where $|s(n,k)|$ is the signless Stirling number of the first type; see \citet[Chapter 2]{Cha(05)} for details. Now, we apply \eqref{gen_fact_stir} to the posterior in \eqref{post_pyp}. For $l=0,1,\ldots,c_{j}$:
\begin{align*}
&\lim_{\alpha\rightarrow0}\text{Pr}[f_{X_{n+1}}=l\,|\,\mathbf{C}_{J}=\mathbf{c},h(X_{n+1})=j]\\
&\quad=\lim_{\alpha\rightarrow0}\frac{\gamma}{J} \binom{c_j}{l} (1 - \alpha)_{(l)} \frac{\sum_{\bm i \in \mathcal S(\bm c, j, -l)} \frac{\Gamma\left(\frac{\gamma + \alpha}{ \alpha} + |\bm i|\right)}{J^{|\bm i|}} \prod_{k=1}^J \mathscr{C}(c_k - l \delta_{k,j}, i_k; \alpha)}{\sum_{\bm i \in \mathcal S(\bm c, j, 1)}\frac{ \Gamma\left(\frac{\gamma}{\alpha} + |\bm i|\right)}{J^{|\bm i|}}  \prod_{k=1}^J \mathscr{C}(c_k + \delta_{k,j}, i_k; \alpha)}\\
&\quad=\lim_{\alpha\rightarrow0}\frac{\gamma}{J} \binom{c_j}{l} (1 - \alpha)_{(l)} \frac{\sum_{\bm i \in \mathcal S(\bm c, j, -l)} \frac{\prod_{t=0}^{|\bm{i}-1|}\left(\gamma+\alpha+\alpha t\right)}{J^{|\bm i|}} \prod_{k=1}^J \frac{\mathscr{C}(c_k - l \delta_{k,j}, i_k; \alpha)}{\alpha^{i_{k}}}}{\sum_{\bm i \in \mathcal S(\bm c, j, 1)}\frac{ \prod_{t=0}^{|\bm{i}-1|}\left(\gamma+\alpha t\right)}{J^{|\bm i|}}  \prod_{k=1}^J \frac{\mathscr{C}(c_k + \delta_{k,j}, i_k; \alpha)}{\alpha^{i_{k}}}}\\
&\quad=\frac{\gamma}{J} \binom{c_j}{l} l! \frac{\sum_{\bm i \in \mathcal S(\bm c, j, -l)} \left(\frac{\gamma}{J}\right)^{|\bm{i}|} \prod_{k=1}^J |s(c_k - l \delta_{k,j},i_{k})|}{\sum_{\bm i \in \mathcal S(\bm c, j, 1)}\left(\frac{ \gamma}{J}\right)^{\bm{i}}  \prod_{k=1}^J |s(c_k + \delta_{k,j},i_{k})|}.
\end{align*}
Recall that signless Stirling numbers of the first type are the coefficients of the series expansion of a rising factorial \citep[Equation 2.4]{Cha(05)}, i.e. for $t>0$, $(t)_{(n)}=\sum_{0\leq k\leq n}^{n}|s(n,k)|t^{k}$. Then,
\begin{align*}
\lim_{\alpha\rightarrow0}\text{Pr}[f_{X_{n+1}}=l\,|\,\mathbf{C}_{J}=\mathbf{c},h(X_{n+1})=j]
& =\frac{\gamma}{J} \binom{c_j}{l} l! \frac{\prod_{k=1}^{J}\left(\frac{\gamma}{J}\right)_{(c_k - l \delta_{k,j})}}{\prod_{k=1}^{J}\left(\frac{\gamma}{J}\right)_{(c_k + \delta_{k,j})}}\\
&=\frac{\gamma}{J}\frac{(c_{j}-l+1)_{(l)}}{\left(\frac{\gamma}{J}+c_{j}-l\right)_{(l+1)}},
\end{align*}
which coincides with the posterior in \eqref{post_dp} by setting $\theta=\gamma$. 

\subsection{Proof of Theorem~\ref{char_pyp}}\label{asypyp}

The first step, proved below, consists of establishing that, with $c_j$ fixed,
\begin{align}\label{post_pyp_asyn}
&\lim_{\mathbf{c}_{-j}\rightarrow+\infty} \prob[f_{X_{n+1}}=l\,|\,\mathbf{C}_{J}=\mathbf{c},h(X_{n+1})=j] \\
&\notag\quad\quad\quad\quad\quad =\gamma{c_{j}\choose l}(1-\alpha)_{(l)}\frac{\Gamma\left(\frac{\gamma}{\alpha}+J\right)(\gamma+J\alpha)_{(c_{j}-l)}}{\Gamma\left(\frac{\gamma}{\alpha}+J-1\right)(\gamma+J\alpha-\alpha)_{(c_{j}+1)}}.
\end{align}
Then, it follows by the Chu-Vandermonde identity \citep[Chapter 2]{Cha(05)} that
\begin{align*}
    & \sum_{l=0}^{c_j} l \gamma{c_{j}\choose l}(1-\alpha)_{(l)}\frac{\Gamma\left(\frac{\gamma}{\alpha}+J\right)(\gamma+J\alpha)_{(c_{j}-l)}}{\Gamma\left(\frac{\gamma}{\alpha}+J-1\right)(\gamma+J\alpha-\alpha)_{(c_{j}+1)}} \\
    & \quad = \gamma \frac{\Gamma\left(\frac{\gamma}{\alpha} + J \right)}{\Gamma\left(\frac{\gamma}{\alpha} + J -1\right)(\gamma + J\alpha - \alpha)_{(c_j + 1)}}\sum_{l=0}^{c_j} l \binom{c_j}{l} (1 - \alpha)_{(l)} (\gamma + J\alpha)_{(c_j - l)} \\
    & \quad = \gamma \frac{\Gamma\left(\frac{\gamma}{\alpha} + J \right)}{\Gamma\left(\frac{\gamma}{\alpha} + J -1\right)(\gamma + J\alpha - \alpha)_{(c_j + 1)}}  \sum_{l=0}^{c_j - 1} c_j \binom{c_j - 1}{l} \frac{\Gamma(1 - \alpha + 1 + l)}{\Gamma(1 - \alpha)} \frac{\Gamma(1 - \alpha + 1)}{\Gamma(1 - \alpha + 1)} (\gamma + J\alpha)_{c_j - 1 - l} \\
    & \quad =  \gamma \frac{\Gamma\left(\frac{\gamma}{\alpha} + J \right)}{\Gamma\left(\frac{\gamma}{\alpha} + J -1\right)(\gamma + J\alpha - \alpha)_{(c_j + 1)}} \frac{1}{\alpha} c_j (1 - \alpha) (2 + \gamma + J\alpha - \alpha)_{(c_j - 1)} \\
    & \quad = \frac{\gamma}{\alpha} c_j \frac{1 - \alpha}{\gamma + J\alpha - \alpha + 1}.
\end{align*}
Recalling the definition of $\hat{f}_{X_{n+1}}$ in~\eqref{eq:post-mean}, this directly leads to the desired result.

To prove~\eqref{post_pyp_asyn}, we rely on an asymptotic property of (normalized) generalized factorial coefficients. That is, from \citet[Lemma 2]{Dol(20)}, for any $z>0$ it holds that:
\begin{equation}\label{asym_gfc_n}
\lim_{n\rightarrow+\infty}\frac{z^{k}\mathscr{C}(n,k;s)}{\sum_{k=1}^{n}z^{k}\mathscr{C}(n,k;s)}=\text{e}^{-z}\frac{z^{k-1}}{(k-1)!}.
\end{equation}
We rewrite numerator and the denominator of \eqref{post_pyp}. For the numerator of \eqref{post_pyp},
\begin{align}\label{num_posterior}
&\frac{\gamma}{J}{c_{j}\choose l}(1-\alpha)_{(l)}\\
&\notag\quad\quad\times\sum_{i_{1}=1}^{c_{1}}J^{-i_{1}}\mathscr{C}(c_{1},i_{1};\alpha)\cdots\sum_{i_{j}=1}^{c_{j}-l}J^{-i_{j}}\mathscr{C}(c_{j}-l,i_{j};\alpha)\cdots\sum_{i_{J}=1}^{c_{J}}J^{-i_{J}}\mathscr{C}(c_{J},i_{J};\alpha)\\
&\notag\quad\quad\quad\times\Gamma\left(\frac{\gamma+\alpha}{\alpha}+i_{1}+\cdots+i_{j}+\cdots+i_{J}\right)\\
&\notag = \frac{\gamma}{J}{c_{j}\choose l}(1-\alpha)_{(l)}\Gamma\left(\frac{\gamma+\alpha}{\alpha}\right)\left(\prod_{1\leq s\neq j\leq J}\sum_{i_{s}=1}^{c_{s}}\frac{1}{J^{i_{s}}}\mathscr{C}(c_{s},i_{s};\alpha)\right)\\
&\notag\quad\quad\times\sum_{i_{j}=1}^{c_{j}-l}\left(\frac{1}{J}\right)^{i_{j}}\mathscr{C}(c_{j}-l,i_{j};\alpha)\sum_{i_{1}=1}^{c_{1}}\left(\frac{\gamma+\alpha}{\alpha}\right)_{(i_{1})}\frac{\left(\frac{1}{J}\right)^{i_{1}}\mathscr{C}(c_{1},i_{1};\alpha)}{\sum_{i_{1}=1}^{c_{1}}\left(\frac{1}{J}\right)^{i_{1}}\mathscr{C}(c_{1},i_{1};\alpha)}\cdots\\
&\notag\quad\quad\quad\cdots\times\left(\frac{\gamma+\alpha}{\alpha}+i_{1}+\cdots+i_{j-1}\right)_{(i_{j})}\cdots\\
&\notag\quad\quad\quad\quad\cdots\times\sum_{i_{J}=1}^{c_{J}}\left(\frac{\gamma+\alpha}{\alpha}+i_{1}+\cdots+i_{J-1}\right)_{(i_{J})}\frac{\left(\frac{1}{J}\right)^{i_{J}}\mathscr{C}(c_{J},i_{J};\alpha)}{\sum_{i_{J}=1}^{c_{J}}\left(\frac{1}{J}\right)^{i_{J}}\mathscr{C}(c_{J},i_{J};\alpha)},
\end{align}
and for the denominator of \eqref{post_pyp}
\begin{align}\label{den_posterior}
&\sum_{i_{1}=1}^{c_{1}}J^{-i_{1}}\mathscr{C}(c_{1},i_{1};\alpha)\cdots\sum_{i_{j}=1}^{c_{j}+1}J^{-i_{j}}\mathscr{C}(c_{j}+1,i_{j};\alpha)\cdots\sum_{i_{J}=1}^{c_{J}}J^{-i_{J}}\mathscr{C}(c_{J},i_{J};\alpha)\\
&\notag\quad\times\Gamma\left(\frac{\gamma}{\alpha}+i_{1}+\cdots+i_{j}+\cdots+i_{J}\right)\\
&\notag\quad=\Gamma\left(\frac{\gamma}{\alpha}\right)\left(\prod_{1\leq s\neq j\leq J}\sum_{i_{s}=1}^{c_{s}}\frac{1}{J^{i_{s}}}\mathscr{C}(c_{s},i_{s};\alpha)\right)\\
&\notag\quad\quad\times\sum_{i_{j}=1}^{c_{j}+1}\left(\frac{1}{J}\right)^{i_{j}}\mathscr{C}(c_{j}+1,i_{j};\alpha)\sum_{i_{1}=1}^{c_{1}}\left(\frac{\gamma}{\alpha}\right)_{(i_{1})}\frac{\left(\frac{1}{J}\right)^{i_{1}}\mathscr{C}(c_{1},i_{1};\alpha)}{\sum_{i_{1}=1}^{c_{1}}\left(\frac{1}{J}\right)^{i_{1}}\mathscr{C}(c_{1},i_{1};\alpha)}\cdots\\
&\notag\quad\quad\quad\cdots\times\left(\frac{\gamma}{\alpha}+i_{1}+\cdots+i_{j-1}\right)_{(i_{j})}\cdots\\
&\notag\quad\quad\quad\quad\cdots\times\sum_{i_{J}=1}^{c_{J}}\left(\frac{\gamma}{\alpha}+i_{1}+\cdots+i_{J-1}\right)_{(i_{J})}\frac{\left(\frac{1}{J}\right)^{i_{J}}\mathscr{C}(c_{J},i_{J};\alpha)}{\sum_{i_{J}=1}^{c_{J}}\left(\frac{1}{J}\right)^{i_{J}}\mathscr{C}(c_{J},i_{J};\alpha)}.
\end{align}
By combining \eqref{num_posterior} with \eqref{den_posterior}, we write the posterior distribution \eqref{post_pyp} as
\begin{align}\label{main_num_den}
&\text{Pr}[f_{X_{n+1}}=l\,|\,\mathbf{C}_{J}=\mathbf{c},h(X_{n+1})=j]=\frac{\gamma}{J}{c_{j}\choose l}(1-\alpha)_{(l)}\frac{\Gamma\left(\frac{\gamma+\alpha}{\alpha}\right)}{\Gamma\left(\frac{\gamma}{\alpha}\right)}\frac{N_{\gamma,\alpha,J}(l;c_{1},\ldots,c_{J})}{D_{\gamma,\alpha,J}(c_{1},\ldots,c_{J})},
\end{align}
where
\begin{align}\label{num_rewrite}
&N_{\gamma,\alpha,J}(l;c_{1},\ldots,c_{J})\\
&\notag\quad=\sum_{i_{j}=1}^{c_{j}-l}\left(\frac{1}{J}\right)^{i_{j}}\mathscr{C}(c_{j}-l,i_{j};\alpha)\sum_{i_{1}=1}^{c_{1}}\left(\frac{\gamma+\alpha}{\alpha}\right)_{(i_{1})}\frac{\left(\frac{1}{J}\right)^{i_{1}}\mathscr{C}(c_{1},i_{1};\alpha)}{\sum_{i_{1}=1}^{c_{1}}\left(\frac{1}{J}\right)^{i_{1}}\mathscr{C}(c_{1},i_{1};\alpha)}\cdots\\
&\notag\quad\quad\cdots\times\left(\frac{\gamma+\alpha}{\alpha}+i_{1}+\cdots+i_{j-1}\right)_{(i_{j})}\cdots\\
&\notag\quad\quad\quad\cdots\times\sum_{i_{J}=1}^{c_{J}}\left(\frac{\gamma+\alpha}{\alpha}+i_{1}+\cdots+i_{J-1}\right)_{(i_{J})}\frac{\left(\frac{1}{J}\right)^{i_{J}}\mathscr{C}(c_{J},i_{J};\alpha)}{\sum_{i_{J}=1}^{c_{J}}\left(\frac{1}{J}\right)^{i_{J}}\mathscr{C}(c_{J},i_{J};\alpha)}
\end{align}
and
\begin{align}\label{den_rewrite}
&D_{\gamma,\alpha,J}(c_{1},\ldots,c_{J})\\
&\notag\quad=\sum_{i_{j}=1}^{c_{j}+1}\left(\frac{1}{J}\right)^{i_{j}}\mathscr{C}(c_{j}+1,i_{j};\alpha)\sum_{i_{1}=1}^{c_{1}}\left(\frac{\gamma}{\alpha}\right)_{(i_{1})}\frac{\left(\frac{1}{J}\right)^{i_{1}}\mathscr{C}(c_{1},i_{1};\alpha)}{\sum_{i_{1}=1}^{c_{1}}\left(\frac{1}{J}\right)^{i_{1}}\mathscr{C}(c_{1},i_{1};\alpha)}\cdots\\
&\notag\quad\quad\cdots\times\left(\frac{\gamma}{\alpha}+i_{1}+\cdots+i_{j-1}\right)_{(i_{j})}\cdots\\
&\notag\quad\quad\quad\cdots\times\sum_{i_{J}=1}^{c_{J}}\left(\frac{\gamma}{\alpha}+i_{1}+\cdots+i_{J-1}\right)_{(i_{J})}\frac{\left(\frac{1}{J}\right)^{i_{J}}\mathscr{C}(c_{J},i_{J};\alpha)}{\sum_{i_{J}=1}^{c_{J}}\left(\frac{1}{J}\right)^{i_{J}}\mathscr{C}(c_{J},i_{J};\alpha)}.
\end{align}
Now, we apply repeatedly \eqref{asym_gfc_n} to both \eqref{num_rewrite} and \eqref{den_rewrite}, starting from the last terms, indexed by $i_{J}$. In particular, for the sum in $i_{J}=1,\ldots,c_{J}$ of \eqref{num_rewrite}, we have that
\begin{align}\label{first_step_num}
&\lim_{c_{J}\rightarrow+\infty}\sum_{i_{J}=1}^{c_{J}}\left(\frac{\gamma+\alpha}{\alpha}+i_{1}+\cdots+i_{J-1}\right)_{(i_{J})}\frac{\left(\frac{1}{J}\right)^{i_{J}}\mathscr{C}(c_{J},i_{J};\alpha)}{\sum_{i_{J}=1}^{c_{J}}\left(\frac{1}{J}\right)^{i_{J}}\mathscr{C}(c_{J},i_{J};\alpha)}\\
&\notag\quad=\sum_{i_{J}\geq1}\left(\frac{\gamma+\alpha}{\alpha}+i_{1}+\cdots+i_{J-1}\right)_{(i_{J})}\text{e}^{-\frac{1}{J}}\frac{\left(\frac{1}{J}\right)^{i_{J}-1}}{(i_{J}-1)!}\\
&\notag\quad=\text{e}^{-\frac{1}{J}}\frac{\frac{\gamma+\alpha}{\alpha}+i_{1}+\cdots+i_{J-1}}{\left(1-\frac{1}{J}\right)^{\frac{\gamma+\alpha}{\alpha}+i_{1}+\cdots+i_{J-1}+1}},
\end{align}
and similarly, for the sum in $i_{J}=1,\ldots,c_{J}$ of \eqref{den_rewrite}, we have:
\begin{align}\label{first_step_den}
&\lim_{c_{J}\rightarrow+\infty}\sum_{i_{J}=1}^{c_{J}}\left(\frac{\gamma}{\alpha}+i_{1}+\cdots+i_{J-1}\right)_{(i_{J})}\frac{\left(\frac{1}{J}\right)^{i_{J}}\mathscr{C}(c_{J},i_{J};\alpha)}{\sum_{i_{J}=1}^{c_{J}}\left(\frac{1}{J}\right)^{i_{J}}\mathscr{C}(c_{J},i_{J};\alpha)}\\
&\notag\quad=\text{e}^{-\frac{1}{J}}\frac{\frac{\gamma}{\alpha}+i_{1}+\cdots+i_{J-1}}{\left(1-\frac{1}{J}\right)^{\frac{\gamma}{\alpha}+i_{1}+\cdots+i_{J-1}+1}}.
\end{align}
Now, consider the sum in $i_{J-1}=1,\ldots,c_{J-1}$ of \eqref{num_rewrite}, combined with \eqref{first_step_num}. That is,
\begin{align*}
&\lim_{c_{J-1}\rightarrow+\infty}\sum_{i_{J-1}=1}^{c_{J-1}}\left(\frac{\gamma+\alpha}{\alpha}+i_{1}+\cdots+i_{J-2}\right)_{(i_{J-1})}\text{e}^{-\frac{1}{J}}\frac{\frac{\gamma+\alpha}{\alpha}+i_{1}+\cdots+i_{J-1}}{\left(1-\frac{1}{J}\right)^{\frac{\gamma+\alpha}{\alpha}+i_{1}+\cdots+i_{J-1}+1}}\\
&\quad\quad\times\frac{\left(\frac{1}{J}\right)^{i_{J-1}}\mathscr{C}(c_{J-1},i_{J-1};\alpha)}{\sum_{i_{J-1}=1}^{c_{J-1}}\left(\frac{1}{J}\right)^{i_{J-1}}\mathscr{C}(c_{J-1},i_{J-1};\alpha)}\\
&\quad=\sum_{i_{J-1}\geq1}\left(\frac{\gamma+\alpha}{\alpha}+i_{1}+\cdots+i_{J-2}\right)_{(i_{J-1})}\text{e}^{-\frac{1}{J}}\frac{\frac{\gamma+\alpha}{\alpha}+i_{1}+\cdots+i_{J-1}}{\left(1-\frac{1}{J}\right)^{\frac{\gamma+\alpha}{\alpha}+i_{1}+\cdots+i_{J-1}+1}}\text{e}^{-\frac{1}{J}}\frac{\left(\frac{1}{J}\right)^{i_{J-1}-1}}{(i_{J-1}-1)!}\\
&\quad=\frac{\text{e}^{-\frac{2}{J}}}{\left(1-\frac{1}{J}\right)^{\frac{\gamma+\alpha}{\alpha}+i_{1}+\cdots+i_{J-2}+2}}\sum_{i_{J-1}\geq1}\left(\frac{\gamma+\alpha}{\alpha}+i_{1}+\cdots+i_{J-2}\right)_{(i_{J-1}+1)}\frac{\left(\frac{1}{J-1}\right)^{i_{J-1}-1}}{(i_{J-1}-1)!}\\
&\quad=\frac{\text{e}^{-\frac{2}{J}}}{\left(1-\frac{1}{J}\right)^{\frac{\gamma+\alpha}{\alpha}+i_{1}+\cdots+i_{J-2}+2}}\frac{\left(\frac{\gamma+\alpha}{\alpha}+i_{1}+\cdots+i_{J-2}\right)_{(2)}}{\left(\frac{J-2}{J-1}\right)^{\frac{\gamma+\alpha}{\alpha}+i_{1}+\cdots+i_{J-2}+2}}\\
&\quad=\text{e}^{-\frac{2}{J}}\frac{\left(\frac{\gamma+\alpha}{\alpha}+i_{1}+\cdots+i_{J-2}\right)_{(2)}}{\left(1-\frac{2}{J}\right)^{\frac{\gamma+\alpha}{\alpha}+i_{1}+\cdots+i_{J-2}+2}}.
\end{align*}
Similarly, consider the sum in $i_{J-1}=1,\ldots,c_{J-1}$ of \eqref{den_rewrite}, together with \eqref{first_step_den}, which leads to
\begin{align*}
&\lim_{c_{J-1}\rightarrow+\infty}\sum_{i_{J-1}=1}^{c_{J-1}}\left(\frac{\gamma}{\alpha}+i_{1}+\cdots+i_{J-2}\right)_{(i_{J-1})}\text{e}^{-\frac{1}{J}}\frac{\frac{\gamma}{\alpha}+i_{1}+\cdots+i_{J-1}}{\left(1-\frac{1}{J}\right)^{\frac{\gamma}{\alpha}+i_{1}+\cdots+i_{J-1}+1}}\\
&\quad\quad\times\frac{\left(\frac{1}{J}\right)^{i_{J-1}}\mathscr{C}(c_{J-1},i_{J-1};\alpha)}{\sum_{i_{J-1}=1}^{c_{J-1}}\left(\frac{1}{J}\right)^{i_{J-1}}\mathscr{C}(c_{J-1},i_{J-1};\alpha)}\\
&\quad=\text{e}^{-\frac{2}{J}}\frac{\left(\frac{\gamma}{\alpha}+i_{1}+\cdots+i_{J-2}\right)_{(2)}}{\left(1-\frac{2}{J}\right)^{\frac{\gamma}{\alpha}+i_{1}+\cdots+i_{J-2}+2}}.
\end{align*}
By proceeding recursively, for the sum in $i_{j+1}=1,\ldots,c_{j+1}$ of \eqref{num_rewrite}, we find:
\begin{align}\label{general_nun}
&\lim_{c_{j+1}\rightarrow+\infty}\sum_{i_{j+1}=1}^{c_{j+1}}\left(\frac{\gamma+\alpha}{\alpha}+i_{1}+\cdots+i_{j}\right)_{(i_{j+1})}\text{e}^{-\frac{J-j-1}{J}}\frac{\left(\frac{\gamma+\alpha}{\alpha}+i_{1}+\cdots+i_{j+1}\right)_{(J-j-1)}}{\left(1-\frac{J-j-1}{J}\right)^{\frac{\gamma+\alpha}{\alpha}+i_{1}+\cdots+i_{j+1}+J-j-1}}\\
&\notag\quad\quad\times\frac{\left(\frac{1}{J}\right)^{i_{j+1}}\mathscr{C}(c_{j+1},i_{j+1};\alpha)}{\sum_{i_{j+1}=1}^{c_{j+1}}\left(\frac{1}{J}\right)^{i_{j+1}}\mathscr{C}(c_{j+1},i_{j+1};\alpha)}\\
&\notag\quad=\sum_{c_{j+1}\geq1}\left(\frac{\gamma+\alpha}{\alpha}+i_{1}+\cdots+i_{j}\right)_{(i_{j+1})}\text{e}^{-\frac{J-j-1}{J}}\frac{\left(\frac{\gamma+\alpha}{\alpha}+i_{1}+\cdots+i_{j+1}\right)_{(J-j-1)}}{\left(1-\frac{J-j-1}{J}\right)^{\frac{\gamma+\alpha}{\alpha}+i_{1}+\cdots+i_{j+1}+J-j-1}}\text{e}^{-\frac{1}{J}}\frac{\left(\frac{1}{J}\right)^{i_{j+1}-1}}{(i_{j+1}-1)!}\\
&\notag\quad=\frac{\text{e}^{-\frac{J-j}{J}}}{\left(1-\frac{J-j-1}{J}\right)^{\frac{\gamma+\alpha}{\alpha}+i_{1}+\cdots+i_{j}+J-j}}\sum_{i_{j+1}\geq1}\left(\frac{\gamma+\alpha}{\alpha}+i_{1}+\cdots+i_{j}\right)_{(i_{j+1}+J-j-1)}\frac{\left(\frac{1}{j+1}\right)^{i_{j+1}-1}}{(i_{j+1}-1)!}\\
&\notag\quad=\frac{\text{e}^{-\frac{J-j}{J}}}{\left(1-\frac{J-j-1}{J}\right)^{\frac{\gamma+\alpha}{\alpha}+i_{1}+\cdots+i_{j}+J-j}}\frac{\left(\frac{\gamma+\alpha}{\alpha}+i_{1}+\cdots+i_{j}\right)_{(J-j)}}{\left(\frac{j}{j+1}\right)^{\frac{\gamma+\alpha}{\alpha}+i_{1}+\cdots+i_{j}+J-j}}\\
&\notag\quad=\text{e}^{-\frac{J-j}{J}}\frac{\left(\frac{\gamma+\alpha}{\alpha}+i_{1}+\cdots+i_{j}\right)_{(J-j)}}{\left(1-\frac{J-j}{J}\right)^{\frac{\gamma+\alpha}{\alpha}+i_{1}+\cdots+i_{j}+J-j}},
\end{align}
and, similarly, the sum in $i_{j+1}=1,\ldots,c_{j+1}$ of \eqref{den_rewrite} it leads to the following:
\begin{align}\label{general_den}
&\lim_{c_{j+1}\rightarrow+\infty}\sum_{i_{j+1}=1}^{c_{j+1}}\left(\frac{\gamma}{\alpha}+i_{1}+\cdots+i_{j}\right)_{(i_{j+1})}\text{e}^{-\frac{J-j-1}{J}}\frac{\left(\frac{\gamma}{\alpha}+i_{1}+\cdots+i_{j+1}\right)_{(J-j-1)}}{\left(1-\frac{J-j-1}{J}\right)^{\frac{\gamma}{\alpha}+i_{1}+\cdots+i_{j+1}+J-j-1}}\\
&\notag\quad\quad\times\frac{\left(\frac{1}{J}\right)^{i_{j+1}}\mathscr{C}(c_{j+1},i_{j+1};\alpha)}{\sum_{i_{j+1}=1}^{c_{j+1}}\left(\frac{1}{J}\right)^{i_{j+1}}\mathscr{C}(c_{j+1},i_{j+1};\alpha)}\\
&\notag\quad=\text{e}^{-\frac{J-j}{J}}\frac{\left(\frac{\gamma}{\alpha}+i_{1}+\cdots+i_{j}\right)_{(J-j)}}{\left(1-\frac{J-j}{J}\right)^{\frac{\gamma}{\alpha}+i_{1}+\cdots+i_{j}+J-j}}.
\end{align}
Now, we consider the sum in $i_{j-1}=1,\ldots,c_{j-1}$ of \eqref{num_rewrite}, together with \eqref{general_nun}. That is,
\begin{align*}
&\lim_{c_{j-1}\rightarrow+\infty}\sum_{i_{j-1}=1}^{c_{j-1}}\left(\frac{\gamma+\alpha}{\alpha}+i_{1}+\cdots+i_{j-2}\right)_{(i_{j-1})}\left(\frac{\gamma+\alpha}{\alpha}+i_{1}+\cdots+i_{j-1}\right)_{(i_{j})}\\
&\quad\quad\times\text{e}^{-\frac{J-j}{J}}\frac{\left(\frac{\gamma+\alpha}{\alpha}+i_{1}+\cdots+i_{j}\right)_{(J-j)}}{\left(1-\frac{J-j}{J}\right)^{\frac{\gamma+\alpha}{\alpha}+i_{1}+\cdots+i_{j}+J-j}}\frac{\left(\frac{1}{J}\right)^{i_{j-1}}\mathscr{C}(c_{j-1},i_{j-1};\alpha)}{\sum_{i_{j-1}=1}^{c_{j-1}}\left(\frac{1}{J}\right)^{i_{j-1}}\mathscr{C}(c_{j-1},i_{j-1};\alpha)}\\
&\quad=\sum_{i_{j-1}\geq1}\left(\frac{\gamma+\alpha}{\alpha}+i_{1}+\cdots+i_{j-2}\right)_{(i_{j-1})}\left(\frac{\gamma+\alpha}{\alpha}+i_{1}+\cdots+i_{j-1}\right)_{(i_{j})}\\
&\quad\quad\times\text{e}^{-\frac{J-j}{J}}\frac{\left(\frac{\gamma+\alpha}{\alpha}+i_{1}+\cdots+i_{j}\right)_{(J-j)}}{\left(1-\frac{J-j}{J}\right)^{\frac{\gamma+\alpha}{\alpha}+i_{1}+\cdots+i_{j}+J-j}}\text{e}^{-\frac{1}{J}}\frac{\left(\frac{1}{J}\right)^{i_{j-1}-1}}{(i_{j-1}-1)!}\\
&\quad=\frac{\text{e}^{-\frac{J-j+1}{J}}}{\left(1-\frac{J-j}{J}\right)^{\frac{\gamma+\alpha}{\alpha}+i_{j}+i_{1}+\cdots+i_{j-2}+J-j+1}}\sum_{i_{j-1}\geq1}\left(\frac{\gamma+\alpha}{\alpha}+i_{1}+\cdots+i_{j-2}\right)_{(i_{j-1}+i_{j}+J-j)}\frac{\left(\frac{1}{j}\right)^{i_{j-1}-1}}{(i_{j-1}-1)!}\\
&\quad=\frac{\text{e}^{-\frac{J-j+1}{J}}}{\left(1-\frac{J-j}{J}\right)^{\frac{\gamma+\alpha}{\alpha}+i_{j}+i_{1}+\cdots+i_{j-2}+J-j+1}}\frac{\left(\frac{\gamma+\alpha}{\alpha}+i_{1}+\cdots+i_{j-2}\right)_{(i_{j}+J-j+1)}}{\left(\frac{j-1}{j}\right)^{\frac{\gamma+\alpha}{\gamma}+i_{j}+i_{1}+\cdots+i_{j-2}+J-j+1}}\\
&\quad=\text{e}^{-\frac{J-j+1}{J}}\frac{\left(\frac{\gamma+\alpha}{\alpha}+i_{1}+\cdots+i_{j-2}\right)_{(i_{j}+J-j+1)}}{\left(\frac{j-1}{J}\right)^{\frac{\gamma+\alpha}{\gamma}+i_{j}+i_{1}+\cdots+i_{j-2}+J-j+1}}.
\end{align*}
Similarly, consider the sum in $i_{j-1}=1,\ldots,c_{j-1}$ of \eqref{den_rewrite}, together with \eqref{general_den}, which leads to
\begin{align*}
&\lim_{c_{j-1}\rightarrow+\infty}\sum_{i_{j-1}=1}^{c_{j-1}}\left(\frac{\gamma}{\alpha}+i_{1}+\cdots+i_{j-2}\right)_{(i_{j-1})}\left(\frac{\gamma}{\alpha}+i_{1}+\cdots+i_{j-1}\right)_{(i_{j})}\\
&\quad\quad\times\text{e}^{-\frac{J-j}{J}}\frac{\left(\frac{\gamma}{\alpha}+i_{1}+\cdots+i_{j}\right)_{(J-j)}}{\left(1-\frac{J-j}{J}\right)^{\frac{\gamma}{\alpha}+i_{1}+\cdots+i_{j}+J-j}}\frac{\left(\frac{1}{J}\right)^{i_{j-1}}\mathscr{C}(c_{j-1},i_{j-1};\alpha)}{\sum_{i_{j-1}=1}^{c_{j-1}}\left(\frac{1}{J}\right)^{i_{j-1}}\mathscr{C}(c_{j-1},i_{j-1};\alpha)}\\
&\quad=\text{e}^{-\frac{J-j+1}{J}}\frac{\left(\frac{\gamma}{\alpha}+i_{1}+\cdots+i_{j-2}\right)_{(i_{j}+J-j+1)}}{\left(\frac{j-1}{J}\right)^{\frac{\gamma}{\gamma}+i_{j}+i_{1}+\cdots+i_{j-2}+J-j+1}}.
\end{align*}
By proceeding recursively, we arrive to the sum in $i_{1}=1,\ldots,c_{1}$ of \eqref{num_rewrite}. That is,
\begin{align}\label{final_num}
&\lim_{c_{1}\rightarrow+\infty}\sum_{i_{1}=1}^{c_{1}}\left(\frac{\gamma+\alpha}{\alpha}\right)_{(i_{1})}\text{e}^{-\frac{J-2}{J}}\frac{\left(\frac{\gamma+\alpha}{\alpha}+i_{1}\right)_{(i_{j}+J-2)}}{\left(\frac{2}{J}\right)^{\frac{\gamma+\alpha}{\alpha}+i_{j}+i_{1}+J-2}}\frac{\left(\frac{1}{J}\right)^{i_{1}}\mathscr{C}(c_{1},i_{1};\alpha)}{\sum_{i_{1}=1}^{c_{1}}\left(\frac{1}{J}\right)^{i_{1}}\mathscr{C}(c_{1},i_{1};\alpha)}\\
&\notag\quad=\sum_{i_{1}\geq1}\left(\frac{\gamma+\alpha}{\alpha}\right)_{(i_{1})}\text{e}^{-\frac{J-2}{J}}\frac{\left(\frac{\gamma+\alpha}{\alpha}+i_{1}\right)_{(i_{j}+J-2)}}{\left(\frac{2}{J}\right)^{\frac{\gamma+\alpha}{\alpha}+i_{j}+i_{1}+J-2}}\text{e}^{-\frac{1}{J}}\frac{\left(\frac{1}{J}\right)^{i_{1}-1}}{(i_{1}-1)!}\\
&\notag\quad=\frac{\text{e}^{-\frac{J-1}{J}}}{\left(\frac{2}{J}\right)^{\frac{\gamma+\alpha}{\alpha}+i_{j}+J-1}}\sum_{i_{1}\geq1}\left(\frac{\gamma+\alpha}{\alpha}\right)_{(i_{1}+i_{j}+J-2)}\frac{\left(\frac{1}{2}\right)^{i_{1}-1}}{(i_{1}-1)!}\\
&\notag\quad=\frac{\text{e}^{-\frac{J-1}{J}}}{\left(\frac{2}{J}\right)^{\frac{\gamma+\alpha}{\alpha}+i_{j}+J-1}}\frac{\left(\frac{\gamma+\alpha}{\alpha}\right)_{(i_{j}+J-1)}}{\left(\frac{1}{2}\right)^{\frac{\gamma+\alpha}{\alpha}+i_{j}+J-1}}\\
&\notag\quad=\text{e}^{-\frac{J-1}{J}}\frac{\left(\frac{\gamma+\alpha}{\alpha}\right)_{(i_{j}+J-1)}}{\left(\frac{1}{J}\right)^{\frac{\gamma+\alpha}{\alpha}+i_{j}+J-1}}
\end{align}
and, similarly, we arrive to the sum in $i_{1}=1,\ldots,c_{1}$ of \eqref{den_rewrite}, which leads to
\begin{align}\label{final_den}
&\lim_{c_{1}\rightarrow+\infty}\sum_{i_{1}=1}^{c_{1}}\left(\frac{\gamma}{\alpha}\right)_{(i_{1})}\text{e}^{-\frac{J-2}{J}}\frac{\left(\frac{\gamma}{\alpha}+i_{1}\right)_{(i_{j}+J-2)}}{\left(\frac{2}{J}\right)^{\frac{\gamma}{\alpha}+i_{j}+i_{1}+J-2}}\frac{\left(\frac{1}{J}\right)^{i_{1}}\mathscr{C}(c_{1},i_{1};\alpha)}{\sum_{i_{1}=1}^{c_{1}}\left(\frac{1}{J}\right)^{i_{1}}\mathscr{C}(c_{1},i_{1};\alpha)}\\
&\notag\quad=\text{e}^{-\frac{J-1}{J}}\frac{\left(\frac{\gamma}{\alpha}\right)_{(i_{j}+J-1)}}{\left(\frac{1}{J}\right)^{\frac{\gamma}{\alpha}+i_{j}+J-1}}.
\end{align}
From \eqref{final_num} and \eqref{final_den}, and \eqref{main_num_den},
\begin{align*}
&\lim_{\mathbf{c}_{-j} \rightarrow +\infty} \frac{\gamma}{J}{c_{j}\choose l}(1-\alpha)_{(l)}\frac{\Gamma\left(\frac{\gamma+\alpha}{\alpha}\right)}{\Gamma\left(\frac{\gamma}{\alpha}\right)}\frac{N_{\gamma,\alpha,J}(l;c_{1},\ldots,c_{J})}{D_{\gamma,\alpha,J}(c_{1},\ldots,c_{J})}\\
&\quad= \frac{\gamma}{J}{c_{j}\choose l}(1-\alpha)_{(l)}\frac{\Gamma\left(\frac{\gamma+\alpha}{\alpha}\right)}{\Gamma\left(\frac{\gamma}{\alpha}\right)}\frac{\sum_{i_{j}=1}^{c_{j}-l}\left(\frac{1}{J}\right)^{i_{j}}\mathscr{C}(c_{j}-l,i_{j};\alpha)\text{e}^{-\frac{J-1}{J}}\frac{\left(\frac{\gamma+\alpha}{\alpha}\right)_{(i_{j}+J-1)}}{\left(\frac{1}{J}\right)^{\frac{\gamma+\alpha}{\alpha}+i_{j}+J-1}}}{\sum_{i_{j}=1}^{c_{j}+1}\left(\frac{1}{J}\right)^{i_{j}}\mathscr{C}(c_{j}+1,i_{j};\alpha)\text{e}^{-\frac{J-1}{J}}\frac{\left(\frac{\gamma}{\alpha}\right)_{(i_{j}+J-1)}}{\left(\frac{1}{J}\right)^{\frac{\gamma}{\alpha}+i_{j}+J-1}}}\\
&\quad=\gamma{c_{j}\choose l}(1-\alpha)_{(l)}\frac{\Gamma\left(\frac{\gamma+\alpha}{\alpha}+J-1\right)}{\Gamma\left(\frac{\gamma}{\alpha}+J-1\right)}\frac{\sum_{i_{j}=1}^{c_{j}-l}\mathscr{C}(c_{j}-l,i_{j};\alpha)\left(\frac{\gamma+\alpha}{\alpha}+J-1\right)_{(i_{j})}}{\sum_{i_{j}=1}^{c_{j}+1}\mathscr{C}(c_{j}+1,i_{j};\alpha)\left(\frac{\gamma}{\alpha}+J-1\right)_{(i_{j})}}\\
&\quad=\gamma{c_{j}\choose l}(1-\alpha)_{(l)}\frac{\Gamma\left(\frac{\gamma+\alpha}{\alpha}+J-1\right)}{\Gamma\left(\frac{\gamma}{\alpha}+J-1\right)}\frac{(\gamma+\alpha+J\alpha-\alpha)_{(c_{j}-l)}}{(\gamma+J\alpha-\alpha)_{(c_{j}+1)}},
\end{align*}
where the sums over $i_{j}$ follows from \citet[Equation 2.49]{Cha(05)}. 

\subsection{Proof of Corollary~\ref{corol_card}}\label{card_proof}
The proof follows the same arguments developed in Appendix~\ref{app:post_pyp}. According to \eqref{card_main_6}, it is sufficient to compute the conditional probability of $h(X_{n+1})$, given $\mathbf{C}_{J}$, as the conditional probability of $f_{X_{n+1}}$, given $\mathbf{C}_{J}$ and $h(X_{n+1})$ is available from \eqref{post_pyp}. For a PK prior,
\begin{displaymath}
\text{Pr}[h(X_{n+1})=j\,|\,\mathbf{C}_{J}=\mathbf{c}]=\frac{\text{Pr}[\mathbf{C}_{J}=\mathbf{c},h(X_{n+1})=j]}{\text{Pr}[\mathbf{C}_{J}=\mathbf{c}]}
\end{displaymath}
for all $j\in[J]$, where, from \eqref{eq:denom_probs},
\begin{align}\label{num_card}
\text{Pr}[\mathbf{C}_{J}=\mathbf{c},h(X_{n+1})=j]&={n\choose c_{1},\ldots,c_{J}}\E\left[P(D_{j})^{c_{j}+1}\prod_{k\neq j}P(D_{k})^{c_{k}}\right]\\
&\notag=Z_T \binom{n}{c_1, \ldots, c_j} \int_{\R_+ } \frac{u^{n+\gamma}}{\Gamma(n + \gamma + 1)} (-1)^{c+1} \frac{\dd^{c_j+1}}{\dd z^{c_j+1}} e^{- \theta \psi(z)/J}|_{(u + \beta)}\\
&\notag\quad\times\prod_{k\neq j} (-1)^{c_k}\frac{\dd^{c_k}}{\dd z^{c_k}} e^{- \theta \psi(z)/J}|_{(u + \beta)} \dd u
\end{align}
and
\begin{align}\label{denom_card}
\text{Pr}[\mathbf{C}_{J}=\mathbf{c}]&={n\choose c_{1},\ldots,c_{J}}\E\left[P(D_{j})^{c_{j}}\prod_{k\neq j}P(D_{k})^{c_{k}}\right]\\
&\notag=Z_T \binom{n}{c_1, \ldots, c_j} \int_{\R_+ } \frac{u^{n+\gamma}}{\Gamma(n + \gamma + 1)}\prod_{k=1}^{J} (-1)^{c_k}\frac{\dd^{c_k}}{\dd z^{c_k}} e^{- \theta \psi(z)/J}|_{(u + \beta)} \dd u.
\end{align}
In the case of a PYP, we have $\beta = 0$, $\psi(u) = u^{\alpha}$ and
\begin{displaymath}
\kappa(u, l) = \alpha u^{\alpha - l} \frac{\Gamma(l - \alpha)}{\Gamma(1 - \alpha)} = \alpha u^{\alpha - l} (1 - \alpha)_{(l-1)}.
\end{displaymath}
Then, by applying to \eqref{num_card} and \eqref{denom_card} the same arguments of Appendix~\ref{app:post_pyp}, we obtain:
\begin{align*}
&\text{Pr}[\mathbf{C}_{J}=\mathbf{c},h(X_{n+1})=j]={n\choose c_{1},\ldots,c_{J}}\sum_{\bm{i}\in\mathcal{S}(\mathbf{c},j,1)}\frac{\frac{\left(\frac{\gamma}{\alpha}\right)_{(|\bm{i}|)}}{J^{|\bm{i}|}}}{(\gamma)_{(n+1)}}\prod_{k=1}^{J}\mathscr{C}(c_{k}+\delta_{k,j},i_{k};\alpha)
\end{align*}
and
\begin{align*}
&\text{Pr}[\mathbf{C}_{J}=\mathbf{c}]={n\choose c_{1},\ldots,c_{J}}\sum_{\bm{i}\in\mathcal{S}(\mathbf{c},j,0)}\frac{\frac{\left(\frac{\gamma}{\alpha}\right)_{(|\bm{i}|)}}{J^{|\bm{i}|}}}{(\gamma)_{(n)}}\prod_{k=1}^{J}\mathscr{C}(c_{k},i_{k};\alpha).
\end{align*}
Then,
\begin{align}\label{ratio_card}
\text{Pr}[h(X_{n+1})=j\,|\,\mathbf{C}_{J}=\mathbf{c}] 
&\quad =\frac{{n\choose c_{1},\ldots,c_{J}}\sum_{\bm{i}\in\mathcal{S}(\mathbf{c},j,1)}\frac{\frac{\left(\frac{\gamma}{\alpha}\right)_{(|\bm{i}|)}}{J^{|\bm{i}|}}}{(\gamma)_{(n+1)}}\prod_{k=1}^{J}\mathscr{C}(c_{k}+\delta_{k,j},i_{k};\alpha)}{{n\choose c_{1},\ldots,c_{J}}\sum_{\bm{i}\in\mathcal{S}(\mathbf{c},j,0)}\frac{\frac{\left(\frac{\gamma}{\alpha}\right)_{(|\bm{i}|)}}{J^{|\bm{i}|}}}{(\gamma)_{(n)}}\prod_{k=1}^{J}\mathscr{C}(c_{k},i_{k};\alpha)}\\
&\notag\quad=\frac{1}{\gamma+n}\frac{\sum_{\bm{i}\in\mathcal{S}(\mathbf{c},j,1)}\frac{\left(\frac{\gamma}{\alpha}\right)_{(|\bm{i}|)}}{J^{|\bm{i}|}}\prod_{k=1}^{J}\mathscr{C}(c_{k}+\delta_{k,j},i_{k};\alpha)}{\sum_{\bm{i}\in\mathcal{S}(\mathbf{c},j,0)}\frac{\left(\frac{\gamma}{\alpha}\right)_{(|\bm{i}|)}}{J^{|\bm{i}|}}\prod_{k=1}^{J}\mathscr{C}(c_{k},i_{k};\alpha)}.
\end{align}
Finally, we combine \eqref{card_main_6} with \eqref{post_pyp} and \eqref{ratio_card}. In particular, for $l\in[n]$,
\begin{align*}
&\text{Pr}[f_{X_{n+1}}=l\,|\,\mathbf{C}_{J}=\mathbf{c}]\\
&\quad=\sum_{j=1}^{J}\frac{\gamma}{J} \binom{c_j}{l} (1 - \alpha)_{(l)} \frac{\sum_{\bm i \in \mathcal S(\bm c, j, -l)} \frac{\Gamma\left(\frac{\gamma + \alpha}{ \alpha} + |\bm i|\right)}{J^{|\bm i|}} \prod_{k=1}^J \mathscr{C}(c_k - l \delta_{k,j}, i_k; \alpha)}{\sum_{\bm i \in \mathcal S(\bm c, j, 1)}\frac{ \Gamma\left(\frac{\gamma}{\alpha} + |\bm i|\right)}{J^{|\bm i|}}  \prod_{k=1}^J \mathscr{C}(c_k + \delta_{k,j}, i_k; \alpha)}\\
&\quad\quad\times\frac{1}{\gamma+n}\frac{\sum_{\bm{i}\in\mathcal{S}(\mathbf{c},j,1)}\frac{\left(\frac{\gamma}{\alpha}\right)_{(|\bm{i}|)}}{J^{|\bm{i}|}}\prod_{k=1}^{J}\mathscr{C}(c_{k}+\delta_{k,j},i_{k};\alpha)}{\sum_{\bm{i}\in\mathcal{S}(\mathbf{c},j,0)}\frac{\left(\frac{\gamma}{\alpha}\right)_{(|\bm{i}|)}}{J^{|\bm{i}|}}\prod_{k=1}^{J}\mathscr{C}(c_{k},i_{k};\alpha)}\\
&\quad=\frac{\frac{\gamma}{J}}{\gamma+n}(1-\alpha)_{(l)}\sum_{j=1}^{J}{c_{j}\choose l}\frac{\sum_{\bm i \in \mathcal S(\bm c, j, -l)} \frac{\Gamma\left(\frac{\gamma + \alpha}{ \alpha} + |\bm i|\right)}{J^{|\bm i|}} \prod_{k=1}^J \mathscr{C}(c_k - l \delta_{k,j}, i_k; \alpha)}{\sum_{\bm{i}\in\mathcal{S}(\mathbf{c},j,0)}\frac{\Gamma\left(\frac{\gamma}{\alpha}+|\bm{i}|\right)}{J^{|\bm{i}|}}\prod_{k=1}^{J}\mathscr{C}(c_{k},i_{k};\alpha)}.
\end{align*}

\subsection{Proof of Equation \eqref{post_card_pd}}\label{card_pyp_pd}
The proof follows the same arguments developed in Appendix~\ref{app:post_pyp_new}, exploiting the behaviour of generalized factorial coefficients as $\alpha\rightarrow0$. 
From  \eqref{post_card}, we write that
\begin{align*}
&\lim_{\alpha\rightarrow0}\text{Pr}[f_{X_{n+1}}=l\,|\,\mathbf{C}_{J}=\mathbf{c}]\\
&\quad=\lim_{\alpha\rightarrow0}\frac{\frac{\gamma}{J}}{\gamma+n}(1 - \alpha)_{(l)}\sum_{j=1}^{J} \binom{c_j}{l}  \frac{\sum_{\bm i \in \mathcal S(\bm c, j, -l)} \frac{\Gamma\left(\frac{\gamma + \alpha}{ \alpha} + |\bm i|\right)}{J^{|\bm i|}} \prod_{k=1}^J \mathscr{C}(c_k - l \delta_{k,j}, i_k; \alpha)}{\sum_{\bm i \in \mathcal S(\bm c, j, 0)}\frac{ \Gamma\left(\frac{\gamma}{\alpha} + |\bm i|\right)}{J^{|\bm i|}}  \prod_{k=1}^J \mathscr{C}(c_k, i_k; \alpha)}\\
&\quad=\lim_{\alpha\rightarrow0}\frac{\frac{\gamma}{J}}{\gamma+n}(1 - \alpha)_{(l)} \sum_{j=1}^{J} \binom{c_j}{l} \frac{\sum_{\bm i \in \mathcal S(\bm c, j, -l)} \frac{\prod_{t=0}^{|\bm{i}-1|}\left(\gamma+\alpha+\alpha t\right)}{J^{|\bm i|}} \prod_{k=1}^J \frac{\mathscr{C}(c_k - l \delta_{k,j}, i_k; \alpha)}{\alpha^{i_{k}}}}{\sum_{\bm i \in \mathcal S(\bm c, j, 0)}\frac{ \prod_{t=0}^{|\bm{i}-1|}\left(\gamma+\alpha t\right)}{J^{|\bm i|}}  \prod_{k=1}^J \frac{\mathscr{C}(c_k , i_k; \alpha)}{\alpha^{i_{k}}}}\\
&\quad=\frac{\frac{\gamma}{J}}{\gamma+n} \sum_{j=1}^{J}\binom{c_j}{l} l! \frac{\sum_{\bm i \in \mathcal S(\bm c, j, -l)} \left(\frac{\gamma}{J}\right)^{|\bm{i}|} \prod_{k=1}^J |s(c_k - l \delta_{k,j},i_{k})|}{\sum_{\bm i \in \mathcal S(\bm c, j, 0)}\left(\frac{ \gamma}{J}\right)^{\bm{i}}  \prod_{k=1}^J |s(c_k ,i_{k})|}.
\end{align*}
Exploiting the definition of signless Stirling numbers of the first type as the coefficients of the series expansion of a rising factorial \citep[Equation 2.4]{Cha(05)}, we obtain
\begin{align*}
\lim_{\alpha\rightarrow0}\text{Pr}[f_{X_{n+1}}=l\,|\,\mathbf{C}_{J}=\mathbf{c}]
& =\frac{\frac{\gamma}{J}}{\gamma+n}\sum_{j=1}^{J} \binom{c_j}{l} l! \frac{\prod_{k=1}^{J}\left(\frac{\gamma}{J}\right)_{(c_k - l \delta_{k,j})}}{\prod_{k=1}^{J}\left(\frac{\gamma}{J}\right)_{(c_k)}}\\
& =\frac{\frac{\gamma}{J}}{\gamma+n}\sum_{l=1}^{J}\frac{(c_{j}-l+1)_{(l)}}{\left(\frac{\gamma}{J}+c_{j}-l\right)_{(l)}},
\end{align*}
which coincides with the posterior in \eqref{post_card_pd} by setting $\theta=\gamma$. 

\subsection{Proof of \Cref{prop:rec_trait_all}}\label{app:prof_trait_all}
Without loss of generality, assume $A_{i,k} \in \mathbb N_0 := \{0, 1, \ldots\}$.
From the Poisson process representation of CRMs and the marking theorem \citep{Kin(93)}, $\tilde N := \{(\omega_k, J_k, (A_{i, k})_{i=1}^{n+1})\}_{k \geq 1}$ is a Poisson process on $\mathbb S \times \R_+ \times \N_0^{n+1}$.
Consider now thinned processes $\tilde N_j$, $j=1, \ldots, J$ obtained from $\tilde N$ by taking only those points for which $\omega_k \in D_j$. By the coloring theorem \citep{Kin(93)}, the $\tilde N_j$'s are independent Poisson processes.

Now observe that the random variables $f_{Y_{n+1,r}}$, $X_{n+1}(Y_{n+1,r})$ depend only on $\tilde N_{h(Y_{n+1,r})}$. Similarly, each $(C_j, B{j})$ depends only on $\tilde N_j$ and the independence is preserved also when marginalizing the $\tilde N_j$'s.
Hence, $(f_{Y_{n+1,r}}$, $X_{n+1}(Y_{r}), C_{h(Y_{n+1,r})}, B_{h(Y_{r})})$ are
independent of all other $C_k$'s ($k \neq h(Y_{r})$) and all of $B_{k}$ ($k \neq h(Y_{r})$), which yields the proof.

\subsection{Proof of Theorem~\ref{teo:trait_general}}\label{app:proof_trait_general}

We evaluate \eqref{eq:trait_rec_single} by writing
\begin{multline*}
    \prob\left[f_{Y_{n+1,r}} = l, X_{n+1}(Y_{n+1,r}) = a \mid h(Y_{n+1,r}) = j,  C_{j} = c, B_{j} = b \right]  \\
    = \frac{\prob\left[f_{Y_{n+1,r}} = l, X_{n+1}(Y_{n+1,r}) = a,  C_{j} = c, B_{j} = b \mid h(Y_{n+1,r}) = j \right]}{\prob\left[C_{j} = c, B_{j} = b \mid h(Y_{n+1,r}) = j \right]}
\end{multline*}
and computing the numerator and denominator separately.
In the following, we will denote by $D_\omega$ the preimage of $h(\omega)$, for $\omega \in \mathbb S$.

\paragraph{Denominator.}

From the Poisson process representation of CRMs and the marking theorem, $N := \left\{(\omega_k, S_k, \{A_{i, k}\}_{i=1}^{n+1} \right\}_{k \geq 1}$ is a Poisson process on $\mathbb S \times \R_+ \times \mathbb N^{n+1}$ with intensity
\[
    \theta G_0(\dd w) \left[\prod_{i=1}^{n+1} G_A(\dd a_i \mid s)\right] \rho(s) \dd s.
\]
Then, by the colouring theorem \citep[see, e.g., Chapter 5 in][]{Kin(93)}, we have that selecting from $N$ only those points for which $\omega_k \in D_j$ (i.e., the points whose features' hashes coincide with the hash of $Y_{n+1, r}$), leads to a point process $N^\prime := \left\{(\omega^\prime_k, S^\prime_k, \{A^\prime_{i, k}\}_{i=1}^{n+1} \right\}_{k \geq 1}$ which is Poisson on $D_j \times \R_+ \times \mathbb N^{n+1}$ with intensity
\begin{equation}\label{eq:int_rescaled}
    \frac{\theta}{J} \bar G_j(\dd w) \left[\prod_{i=1}^{n+1} G_A(\dd a_i \mid s)\right] \rho(s) \dd s,
\end{equation}
where $\bar G_j(\dd w)$ is $G_0$ truncated on $D_{\omega^*}$ and re-normalized.
Then,
\begin{align*}
     \prob\left[ C_{j} = c, B_j = b \right] 
    & = \prob\left[ \sum_{k \geq 1} \indicator[\omega_k \in D_{j}] \sum_{i=1}^n A_{i, k} = c, \sum_{k \geq 1} \indicator[\omega_k \in D_{j}]  A_{i, n+1} = b \right] \\
    & = \prob\left[ \sum_{k \geq 1} \sum_{i=1}^m A^\prime_{i, k} = c, \sum_{k \geq 1}  A^\prime_{n+1, k} = b\right].
\end{align*}
Observe that $\omega_k$ is not involved in the last probability, so we can marginalize with respect to it and consider the point process $\{J^\prime_k, (A^\prime_{i,k})_{i=1}^{n+1}\}_{k \geq 1}$ on $\R_+ \times \mathbb N_0^{n+1}$ with intensity $\theta J^{-1} \left[\prod_{i=1}^{n+1} G_A(\dd a_i \mid s)\right] \rho(s) \dd s$.

\paragraph{Numerator.}

We start by considering $Y_{n+1, r} = \omega^*$ fixed.
Let $B_{\omega^*}$ be a ball of radius $\varepsilon$ centered in $\omega^*$.
We compute
\[
  \prob\left[\sum_{i=1}^n X_i(B_{\omega^*}) = l, X_{n+1}(B_{\omega^*}) = a, C_j = c, B_j = b \mid h(\omega^*)=j \right],
\]
which converges to the numerator in \eqref{eq:trait_rec_single} as $\varepsilon \rightarrow 0$.
By the definition of $B_j$ and $C_j$ we have
\begin{align}
  &\prob\left[\sum_{i=1}^n X_i(B_{\omega^*}) = l, X_{n+1}(B_{\omega^*}) = a, C_j = c, B_j = b \right] \label{eq:trait_num_temp} \\
  & \qquad = \prob\Big[\sum_{i=1}^n X_i(B_{\omega^*}) = l, X_{n+1}(B_{\omega^*}) = a, \nonumber \\
  & \hspace{1cm}
  \sum_{i=1}^n X_i(D_{j} \setminus B_{\omega^*}) = (c-l), X_{n+1}(D_{j} \setminus B_{\omega^*}) = (b-a)\Big].\nonumber
\end{align}
Further, $(X_i(B_{\omega^*}))_{i \geq 1}$ is a collection of random variables independent of $(X_i(D_{j} \setminus B_{\omega^*}) )_{i \geq 1}$. Therefore, we can consider the first two events and the last two events separately.

For the last two, arguing as in the denominator case, we have that as $\varepsilon \rightarrow 0$
\begin{multline*}
\prob\Big[
  \sum_{i=1}^n X_i(D_{j} \setminus B_{\omega^*}) = (c-l), X_{n+1}(D_{j} \setminus B_{\omega^*}) = (b-a)\Big] \\
\longrightarrow \prob\left[ \sum_{k \geq 1} \sum_{i=1}^n A^\prime_{i, k} = c - l, \sum_{k \geq 1}  A^\prime_{n+1, k} = b - a\right],
\end{multline*}
where the $A^\prime_{ik}$'s come from the Poisson process $\{J^\prime_k, (A^\prime_{i,k})_{i=1}^{n+1}\}_{k \geq 1}$ on $\R_+ \times \mathbb N_0^{n+1}$ with intensity $\theta J^{-1} \left[\prod_{i=1}^{n+1} G_A(\dd a_i \mid s)\right] \rho(s) \dd s$.

For the first two events instead, an application of Campbell's theorem yields
\begin{align*}
  & \prob\Big[\sum_{i=1}^n X_i(B_{\omega^*}) = l, X_{n+1}(B_{\omega^*}) = a \Big] \\
  & \qquad = \E\left[\int I\left[\left\{\sum_{i=1}^n a_i\right\} = l\right] I[a_{n+1} = a] I[\omega \in B_{\omega^*}] N(\dd s \, \dd \bm a \, \dd \omega) \right] \\
  & \qquad = \theta G_0(B_{\omega^*}) \int \prob \left[ \sum_{i=1}^n \tilde A_{i} = l, \tilde A_{n+1} = a \mid s \right] \rho(s) \dd s.
\end{align*}
The proof concludes by letting $\varepsilon \rightarrow 0$ and noting that $\int_{D_j} G_0(\dd \omega^*) = J^{-1}$.

\subsection{Proof of \Cref{prop:poisson_traits}}\label{app:proof_poisson}

The closeness of the Poisson distribution under convolution entails that $\sum_{i=1}^n \tilde A_{i} \mid s \sim \mbox{Poi}(n \lambda s)$.
Similarly, $\sum_{k \geq 1}\sum_{i=1}^n A^\prime_{i, k} \sim \mbox{Poi}(n \lambda T^\prime)$ where $T^\prime = \sum_{k \geq 1} J^\prime_k$.
Then,
\begin{align*}
  \int_{\R_+} \prob \left[ \sum_{i=1}^n \tilde A_{i} = l, \tilde A_{n+1} = a \mid s \right] \rho(s) \dd s 
  & = \int_{\R_+} \frac{1}{l!} \frac{1}{a!} (nr s)^l e^{-nrs} (\lambda s)^a e^{- \lambda s} \rho(s) \dd s \\
  &  = \frac{n^{l} \lambda^{l + a}}{l! a!} \kappa(l+a, (n+1) \lambda).
\end{align*}
Moreover
\begin{align*}
  \prob\left[ \sum_{k \geq 1} \sum_{i=1}^n A^\prime_{i, k} = c, \sum_{k \geq 1}  A^\prime_{n+1, k} = b \right]
  & = \frac{n^c \lambda^{c+b}}{c! b!}  \E_{T^\prime} \left[e^{-(n+1) \lambda T^\prime} (T^\prime)^{c+b} \right]\\
  & =\frac{n^c \lambda^{c+b}}{c! b!}  (-1)^{c+b} \frac{\dd^{c+b}}{\dd ((n+1)\lambda)^{c+b}} \E_{\mu^\prime}[e^{-(n+1)\lambda T^\prime}] \\
  &  = \frac{n^c \lambda^{c+b}}{c! b!} (-1)^{c+b} \frac{\dd^{c+b}}{\dd z^{c+b}} \exp\left(- \theta/J \psi(z) \right)\big|_{(n+1)\lambda}.
\end{align*}
In a similar fashion,
\begin{equation}\label{eq:proof_poi_num}
  \begin{aligned}
  & \prob\left[ \sum_{k \geq 1} \sum_{i=1}^n A^\prime_{i, k} = c - l, \sum_{k \geq 1}  A^\prime_{n+1, k} = b - a\right] \\
  & \qquad \qquad \frac{n^{c-l} \lambda^{c - l+b - a}}{(c-l)! (b-a)!} (-1)^{c-l+b-a} \frac{\dd^{c-l+b-a}}{\dd z^{c-l+b-a}} \exp\left(- \theta/J \psi(z) \right)\big|_{(n+1)\lambda}.
  \end{aligned}
\end{equation}
Combining these expressions together yields the proof.

\subsection{Proof of \Cref{prop:bern_traits}}

Observe that $(J^\prime_k)_{k \geq 1}$ in \Cref{teo:trait_general} is a Poisson process on $\R_+$ with intensity $\theta/J \rho(s) \dd s$.
Consider now the numerator in \eqref{eq:post_trait}.
Let $S_{n} := \sum_{k \geq 1} \sum_{i=1}^n A^\prime_{i,k}$ and $Z:= \sum_{k \geq 1} A^{\prime}_{n+1, k}$.
Conditional on $(J^\prime_k)_{k \geq 1}$, $S_{n,k}$ is Poisson-binomial with parameters $J^\prime_1, \ldots, J^\prime_1, J^\prime_2, \ldots, J^\prime_2, \ldots$, where each $J^\prime_k$ appears exactly $n$ times. Similarly, $Z \mid (J^\prime_k)_{k \geq 1}$ is Poisson-binomial with parameters $J_1, J_2, \ldots$.
Let $\tilde S_n \mid (J^\prime_k)_{k \geq 1} \sim \mbox{Poi}(n T^\prime)$ and $\tilde Z \mid (J^\prime_k)_{k \geq 1} \sim \mbox{Poi}(T^\prime)$ where $T^\prime = \sum_{k} J^\prime_k$.
Then, by \cite{LeCam(60)}, conditionally to $(J^\prime_k)_{k \geq 1}$, $S_n \approx \tilde S_n$, $Z \approx \tilde Z$. Hence
\[
  \prob[S_n = c-l, Z = b- 1] \approx \prob[\tilde S_n = c-l, \tilde Z = b- 1].
\]
Then, from \eqref{eq:proof_poi_num} we get
\begin{multline*}
  \prob[\tilde S_n = c-l, \tilde Z = b - 1] = \\ \frac{n^{c-l} \lambda^{c - l+ b - 1}}{(c-l)! (b-1)!} (-1)^{c-l+b-1} \frac{\dd^{c-l+b-1}}{\dd z^{c-l+b-1}} \exp\left(- \theta/J \psi(z) \right)\big|_{(n+1) \lambda }.
\end{multline*}

Moreover,
\begin{align*}
   \prob \left[ \sum_{i=1}^n \tilde A_{i} = l, \tilde A_{n+1} = 1 \mid s \right]
  &  = \binom{n}{l} \prob[\tilde A_{1} = 1, \ldots, \tilde A_{l} = 1, \tilde A_{n+1} = 1, \tilde A_{l+1} = 0, \ldots, \tilde A_{n}=0] \\
  & = \binom{n}{l} s^{n+1} (1 - s)^{m - l} = \binom{n}{l} \left(s^{l+1} - s^{n+1} \right).
\end{align*}
Integrating this with respect to $\rho(s) \dd s$ and ignoring multiplicative terms yields \eqref{eq:post_trait}.

To prove the error bound, for ease of notation, let $f_Y \equiv f_{Y_{n+1, r}}$, $X \equiv X_{n+1}(Y_{n+1, r})$ and let $M^\prime= (J_k)_{k\geq 1}$.
Further, note that $(f_Y, X)$ is independent of $(S_n, Z)$ and $(\tilde S_n, \tilde Z)$.
Then,
\begin{align*}
  & TV\left(\left[f_Y, X, S_n, Z\right], \left[f_Y, X, \tilde S_n, \tilde Z\right]\right) \\
  & \quad = \sum_{l, x, s, z} \Big| \prob\left[f_Y = l, X = x, S_n = s, Z = z\right] - \prob\left[f_Y = l, X = x, \tilde S_n = s, \tilde Z = z\right] \Big| \\
  &\quad = \sum_{l, x, s, z} \prob\left[f_Y = l, X= x \right] \Big| \prob\left[S_n = s, Z = z\right] -  \prob\left[\tilde S_n = s, \tilde  Z = z\right]\Big| \\
  & \quad \leq  \sum_{s, z}  \Big| \E \prob\left[S_n = s, Z = z \mid M^\prime \right] -  \E \prob\left[\tilde S_n = s, \tilde  Z = z \mid M^\prime \right]\Big| \\
  & \quad \leq \E \sum_{s, z} \Big| \prob\left[S_n = s, Z = z \mid M^\prime \right] -  \prob\left[\tilde S_n = s, \tilde  Z = z \mid M^\prime \right]\Big|,
\end{align*}
where the last inequality follows from Jensen's inequality and an application of the Fubini theorem.
Then by the conditional independence between $S_n$ and $Z$, and $\tilde S_n$ and $\tilde Z$ we get
\begin{align*}
  &  \E \sum_{s, z} \Big| \prob\left[S_n = s, Z = z \mid M^\prime \right] -  \prob\left[\tilde S_n = s, \tilde  Z = z \mid M^\prime \right]\Big| \\
& \qquad =  \E\left[TV(S_n \mid M^\prime, \tilde S_n \mid M^\prime) + TV(Z \mid M^\prime, \tilde Z \mid M^\prime) \right].
\end{align*}
To upper bound the error, we start from Eq. (5.5) in \cite{Ste(94)} so that
\begin{align*}
  TV(S_n \mid M^\prime, \tilde S_n \mid M^\prime) &\leq \frac{(1 - e^{-nT^\prime})}{nT^\prime} n \sum_{k \geq 1} J_k^{\prime 2} \\
  TV(Z \mid M^\prime, \tilde Z \mid M^\prime) &\leq \frac{(1 - e^{-T^\prime})}{T^\prime}  \sum_{k \geq 1} J_k^{\prime 2},
\end{align*}
so that
\begin{align*}
 & \E\left[TV(S_n \mid M^\prime, \tilde S_n \mid M^\prime) + TV(Z \mid M^\prime, \tilde Z \mid M^\prime) \right] \\
& \qquad \leq \E \left[\frac{2}{T^\prime}  \sum_{k \geq 1} J_k^{\prime 2} \right] = 2 \int_{\R_+} \E \left[e^{-u T^\prime} \sum_{k \geq 1} J_k^{\prime 2}\right] \dd u \\
& \qquad = 2 \int_{\R_+} \E \left[ \int_{\R_+} e^{-u \int_{\R_+} z M^\prime(\dd z)} s^2 M^\prime(\dd z)\right] \dd u,
\end{align*}
where the last equality follows by identifying $M^\prime$ with the random counting measure $\sum_{k \geq 1} \delta_{J^\prime_k}$.
Then, Mecke's equation \citep[Theorem 4.1]{Las(18)} leads to:
\begin{align*}
  2 \int_{\R_+} \E \left[ \int_{\R_+} e^{-u \int_{\R_+} z M^\prime(\dd z)} s^2 M^\prime(\dd z)\right] \dd u
  &  = \frac{2\theta}{J} \int_{\R_+} \int_{\R_+} \E[e^{- u T^\prime}] e^{-us} s^2 \rho(s) \dd s.
\end{align*}
The proof follows by noticing that, by the L\'evy-Kintchine representation, $\E[e^{- u T^\prime}] = e^{-\psi(u)}$ and from the definition of $\kappa(u, n)$.


\subsection{Monte Carlo Estimation of \eqref{post_pyp}}\label{app:computations_pyp}

As in \cite{Dol(23)}, we note that \eqref{post_pyp} can be equivalently expressed as:
\begin{multline}\label{eq:rec_nrm_py_mc}
    \prob\left[ f_{X_{n+1}} = l \mid \bm C_{J} = \bm c, h(X_{n+1}) = j\right] = \\ \frac{\gamma}{J} \binom{c_j}{l} (1 - \alpha)_{(l)} \frac{(\gamma)_{(c_j -l)}}{(\gamma)_{(c_j + 1)}} \frac{\E\left[\frac{ \left(\frac{\gamma + \alpha}{\gamma} \right)_{K^{\mathcal S(\bm c, j, -l) }_{\boldsymbol \cdot}}}{J^{K^{\mathcal S(\bm c, j, -l) }_{\boldsymbol \cdot}} \prod_{k=1}^J \left(\frac{\gamma}{\alpha}\right)_{K_{c_k - l \delta_{k,j}}}} \right]}{\E\left[\frac{ \left(\frac{\gamma }{\gamma} \right)_{K^{\mathcal S(\bm c, j, 1) }_{\boldsymbol \cdot}}}{J^{K^{\mathcal S(\bm c, j, 1) }_{\boldsymbol \cdot}} \prod_{k=1}^J \left(\frac{\gamma}{\alpha}\right)_{K_{c_k + \delta_{k,j}}}} \right]},
\end{multline}
where $K_m$ is number of distinct values in a sample of size $m$ from a PYP and $K^{\mathcal S(\bm c, j, -l) }_{\boldsymbol \cdot} := \sum_{1\leq k\leq J} K_{c_k - l \delta_{k, j}}$. To prove \eqref{eq:rec_nrm_py_mc}, recall that the number of distinct elements $K_m$ in a sample of size $m$ from a PYP has distribution
\[
    \prob\left[K_m = k\right] = \frac{\left(\frac{\gamma}{\alpha} \right)_{(k)}}{(\gamma)_{(m)}} \mathscr{C}(m, k; \alpha).
\]

Consider now \eqref{post_pyp} and note that $\Gamma\left(a + b\right) = (a)_{(b)} \Gamma(a)$.
Then
\begin{align*}
    &\prob\left[ f_{X_{n+1}} = l \mid \bm C_{J} = \bm c, h(X_{n+1}) = j\right] = \\
    & \qquad = \frac{\gamma}{J} \binom{c_j}{l} (1 - \alpha)_{(l)} \frac{(\gamma)_{(c_j -l)}}{(\gamma)_{(c_j + 1)}}  \frac{\sum_{\bm i \in \mathcal S(\bm c, j, -l)} \frac{\left(\frac{\gamma + \alpha}{\alpha}\right)_{(|\bm i|)}}{J^{|\bm i|} \prod_{h=1}^J \left(\frac{\gamma}{\alpha} \right)_{(i_h)}}  \prod_{k=1}^J \frac{\left(\frac{\gamma}{\alpha}\right)_{(i_k)}}{(\gamma)_{(c_k - l \delta_{k, j})}}\mathscr{C}(c_k - l \delta_{k,j}, i_k; \alpha)
    }{
    \sum_{\bm i \in \mathcal S(\bm c, j, 1)} \frac{\left(\frac{\gamma}{\alpha}\right)_{(|\bm i|)}}{J^{|\bm i|} \prod_{h=1}^J \left(\frac{\gamma}{\alpha} \right)_{(i_h)}}  \prod_{k=1}^J \frac{\left(\frac{\gamma}{\alpha}\right)_{(i_k)}}{(\gamma)_{(c_k + \delta_{k, j})}}\mathscr{C}(c_k - l \delta_{k,j}, i_k; \alpha)
    } \\
     & \qquad = \frac{\gamma}{J} \binom{c_j}{l} (1 - \alpha)_{(l)} \frac{(\gamma)_{(c_j -l)}}{(\gamma)_{(c_j + 1)}}  \frac{\sum_{\bm i \in \mathcal S(\bm c, j, -l)} \frac{\left(\frac{\gamma + \alpha}{\alpha}\right)_{(|\bm i|)}}{J^{|\bm i|} \prod_{h=1}^J \left(\frac{\gamma}{\alpha} \right)_{(i_h)}}  \prod_{k=1}^J \prob\left[K_{c_k - l \delta_{k,j}} = i_k\right]
    }{
    \sum_{\bm i \in \mathcal S(\bm c, j, 1)} \frac{\left(\frac{\gamma}{\alpha}\right)_{(|\bm i|)}}{J^{|\bm i|} \prod_{h=1}^J \left(\frac{\gamma}{\alpha} \right)_{(i_h)}}  \prod_{k=1}^J
    \prob\left[K_{c_k + \delta_{k,j}} = i_k\right]
    }\\
    & \qquad = \frac{\gamma}{J} \binom{c_j}{l} (1 - \alpha)_{(l)} \frac{(\gamma)_{(c_j -l)}}{(\gamma)_{(c_j + 1)}} \frac{\E\left[\frac{ \left(\frac{\gamma + \alpha}{\gamma} \right)_{K^{\mathcal S(\bm c, j, -l) }_{\boldsymbol \cdot}}}{J^{K^{\mathcal S(\bm c, j, -l) }_{\boldsymbol \cdot}} \prod_{k=1}^J \left(\frac{\gamma}{\alpha}\right)_{K_{c_k - l \delta_{k,j}}}} \right]}{\E\left[\frac{ \left(\frac{\gamma }{\gamma} \right)_{K^{\mathcal S(\bm c, j, 1) }_{\boldsymbol \cdot}}}{J^{K^{\mathcal S(\bm c, j, 1) }_{\boldsymbol \cdot}} \prod_{k=1}^J \left(\frac{\gamma}{\alpha}\right)_{K_{c_k + \delta_{k,j}}}} \right]}.
\end{align*}

\section{Details about the Poisson-IBP}\label{app:poisson}

\paragraph{Proof of Equation \eqref{eq:pois_gamma}}
It follows from $\psi(u) = \log(1 + u)$ and $\kappa(u, n) = (n-1)! / (u + 1)^n$.

\paragraph{Proof of \eqref{eq:pois_gg}}

We have
\begin{align*}
    \psi(u) &= \frac{\alpha}{\Gamma(1 - \alpha)} \int_{\R_+}(1 - e ^{-us}) \rho(s) \dd s =  \frac{\alpha}{\Gamma(1 - \alpha)} \int_{\R_+}(1 - e ^{-us}) s^{-1 - \alpha} e^{-\tau s}\dd s \\
    &=  \frac{\alpha}{\Gamma(1 - \alpha)} \int_{\R_+} \int_{\R_+} (1 - e ^{-us}) \frac{t^\alpha e^{-ts}}{\Gamma(\alpha + 1)} e^{-s} \dd t \, \dd s \\
    &=  \frac{\alpha}{\Gamma(1 - \alpha)} \int_{\R_+} \frac{t^\alpha u}{(t+1)(t+u+1)} \dd t = \frac{\alpha \Gamma(\alpha) \Gamma(1 - \alpha)}{\Gamma(1 - \alpha) \Gamma(\alpha + 1)} \left[(\tau + u)^\alpha - \tau^\alpha \right] \\
    &= \left[(\tau + u)^\alpha - \tau^\alpha \right],
\end{align*}
and
\[
    \kappa(l, u) = \alpha \frac{\Gamma(l - \alpha)}{\Gamma(1 - \alpha)} (\tau + u)^{\alpha - l} = \alpha (1 - \alpha)_{(l-1)}(\tau + u)^{\alpha - l}.
\]
Denoting by $(*)$ the summation over positive integers $(k_1, \ldots, k_i)$ such that $\sum_{j=1}^i k_j = n$,
\begin{align*}
    \frac{\dd^n}{\dd u^n} e^{- \theta \psi(u) / J} &= e^{\theta/J \tau^\alpha} \frac{\dd^n}{\dd u^n} e^{- \theta/J (\tau + u)^\alpha} \\
    &= e^{-\theta/J \psi(u)} \sum_{i=1}^n \left(\frac{\theta}{J}\right)^{i} \sum_{(*)} \frac{1}{i!} \binom{n}{k_1, \ldots, k_i} \prod_{j=1}^i \alpha (\tau + u)^{\alpha - k_j} (1 - \alpha)_{k_j - 1} \\
    &= e^{-\theta/J \psi(u)} \sum_{i=1}^n \left(\frac{\theta}{J}\right)^{i} \frac{\calC(n, i; \alpha)}{(\tau + u)^{n - \alpha i}}.
\end{align*}
Plugging these in the expression of \Cref{prop:poisson_traits} leads to \eqref{eq:pois_gg}.

\FloatBarrier

\section{Additional Numerical Results}\label{app:numerical}

\subsection{Cardinality recovery} \label{app:numerical-cardinality}

\begin{figure}[!htb]
\centering
\includegraphics[width=0.8\linewidth]{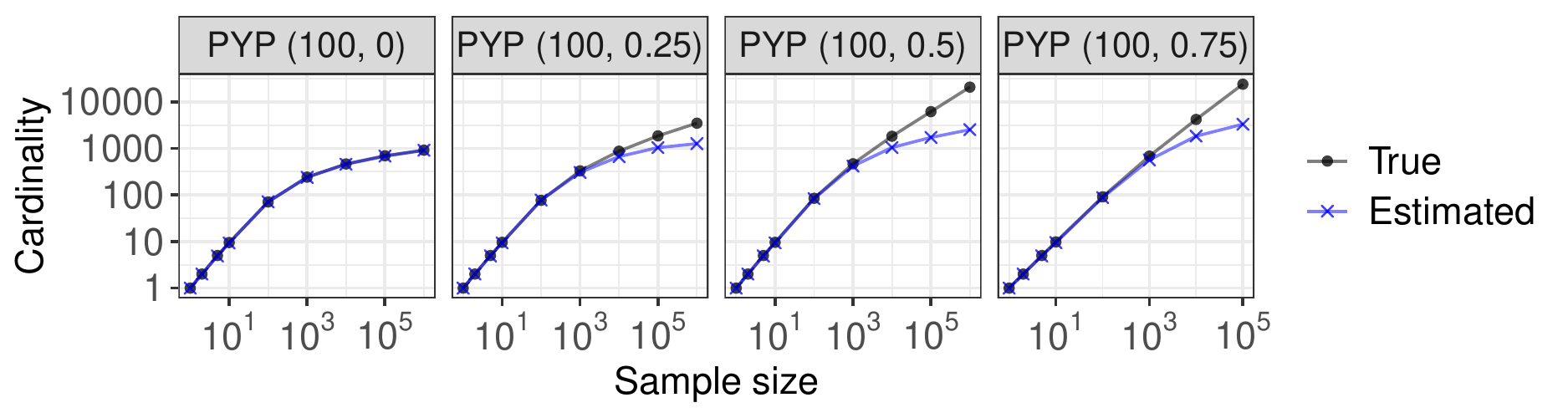}
\caption{\small{True and estimated cardinality in synthetic data from PYP prior models, as a function of the sample size. The estimates assume a mis-specified DP prior fitted via maximum marginal likelihood. Other details are as in Figure~\ref{fig:exp-dp-dp-distinct}.}}
\label{fig:exp-dp-pyp-distinct}
\end{figure}

\begin{figure}[!htb]
\centering
\includegraphics[width=0.8\linewidth]{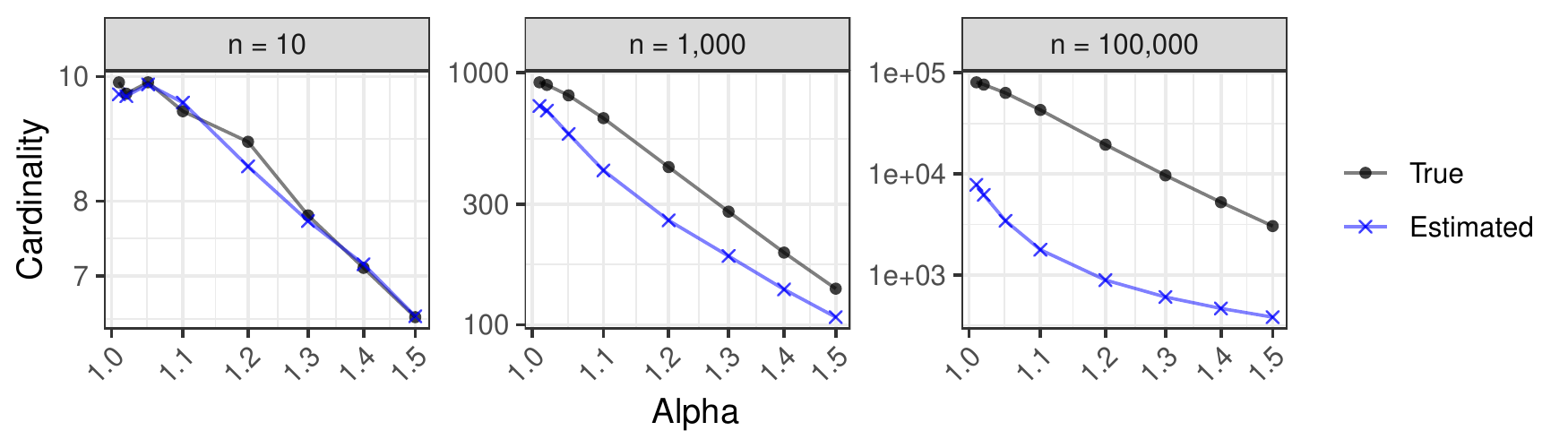}
\caption{\small{True and estimated cardinality for data sampled from a Zipf distribution, as a function of the tail parameter $\alpha$. The estimates assume a mis-specified DP prior fitted via maximum marginal likelihood. Other details are as in Figure~\ref*{fig:exp-dp-dp-distinct}.}}
\label{fig:exp-dp-zipf-distinct}
\end{figure}

\FloatBarrier


\end{document}